\documentclass[11pt]{amsart}

\usepackage{amssymb}
 \usepackage[all]{xy}

 \newtheorem{theorem}{Theorem}[section]
 \newtheorem{proposition}[theorem]{Proposition}
 \newtheorem{lemma}[theorem]{Lemma}
 \newtheorem{corollary}[theorem]{Corollary}

 \newtheorem{lemma/definition}[theorem]{Lemma/Definition}
 \newtheorem{proposition/definition}[theorem]{Proposition/Definition}
 \newtheorem{lem/not}[theorem]{Lemma/Notations}

 \newtheorem{definition}[theorem]{Definition}
 \newtheorem{example}[theorem]{Example}

 \newtheorem{remark}[theorem]{Remark}

\numberwithin{equation}{section}

 \newcommand{\Rr}{\mathbb R}

 \newcommand{\rmap}{\longrightarrow}

 \newcommand{\F}{\ensuremath{\mathcal{F}}}

\newcommand{\hhstar}{\mathfrak{h}^*}

\newcommand{\Dir}{\text{\rm Dir}\,}        % set of Dirac structures 
\newcommand{\For}{\mathfrak{F}}            % forward map
\newcommand{\Ran}{\mathcal{R}}             % range of Dirac structure

\newcommand{\SP} [1]     {{\left\langle {{#1}} \right\rangle}}

 \newcommand{\s}{\mathbf{s}}             % source map
 \renewcommand{\t}{\mathbf{t}}           % target map
 \newcommand{\Lie}{\mathcal{L}}          % Lie derivative
 \newcommand{\Ker}{\text{\rm Ker}\,}     % Kernel
 \renewcommand{\Im}{\text{\rm Im}\,}     % Image
 \newcommand{\Ad}{\text{\rm Ad}\,}       % Adjoint
 \newcommand{\proj}{p}                    % general projection 
 
\newcommand{\lcart} {\lambda}           %Maurer-Cartan forms
\newcommand{\rcart} {\overline{\lambda}} 

\newcommand{\R}{\ensuremath{\mathbb R}}

\begin{document}

 \title[Integration of twisted Dirac brackets]{Integration of twisted Dirac brackets}

% author one information
\author{Henrique Bursztyn}
\address{Department of Mathematics, University of Toronto, Toronto, Ontario M5S3G3, Canada}
\author{Marius Crainic}
\address{Department of Mathematics, University of California, Berkeley, CA 94720, USA, and Utrecht University, P.O. Box 80.010, 3508 TA, Utrecht, The Netherlands}
\author{Alan Weinstein}
\address{Department of Mathematics, University of California, Berkeley, CA 94720, USA}
\author{Chenchang Zhu}
\address{Department of Mathematics, University of California, Berkeley, CA 94720, USA}

% Use this \subjclass if you are using amsproc version 2.0 (December 1999).
\date{}

%%%

\begin{abstract}
Given a Lie groupoid $G$ over a manifold
$M$, we show that multiplicative 2-forms on $G$ relatively closed with
respect to a closed 3-form $\phi$ on $M$ correspond to maps from the Lie algebroid
of $G$ into $T^*M$ satisfying an algebraic condition and a differential
condition with respect to the $\phi$-twisted Courant bracket. This
correspondence describes, as a special case, the global objects associated
to $\phi$-twisted Dirac structures. As applications, we relate our results
to equivariant cohomology and foliation theory, and we give a new description
of quasi-hamiltonian spaces and group-valued momentum maps.
\end{abstract}
%%%

\maketitle
%%%

% \tableofcontent
\tableofcontents

\section{Introduction}
The correspondence between Poisson structures and symplectic groupoids \cite{CaFe,CDW,CrFe2},
analogous to the one of Lie algebras and Lie groups, plays an important role in Poisson
geometry; it offers, in particular, a unifying framework for the study of hamiltonian
and Poisson actions (see e.g. \cite{WeGeom}). In this paper, we extend this correspondence to the
context of Dirac structures twisted by a closed 3-form.

Dirac structures were introduced by Courant in \cite{Cour,Co-We},
motivated by the study of constrained mechanical systems.
Examples of Dirac structures include Poisson structures, presymplectic forms and regular foliations.
Connections between Poisson geometry and topological sigma models \cite{Strobl,Park} have led to
the notion of Poisson structures twisted by a closed 3-form, which were described by Weinstein and Severa
in \cite{WeSe} as special cases of twisted Dirac structures. The diagram below describes these various
generalizations of Poisson structures and the corresponding global objects.

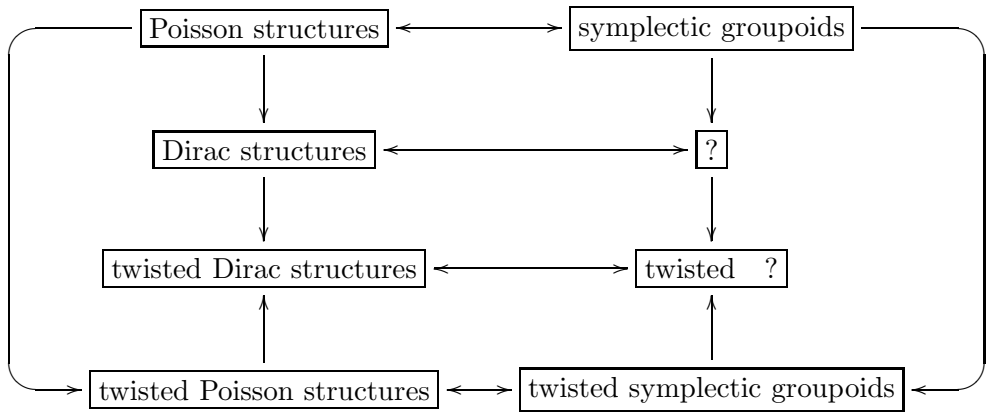
\[ \xymatrix{
& \fbox{\parbox{.3\linewidth}\center{Poisson structures}} \ar@{<->}[r]\ar[d]\ar@{->}`l[ldd] `[ddd] [ddd] & \fbox{\parbox{.3\linewidth}\center{symplectic groupoids}} \ar[d]\ar@{->}`r[rdd] `[ddd] [ddd] & \\
& \fbox{\parbox{.3\linewidth}\center{Dirac structures}} \ar@{<->}[r] \ar[d] & \fbox{\parbox{.3\linewidth}\center{?}}\ar[d] & \\
& \fbox{\parbox{.3\linewidth}\center{twisted Dirac structures}} \ar@{<->}[r] & \fbox{\parbox{.3\linewidth}\center{twisted \ \ ?}} & \\
& \fbox{\parbox{.3\linewidth}\center{twisted Poisson structures}} \ar@{<->}[r] \ar[u] & \fbox{\parbox{.3\linewidth}\center{twisted symplectic groupoids}}\ar[u] & 
} \]

In this paper, we fill the gaps in the diagram above by introducing
 the notion of  {\bf twisted presymplectic
groupoid} relative to a closed 3-form $\phi$ on a manifold $M$ and
establishing a bijective correspondence (up to natural isomorphisms)
between source-simply-connected $\phi$-twisted presymplectic
groupoids over $M$ and $\phi$-twisted Dirac structures.
Our results complete,
and were inspired by, earlier work of Bursztyn and Radko \cite{BuRa}
on gauge transformations of Poisson structures, and of Cattaneo and Xu
\cite{CaXu}, who used the method of sigma models \cite{CaFe}
to construct global objects attached to twisted Poisson structures.
We remark that the  degeneracies of presymplectic forms
force our techniques to be different in nature
from those used to establish the correspondence of Poisson structures and
symplectic groupoids (as in \cite{CDW, Wein3, CaXu}).
In our search for the right notion of nondegeneracy that
2-forms on presymplectic groupoids must satisfy (Definition \ref{def:twistpresymp}),
the examples provided by \cite{BuRa, CrFe2}
were essential. (In fact, our presymplectic groupoids are closely related to those introduced in [14]).

Twisted presymplectic groupoids turn out to be best approached through
a general study of multiplicative 2-forms on groupoids, which 
turn out to be extremely rigid.  (On groups, they are all zero.)
Given a closed 3-form $\phi$ on the base $M$ of a groupoid $G$, we
call a 2-form $\omega$ on $G$ {\bf relatively} $\phi$-{\bf closed}
if $d\omega=s^*\phi-t^*\phi$, where $s$ and $t$ are the source and
target maps of $G$.  If $G$ is source-simply-connected, we
identify the infinitesimal counterparts of relatively
$\phi$-closed multiplicative 2-forms on $G$ as those maps to $T^*M$ from the Lie
algebroid $A$ of $G$ which satisfy an algebraic condition and a
differential condition related to the $\phi$-twisted Courant bracket \cite{WeSe}.
In fact, the [twisted] Courant bracket itself could be rediscovered
from the properties of [relatively] closed 2-forms on groupoids.
The reconstruction of multiplicative 2-forms out of the infinitesimal
data is based on the constructions of \cite[Sec.~4.2]{CrFe}.

Motivated by the relationship between symplectic realizations
of Poisson manifolds and hamiltonian actions \cite{CDW},
we study {\bf presymplectic realizations} of twisted Dirac structures.
Just as in usual Poisson geometry, these presymplectic realizations carry
natural actions of  presymplectic groupoids.
In fact, it is this property that determines our definition of
presymplectic realizations. An important example of  twisted Dirac structure
is described in \cite[Example~4.2]{WeSe}: any 
nondegenerate invariant inner product on the Lie algebra
$\mathfrak{h}$  of a Lie group $H$ induces a natural 
Dirac structure on $H$, twisted
by the invariant Cartan 3-form; we call such structures 
{\bf Cartan-Dirac structures}. We show that 
presymplectic realizations of Cartan-Dirac structures are  equivalent to the  
quasi-hamiltonian $\mathfrak{h}$-spaces of Alekseev, 
Malkin and Meinrenken \cite{AMM} in such a way that
the realization maps  are the associated group-valued momentum maps.
It also follows from our results that the transformation
groupoid $H\ltimes H$ corresponding to the conjugation
action carries a canonical twisted presymplectic structure,
which we obtain explicitly by ``integrating''
the Cartan-Dirac structure. As a result,  we recover the 2-form on the ``double'' $D(H)$
of \cite{AMM}
and the AMM groupoid of \cite{bxz}. (Closely related forms were
introduced  earlier in  \cite{hj,WeFloer}.) A unifying approach to momentum map
theories based on Morita equivalence of presymplectic groupoids has been developed by Xu in
\cite{Xu}; much of our motivation for considering quasi-hamiltonian spaces comes
from his work. Our results indicate that Dirac structures provide a natural framework
for the common description of various notions of momentum maps 
(as e.g. in \cite{AMM,Lu}, see also \cite{WeGeom}).

We illustrate our results on multiplicative 2-forms and Dirac structures in many examples.
In the case of action groupoids,
we obtain an explicit formula for the natural map 
from the cohomology of the Cartan model of an $H$-manifold \cite{AtBo, Vergne}
to (Borel) equivariant cohomology \cite{AtBo, BSS} in degree three;
for the monodromy groupoid of a foliation $\mathcal{F}$,
we show that multiplicative 2-forms are closely related
to the usual cohomology and spectral sequence of $\mathcal{F}$ \cite{KmTd}.

The entire discussion of relatively closed multiplicative 2-forms on
groupoids may be embedded in the more general context of a van Est
theorem for the ``bar-de
Rham'' double
complex of forms on the simplicial space of composable sequences in a
groupoid $G$, whose total complex computes the cohomology of the classifying
space $BG$.  We reserve this more general discussion for a future paper.\\

{\bf Outline of the paper:}
In Section \ref{sec-basic}, we review basic definitions concerning
Dirac structures and groupoids, we introduce various notions of
``adapted'' 2-forms on groupoids, and we state our two main results,
which allow one to go back and forth between [twisted] presymplectic
groupoids and their infinitesimal counterparts, [twisted] Dirac manifolds.
Section \ref{sec-mult} begins the work of proving the theorems with a
study of the global objects, namely multiplicative 2-forms on
groupoids.  In Section \ref{sec:nondeg}, we pass from the groupoids to the
infinitesimal objects, and in Section \ref{Reconstructing
  multiplicative forms}, we go in the opposite direction, completing
the proofs of the main theorems.  Section \ref{examples} begins with
simple examples.  Then, after a discussion of 2-forms on groupoids which
become presymplectic only after one passes to the quotient by a
foliation, we show how twisted-multiplicative forms on action
groupoids are related to equivariant cohomology.  This leads us to
Section \ref{sec-realizations}, where we study the AMM-groupoid and
apply our results to quasi-hamiltonian actions. Finally, Section
\ref{sec-foliations} is devoted to multiplicative 2-forms on foliation
groupoids, with applications to Dirac structures whose presymplectic
leaves all have the same dimension.\\

{\bf Acknowledgments:}
We would like to thank many people for discussions which were
essential to our work, including Anton Alekseev, Alberto Cattaneo, Rui Loja Fernandes, Johannes Huebschmann,
Yvette Kosmann-Schwarzbach, Kirill Mackenzie, Eckhard Meinrenken, Pavol
\v{S}evera, Marco Zambon, and especially Ping Xu, who 
should be considered as a ``virtual co-author'' of this paper. 
Many institutions were our hosts during part of the work, including
MSRI and Universit\'e Libre de Bruxelles (Bursztyn), Institut
Math\'ematique de Jussieu (Weinstein), the Erwin Schr\"odinger Institute (Crainic).
Our collaboration was
facilitated by invitations to the conferences Poisson 2002 (Institute
Superior Tecnico, Lisbon) and Groupoidfest 2002 (University of Nevada, Reno).
For financial support, we would like to thank the National Science
Foundation (Bursztyn, grant DMS-9810361, Weinstein and Zhu, grant
DMS-0204100) and the Miller Institute for Basic Research in Science and the Dutch Royal Academy (Crainic).

\section{Basic definitions and the main results}
\label{sec-basic}
\subsection{Twisted Dirac structures}
\label{Twisted Dirac structures}

We recall some basic concepts in Dirac geometry \cite{Cour}. 
Let $V$ be a finite dimensional
vector space, and equip $V \oplus V^*$ with the symmetric pairing
\begin{equation}
\SP{(x,\xi),(y,\eta)}_{+} = \xi(y) + \eta(x).
\end{equation}
A {\bf linear Dirac structure} on $V$ is a subspace
$L \subset V \oplus V^*$ which is maximally isotropic with respect to $\SP{\; , \;}_{+}$.
The set of all Dirac structures on $V$ is a smooth submanifold of a Grassmann manifold;
we denote it by $\Dir(V)$.

Natural examples of vector spaces carrying linear Dirac structures are
presymplectic and Poisson vector spaces (i.e, spaces equipped with a skew-symmetric
bilinear form and a bivector, respectively). More precisely, a skew-symmetric
bilinear form $\theta$ (resp. bivector $\pi$) corresponds to the Dirac structure
$L_{\theta}$ (resp. $L_{\pi}$) given by the graph of the map $\widetilde{\theta}: V \rmap V^*$,
$\widetilde{\theta}(x)(y) = \theta(x,y)$ (resp. $\widetilde{\pi}: V^* \rmap V$, 
$\widetilde{\pi}(\alpha)(\beta) = \pi(\beta, \alpha)$). 

General Dirac structures
can be described
either in terms of bilinear forms or bivectors.
Let $pr_1:V\oplus V^* \longrightarrow V$ and $pr_2:V\oplus V^* \longrightarrow V^*$ be the natural
projections.
A linear Dirac structure $L$ has an associated range,
\[ 
\Ran(L)= pr_1(L)= \{ v\in V: (v, \xi)\in L\ \text{for\ some}\ \xi\in V^*\}\subset V,
\]
and a skew-symmetric bilinear form $\theta_{L}$ on $\Ran(L)$,
\begin{equation}\label{eq:thetaL}
\theta_L(v_1, v_2)= \xi_1(v_2), \ \ \text{where}\ \xi_1\in V^*\ \text{is\ such\ that}\ (v_1, \xi_1)\in L.
\end{equation}
For the description of $L$ in terms of  a bivector, let us define the {\bf kernel} of $L$ as the kernel
of $\theta_L$:
\[ 
\Ker(L) := \Ker(\theta_L) = pr_2(L)^{\circ} = \{ v\in V: (v, 0)\in L\}\subset V .
\]
(Here $^\circ$ stands for the annihilator.)
The bivector $\pi_L$, defined on $V/\Ker(L)$, is the one induced by a form analogous to 
\eqref{eq:thetaL} on $pr_2(L)\subset V^*$.
It is not difficult to see that $L$ is completely characterized by the pair $(\Ran(L), \theta_L)$,
or, analogously, by the pair $(\Ker(L), \pi_L)$. 
We observe that

\begin{enumerate}
\item[(i)] $R(L)= V$ if and only if $L= L_{\theta}$ for some skew-symmetric
bilinear form $\theta$ on $V$;
\item[(ii)] $\Ker(L)= 0$ if and only if $L= L_{\pi}$ for some bivector $\pi$ on $V$.
\end{enumerate}

If $V$ and $W$ are vector spaces, any linear map $\psi: V \longrightarrow W$
induces a push-forward map $\For\psi : \Dir(V) \longrightarrow \Dir(W)$ by
\begin{equation} \label{eq:forward}
\For\psi (L_{V}) = \{ (\psi(x),\eta) \; |\; x \in V, \, \eta \in W^*, \, (x,\psi^*(\eta)) \in
                      L_{V} \},
\end{equation}
where $L_{V} \in \Dir(V)$.  We note that the map $\For \psi$ is
{\it not} continuous at every $L_{V}$. 

Given $L_{V} \in \Dir(V)$ and $L_{W} \in \Dir(W)$, we call
a linear map $\psi: V \rmap W$ {\bf forward Dirac} if $\For \psi(L_{V}) = L_W$. There
is a corresponding concept of a backward Dirac map (see e.g. \cite{BuRa}), but we will not 
deal with it in this paper.
Hence, for simplicity, we will refer to forward Dirac maps just as {\bf Dirac maps}. 
As an example, we recall that a Dirac map between Poisson vector spaces is just a Poisson
map.

An almost Dirac structure on a smooth manifold $M$ is a subbundle $L \subset TM \oplus T^*M$
defining a linear Dirac structure on each fiber. Note that the dimensions of the range 
$\Ran(L)$ and the kernel
$\Ker(L)$ (defined fiberwise) may vary from a point to another. 

A {\bf Dirac structure} 
is an almost Dirac structure whose sections are closed under
the Courant bracket\footnote{This is the non-skew-symmetric version of Courant's \cite{Cour}
original bracket, as introduced in \cite{LiuWeXu} and used in \cite{WeSe}.}
$[\; , \;] :
\Gamma(TM \oplus T^*M) \times \Gamma(TM \oplus T^*M) \rmap \Gamma(TM \oplus T^*M)$,
\begin{equation}\label{eq:ordcourant}
[(X,\xi),(Y,\eta)] = ([X,Y], \Lie_{X}\eta - i_{Y}d\xi).
\end{equation}
For example, a bivector field $\pi$ on $M$ corresponds to a Dirac
structure if and only if it defines a  
Poisson structure.

As observed in \cite{WeSe}, one can use a closed 3-form $\phi$ on $M$ to modify the 
standard Courant bracket as follows:
\begin{equation}\label{eq:twistcourant}
[(X,\xi),(Y,\eta)]_{\phi} = ([X,Y], \Lie_{X}\eta - i_{Y}d\xi + \phi(X,Y, \cdot)).
\end{equation}
A Dirac structure twisted by $\phi$, or simply a {\bf $\phi$-twisted Dirac structure}, is
an almost Dirac structure whose sections are closed under $[\cdot , \cdot ]_{\phi}$.
A bivector $\pi$ defines a $\phi$-twisted Dirac structure if and only if
$$
[\pi,\pi] = 2 (\wedge^3)\tilde{\pi}(\phi),
$$
where $[\, , \,]$ is the Schouten bracket. Similarly, a 2-form $\omega$ defines
a $\phi$-twisted Dirac structure if and only if $d\omega+ \phi= 0$.
These two kinds of Dirac structures are called {\bf $\phi$-twisted
  Poisson structures} and {\bf $\phi$-twisted presymplectic
  structures}.

Let $(M,L_{M},\phi_{M})$ and $(N, L_{N}, \phi_{N})$ be twisted Dirac manifolds.
A smooth map $\psi: M \rmap N$ is called a (forward) {\bf Dirac map} if
$\For (d\psi)_x ((L_M)_x) = (L_N)_{\psi(x)}$, for all $x \in M$.
(Here, and throughout this paper,
we denote by $(df)_{x}:T_xM\rmap T_{f(x)}N$ the differential of a smooth
function $f:M\rmap N$ at the point $x\in M$.) 

We observe that the push-forward operation between Dirac structures
defined in the linear case (\ref{eq:forward}) is not well defined for manifolds
in general. For instance, if $\psi: M \rmap N$ is a surjective submersion and $L_M$ is
a Dirac structure on $M$, the pointwise push-forward structures $\For d_x\psi ((L_M)_x)$
may differ for points $x$ along the same $\psi$-fiber; even if they coincide,
the resulting family of vector spaces might not be a subbundle, because
of discontinuities of the map $\For$.

\subsection{Groupoids and algebroids}
\label{Groupoids and algebroids}

Throughout the text, $G$ will denote a Lie group\-oid. We denote the
unit map by $\varepsilon: M\rmap G$, the inversion by $i: G \rmap G$
and the source (resp. target) map by $s: G\rmap M$ (resp. $t: G\rmap
M$). We denote the set of composable pairs of groupoid elements by
$G_2$ (adopting the convention that $(g, h) \in G\times G$ is
composable if $s(g)=t(h)$), and we write $m: G_2 \rmap G$ for the
multiplication operation.  We will often identify $M$ with the
submanifold of $G$ of identity arrows . In particular, given $x\in M$
and $v\in T_xM$, $\epsilon(x)= 1_x\in G$ will be
identified with $x$, and the tangent vector $(d\epsilon)_{x}(v)\in
T_xG$ with $v\in T_xM$.

We emphasize that the Lie groupoids we consider may be non-Hausdorff. 
Basic important examples come from foliation theory and from  integration of
bundles of Lie algebras. However, $M$ and the fibers of $s$ are
always assumed to be Hausdorff manifolds, and $s$-fibers will be assumed to be connected. 
We say that $G$
is $s$-simply connected if the $s$-fibers are simply connected. We will
use the notation $G(-, x)= s^{-1}(x)$ (arrows starting at $x \in M$), and, for 
$ y \stackrel{g}{\leftarrow} x$, we write the corresponding right multiplication map as
$R_g: G(-, y) \longrightarrow G(-, x)$, $R_g(a)= ag$.

The infinitesimal version of a Lie groupoid is a Lie algebroid. To fix our notation, we
recall that a Lie algebroid $A$ over $M$ is a vector bundle  $A$ over $M$ together with
a Lie bracket $[\cdot, \cdot ]$ on the space of sections $\Gamma(A)$, and a bundle
map $\rho: A\rmap TM$ satisfying the Leibniz rule
\[ 
[\alpha, f\beta]= f[\alpha, \beta]+ \Lie_{\alpha}(f)\beta.
\]
Here and elsewhere in this paper, we use the notation $\Lie_{\alpha}$ for the Lie derivative
with respect to the vector field $\rho(\alpha)$.

Given a Lie groupoid $G$, the associated Lie algebroid $A= \mathrm{Lie}(G)$
has fibers $A_x= \Ker(ds)_x= T_x(G(-, x))$, for $x \in M$. Any $\alpha\in \Gamma(A)$ 
extends to a unique right-invariant vector field on $G$, which will be 
denoted by the same letter $\alpha$. This correspondence identifies $\Gamma(A)$ 
with the space $\mathcal{X}^{s}_{inv}(G)$ of vector fields on $G$ which are tangent to the $s$-fibers and 
right invariant. The usual Lie bracket on vector fields induces  the bracket on $\Gamma(A)$,
and the anchor is given by $\rho = dt: A \rmap TM$.

Given a Lie algebroid $A$, an {\bf integration} of $A$ is a Lie groupoid $G$ {\it together
with} an isomorphism $A\cong \mathrm{Lie}(G)$. By abuse of language, we will often call $G$
alone an integration of $A$.  If such a $G$ exists, we say that $A$ is {\bf integrable}.
In contrast with the case of Lie algebras, not all
Lie algebroids are integrable (see \cite{CrFe} and references therein). However, as 
with Lie algebras, integrability implies the existence  
of a canonical $s$-simply connected integration $G(A)$ of $A$,
and any other $s$-simply connected integration of $A$ will be isomorphic to
$G(A)$. Roughly speaking, $G(A)$ consists of $A$-homotopy classes of $A$-paths,
where an $A$-path consists of a path $\gamma: I\rmap M$ together 
with an ``$A$-derivative of $\gamma$'', i.e. a path $a: I\rmap A$ above $\gamma$
with the property that $\rho(a)$ is the usual derivative $\frac{d\gamma}{dt}$.
In general, $G(A)$ is a topological groupoid carrying an additional
``smooth structure'' as the leaf space of a foliation, 
called in \cite{CrFe} the {\it  Weinstein groupoid of} $A$, and 
its being a manifold is equivalent to the integrability of $A$. 
Further details and the precise obstructions to integrability can be found in \cite{CrFe},
though some facts about $G(A)$ will be recalled in Section \ref{Reconstructing multiplicative forms}.

A $\phi$-twisted Dirac structure $L$ has an induced Lie algebroid structure \cite{Cour,LiuWeXu,WeSe}:
the bracket on $\Gamma(L)$ is the restriction of the Courant bracket $[\,\cdot,\cdot]_{\phi}$, 
and the anchor is
the restriction of the projection $pr_1:TM\oplus T^*M \rmap TM$.
We denote by $G(L)$ the groupoid associated
to this Lie algebroid structure. The main theme in this paper is the
description of the extra structure 
on $G(L)$ induced by the Dirac structure. For instance, if $L= L_{\pi}$ is
the Dirac structure coming from a Poisson tensor $\pi\in \Gamma(\Lambda^2TM)$, 
then $G(L_{\pi})$ carries a canonical symplectic structure making it into a {\it symplectic groupoid}
(also interpreted as the phase space of an associated Poisson
sigma-model, see \cite{CaFe, CrFe2}). If $\pi$ is a twisted Poisson
structure, $G(L_{\pi})$ becomes what Cattaneo and Xu \cite{CaXu} call
a non-degenerate quasi-symplectic groupoid (and we call a twisted
symplectic groupoid).

%%%%%%%%%%%%%%%%%%%%%%%%%%%%%%%%%%%%%%%%%%%%%%%%%%%%%%%%%%%%%%%%%%%%%%%%%%%%%%%%%%%%%%%%%%%%%%%%%%
\subsection{Dirac $\leftrightarrow$ presymplectic}
\label{Dirac-pre}
Symplectic structures appear in several ways in connection with Poisson manifolds,
e.g., as symplectic leaves, symplectic realizations and symplectic groupoids.
With the relationship (Poisson manifolds) $\leftrightarrow$ (symplectic structures) in mind,
we briefly discuss in the remainder of this section
 the analogous correspondence for (twisted)
Dirac manifolds.

%	\begin{center}
%	\begin{tabular}{|p{2.6cm}||p{2.6cm}|p{2.6cm}|p{2.6cm}|}
%\hline
%Poisson manifolds & symplectic leaves & symplectic
%\mbox{realizations} & symplectic groupoids \\
%\hline 
%Dirac manifolds &  presymplectic leaves & presymplectic realizations & presymplectic groupoids \\
%\hline
%	\end{tabular}
%	\end{center}

Recall that a {\bf presymplectic manifold} is a manifold $M$
equipped with a closed 2-form $\omega$. When $\Ker(\omega)$ has constant rank, it defines
a foliation on $M$; if this foliation is {\it simple} (i.e., the space of leaves $M/\Ker(\omega)$ is smooth
and the quotient map is a submersion), then $M/\Ker(\omega)$ becomes a
symplectic manifold, with symplectic form induced by $\omega$. Hence, 
 modulo global regularity issues, presymplectic manifolds can be reduced to symplectic
manifolds. In the case of a  $\phi$-twisted presymplectic manifold,
the reduction mentioned above works only when $\Ker(\omega)\subset \Ker(\phi)$. 

Just as Poisson structures can be viewed as singular foliations whose leaves are
symplectic manifolds, $\phi$-twisted Dirac structures $L$ are singular foliations
whose leaves are $\phi$-twisted presymplectic manifolds. Given $L$,
the foliation is defined at each $x \in M$ by its range $\Ran(L_x)$.
Note that this (singular) distribution coincides with the image of the anchor map of the  Lie algebroid
structure of $L$, and hence is necessarily integrable. 
In particular, the leaves of the foliation
are the orbits of this Lie algebroid. For each such leaf $S$, the $\phi|_{S}$-presymplectic form
$\theta_{S}\in \Omega^2(S)$ is, at each point,  just the 2-form $\theta_{L_{x}}$ associated to the linear Dirac space
$L_{x}$. As above, under certain regularity 
conditions (e.g. if $\Ker(\omega)\subset \Ker(\phi)$ and $\Ker(L)$ is simple) 
one can quotient out $M$ by $\Ker(L)$
and reduce $(M, L)$ to a twisted Poisson manifold $M/\Ker(L)$ whose symplectic leaves are precisely
the reductions of the presymplectic leaves of $L$. 
This suggests that Dirac structures could very well be called ``pre-Poisson structures''.

The notion of a {\bf presymplectic realization} of a $\phi$-twisted Dirac structure 
$L$ on $M$ is more subtle:
it is a Dirac map
\[ 
\mu: (P, \eta)\rmap (M, L) ,
\]
where $\eta\in \Omega^2(P)$ is a $\mu^*\phi$-twisted presymplectic form ($d\eta +\mu^*\phi=0$), 
with the extra property that $\Ker(d\mu)\cap \Ker(\eta)=\{0\}$. 
This ``non-degeneracy condition'' for $\eta$  will be explained in more detail
in  Section \ref{examples}. As we will see, presymplectic realizations are the 
{\it infinitesimal} counterpart
of the actions studied in \cite{Xu}, from where much of our motivation comes.

The only thing still to be explained in the correspondence between Dirac structures and presymplectic structures 
are presymplectic groupoids, and this leads us to the main results of the paper.

\subsection{The main results}

A 2-form $\omega$ on a Lie groupoid $G$ is called {\bf multiplicative} if the
graph of $m: G_2 \rmap G$ is an isotropic submanifold of 
$(G,\omega)\times (G, \omega) \times (G, -\omega)$, or, equivalently, if
\begin{equation}\label{multiplicative}
m^*\omega = pr_1^*\omega + pr_2^*\omega,
\end{equation}
where $pr_i:G_2 \rmap G$, $i=1,2$, are the natural projections.

Let $G$ be a Lie groupoid over $M$, $\phi$ a closed $3$-form on $M$, and
$\omega$ a multiplicative 2-form on $G$. We call $\omega$
{\bf relatively $\phi$-closed} if $d\omega = s^*\phi - t^*\phi$.

\begin{definition}\label{def:twistpresymp}
We call $(G,\omega,\phi)$ a  presymplectic groupoid twisted by $\phi$, or a
{\bf $\phi$-twisted presymplectic groupoid}, if $\omega$ is relatively $\phi$-closed,
$\mathrm{dim}(G) = 2~\mathrm{dim}(M)$, and if
\begin{equation}\label{eq:nondeg}
\Ker(\omega_x)\cap \Ker(ds)_x\cap \Ker(dt)_x= \{ 0\},
\end{equation}
for all $x \in M$.
\end{definition}

Condition (\ref{eq:nondeg}) provides a restriction on how degenerate $\omega$ can be;
we will discuss it further in Section \ref{sec:nondeg}.
If $\phi=0$ and $\omega$ is nondegenerate, the groupoid in Definition~\ref{def:twistpresymp}
is just a symplectic groupoid in the usual sense \cite{Wein}.

We  find that, modulo integrability issues,
(twisted) Dirac structures on $M$ are basically the same thing as (twisted) presymplectic
groupoids over $M$. This extends the relationship between
integrable Poisson manifolds and symplectic groupoids (see e.g. \cite{CDW}). 
More precisely, we have two main results on
this correspondence.  The first starts with the groupoid.

\begin{theorem}
\label{theorem1}
Let $(G,\omega,\phi)$ be a $\phi$-twisted presymplectic groupoid. Then
\begin{enumerate}
% \item[(i)] $(\Ker(dt)_g)^{\perp}= \Ker (ds)_g+ \Ker(\omega_g)$ for all $g\in G$;
% \item[(ii)] $(TM)^{\perp}= TM+ \Ker(\omega)$;
\item[(i)] There is a (canonical, and unique) $\phi$-twisted Dirac structure $L$ on $M$, 
such that $t: G\rmap M$ is a Dirac map, while $s$ is anti-Dirac; 
\item[(ii)] There is an induced isomorphism between 
the Lie algebroid $\mathrm{Lie}(G)$ of $G$ and the Lie algebroid of $L$.
\end{enumerate}
\end{theorem}

We will also see that known properties of symplectic groupoids naturally  extend to our setting.
For instance, $(\Ker(dt)_g)^{\perp}= \Ker (ds)_g+ \Ker(\omega_g)$ for all $g\in G$, and 
$(TM)^{\perp}= TM+ \Ker(\omega)$ (where $\perp$ denotes the orthogonal with respect to $\omega$). 

\begin{remark}\rm
The alternative convention of identifying the Lie algebroid of $G$ with left-invariant
vector fields would lead to $s$ being Dirac and $t$ being anti-Dirac. This is
the convention adopted in \cite{CDW,Wein}.
\end{remark}

In the situation of the theorem, we say that $(G, \omega)$ is an {\bf integration} of
the $\phi$-twisted Dirac structure $L$. Note that such an integration immediately gives rise
to an integration of the Lie algebroid associated with $L$ (in the sense of  
subsection \ref{Groupoids and algebroids}), namely the Lie groupoid $G$
together with the isomorphism insured by Theorem~\ref{theorem1}, part (ii) 
(which does depend on $\omega$!).

Our second result starts with a Dirac structure.

% \begin{theorem}
% \label{theorem2}
% Let $L$ be a $\phi$-twisted Dirac structure on $M$ whose associated Lie algebroid is integrable.
% Then any $s$-simply connected integration $G$ of this Lie algebroid (which does exists,
% and is unique up to isomorphism), can be extended uniquely to an integration
% $(G, \omega)$ of the $\phi$-twisted Dirac structure $L$.
% \end{theorem}

\begin{theorem}
\label{theorem2}
Let $L$ be a $\phi$-twisted Dirac structure on $M$ whose associated Lie algebroid is integrable,
and let $G(L)$ be its $s$-simply connected integration. Then there exists a unique 2-form $\omega_L$
such that $(G(L), \omega_{L})$ is an integration of the $\phi$-twisted Dirac structure $L$.
\end{theorem}

Hence we obtain a one-to-one correspondence between integrable $\phi$-twisted Dirac structures on $M$
and $\phi$-twisted presymplectic groupoids over $M$ by
\[ 
L \leftrightarrow (G(L), \omega_{L}). 
\]

In order to prove these two theorems, we have to understand the intricacies of
multiplicative 2-forms, starting with their infinitesimal 
counterpart. 
We will prove the following result, which we expect to be useful in other settings as well.

\begin{theorem}
\label{theorem3}
Let $G$ be an $s$-simply connected Lie groupoid over $M$, with Lie algebroid $A$,
and let $\phi\in \Omega^3(M)$ be a closed 3-form. Then there is a one-to-one correspondence between 
\begin{enumerate}
\item[(i)] multiplicative 2-forms $\omega\in \Omega^2(G)$, with $d\omega = s^*\phi - t^*\phi$.
\item[(ii)] bundle maps $\rho^*:A\rmap T^*M$ with the properties:
\[  \begin{split}
\SP{\rho^{*}(\alpha),\rho(\beta)} &  =  -\SP{\rho^{*}(\beta),\rho(\alpha)};\\
\rho^{*}([\alpha, \beta])  =   \mathcal{L}_{\alpha}(\rho^{*}(\beta)) - 
                                      & \mathcal{L}_{\beta}(\rho^{*}(\alpha))+ 
                                      d\SP{\rho^*(\alpha),\rho(\beta)} + 
                                      i_{\rho(\alpha)\wedge \rho(\beta)}(\phi).
\end{split}  \]
\end{enumerate}
In fact, for a given $\omega$, the corresponding $\rho^*$ is $\rho^{*}_{\omega}(\alpha)(X)= \omega(\alpha, X)$.
\end{theorem}

%%%%%%%%%%%%%%%%%%%%%%%%%%%%%%%%%%%%%%%
%%%%%%%%%%%%%%%%%%%%%%%%%%%%%%%%%%%%%%%
\section{Multiplicative 2-forms on groupoids}
\label{sec-mult}
%%%%%%%%%%%%%%%%%%%%%%%%%%%%%%%%%%%%%%%
%%%%%%%%%%%%%%%%%%%%%%%%%%%%%%%%%%%%%%%
%%%%%%%%%%%%%%%%%%%%%%%%%%%%%%%%%%%%%%%

In this section we discuss  general properties of multiplicative 
2-forms on Lie groupoids.

\begin{lemma} 
\label{lemma}
If $\omega\in \Omega^2(G)$ is multiplicative, then 
\begin{enumerate} 
\item[(i)] $\varepsilon^*\omega= 0$, and $i^*\omega= -\omega$;
\item[(ii)] $\Ker(ds)_g+ \Ker(\omega_g)\subset (\Ker(dt)_{g})^{\perp}$ for all $g\in G$;
\item[(iii)] For all arrows $g: y \longleftarrow x $, the map $(di)_{g}$ induces an isomorphism
 \[ \Ker(\omega_g) \longrightarrow \Ker(\omega_{g^{-1}}) ,\]
and $(dR_{g})_{y}$ induces isomorphisms
 \[ 
 \Ker(\omega_y)\cap \Ker(ds)_y \longrightarrow \Ker(\omega_g)\cap 
 \Ker(ds)_g,
 \] 
 \[ 
 \Ker(\omega_y)\cap \Ker(ds)_y\cap \Ker(dt)_y \longrightarrow \Ker(\omega_g)\cap 
 \Ker(ds)_g\cap \Ker(dt)_g, 
 \]
 \[
 \Ker(ds)_y\cap(\Ker(ds)_y)^\perp \longrightarrow \Ker(ds)_g\cap(\Ker(ds)_g)^\perp.
 \]
\item[(iv)] On each orbit $S$ of $G$, there is an induced 2-form $\theta_S$,
uniquely determined by the formula $\omega|_{G_S}= t^*\theta_S- s^*\theta_S$, 
where $G_S= s^{-1}(S)= t^{-1}(S)$ is the restriction of $G$ to $S$. Moreover, if $\omega$
is relatively $\phi$-closed, then $d\theta=-\phi|_{S}$; hence each orbit of $G$ 
becomes a $(\phi|_S)$-twisted presymplectic manifold.
\end{enumerate}
\end{lemma}
Note that, in (iii), $\Ker(\omega_y)$ and $\Ker(\omega_g)$ are not isomorphic in general.

\begin{remark}\label{rem:st}
The statements in (ii) and (iii) clearly hold with the roles of $s$ and $t$ interchanged (and right multiplication
replaced by left multiplication). 
\end{remark}

\begin{proof}
We will need the following simple identities:
\begin{equation}
\label{eq1} v_{g}= (dm)_{y, g}( (dt)_g(v_g), v_g)= (dm)_{g,x}(v_g, (ds)_g(v_g)), \ \text{for\ all}\ v_g\in T_gG,
\end{equation}
\begin{equation}
\label{eq2} (dR_g)_y(\alpha_y)= (dm)_{y, g}(\alpha_y, 0), \ \text{for\ all}\ \alpha_y\in \Ker(ds)_y.
\end{equation}
\begin{equation}
\label{eq3} (dt)_g(v_g)= (dm)_{g, g^{-1}}(v_g, (di)_g(v_g)), \ \text{for\ all}\ v_g\in T_gG.
\end{equation}

The identity in (\ref{eq1}) is obtained by differentiating  $id_{G}= m\circ (t, id_G)= m\circ (id_{G}, s)$.
We verify the other two formulas similarly: for (\ref{eq2}), we write
$R_{g}: G(-,y)\rmap G$ as $a\mapsto (a, g)\stackrel{m}{\mapsto} ag$, while for (\ref{eq3})  we write
the target map $t$ as the composition of $G\rmap G\times_{M}G$, $g\mapsto (g, g^{-1})$ with $m$.

We now prove the lemma.  
Consider the map $(id\times i): G\rmap G\times_{M}G$, $\, g\mapsto (g, g^{-1})$.
Applying $(id\times i)^*$ to (\ref{multiplicative}) and using (\ref{eq3}), we deduce that 
$t^*\varepsilon^*\omega= \omega+ i^*\omega$. Now applying $\varepsilon^*$ to this equation, we get 
$\varepsilon^*\omega= \varepsilon^*\omega+ \varepsilon^*\omega = 0$, and therefore 
$\omega+ i^*\omega= 0$. This proves (i).

Using (\ref{eq1}), (\ref{eq2}), and the multiplicativity of $\omega$,
we get
\begin{equation}
\label{eq4} 
\omega_g( (dR_g)_y(\alpha_y), v_{g})= \omega_y( \alpha_y,(dt)_g(v_g)) 
\end{equation}
for all $\alpha_y\in \Ker(ds)_y$ and $v_g\in T_gG$. When $v_g\in \Ker(dt)_g$,
since $(dR_g)_y$ maps $\Ker (ds)_y$ isomorphically into $\Ker (ds)_g$, it 
follows that $\Ker (ds)_g\subset (\Ker(dt)_{g})^{\perp}$. Since $\Ker(\omega)$ is 
inside all orthogonals, (ii) follows.

The first isomorphism in (iii) follows from $i^*\omega= -\omega$. In order to check the following two, 
note that $(dR_g)_y$ maps $\Ker(ds)_y$ isomorphically onto $\Ker(ds)_g$, and 
$\Ker(ds)_y\cap \Ker(dt)_y$ isomorphically onto $\Ker(ds)_g\cap \Ker(dt)_g$. 
So it suffices to prove the first isomorphism induced by $(dR_g)_y$ (as the second follows from it).
Let $\alpha_y\in \Ker(ds)_y$. Using (\ref{eq2}) and the multiplicativity of $\omega$, we have
\[ 
\omega_g((dR_g)_y(\alpha_y), (dm)_{y, g}(v_y, w_g))= \omega_y(\alpha_y, v_y) 
\]
for all $(v_y, w_g)$ tangent to the graph of the multiplication, 
which shows that $\alpha_y\in \Ker(ds)_y\cap \Ker(\omega_y)$
if and only if $(dR_g)_y(\alpha_y)\in \Ker(\omega_g)$. 
The last isomorphism in (iii) is implied by
$$
\omega_g((dR_g)_y(\alpha_y), (dR_g)_y(u_y)) = \omega_y(\alpha_y,u_y)
$$ 
for all $\alpha_y, u_y \in \Ker(ds)_y$, which follows from (\ref{eq2}) and the multiplicativity of $\omega$. 

Part (iv) is a statement about transitive groupoids (namely $G|_{S}$), so
we may assume that $S= M$ and $G$ is transitive.
If $G$ is the pair groupoid $M\times M$ 
over $M$, the proof is straightforward.
Indeed, a multiplicative 2-form $\omega$ satisfies
\begin{equation}\label{eq:pairmult}
\omega_{(x,y)}((u_x,u_y),(v_x,v_y))=
\omega_{(x,z)}((u_x,w_z),(v_x,w'_z)) +
\omega_{(z,y)}((w_z,u_y),(w'_z,v_y)),
\end{equation}
for $u_x,v_x \in T_xM$, $u_y,v_y \in T_yM$ and $w_z,w'_z \in T_zM$.
In particular, we can write 
\begin{equation}\label{eq:pairmult2}
\omega_{(x,y)}((u_x,u_y),(v_x,v_y))=
\omega_{(x,z)}((u_x,0_z),(v_x,0_z)) +
\omega_{(z,y)}((0_z,u_y),(0_z,v_y)).
\end{equation}
A direct application of \eqref{eq:pairmult} shows that the first factor
in the r.h.s. of \eqref{eq:pairmult2} is of the form $t^*\theta$, and the second is of the form $s^*\theta'$, for $\theta, \theta' \in \Omega^2(M)$.
The fact that $\theta = -\theta'$ follows from the first assertion of part (i).

In general, a transitive 
groupoid can be written as $G= (P\times P)/K$, the quotient of
the pair groupoid $G(P)= P\times P$ by the action of a Lie group $K$,
where $P\rmap M$ is a principal $K$-bundle (fix $x\in M$ and take $P$ as the set of arrows
starting at $x$). Let $p_1: G(P)\rmap G$ and $p_2: P \rmap M$ be the natural projections.
(By abuse of language, we denote source and target maps on either $G(P)$ or $G$ by $s$ and $t$,
since the context should avoid any confusion.)
If $\omega\in \Omega^2(G)$ is multiplicative, then so is $p_1^*\omega$, and we can write 
$p_1^*\omega= t^*\theta_0- s^*\theta_0$, for some $\theta_0 \in \Omega^2(P)$. 
Since the correspondence $\theta_0 \mapsto 
t^*\theta_0- s^*\theta_0$ is injective, the fact that $p_1^*\omega$ is 
basic implies that $\theta_0$ is basic. Hence $\theta_0= 
p_2^*\theta$ for some $\theta \in \Omega^2(M)$. Since $p_1$ is a 
submersion, it follows that $\omega= t^*\theta - s^*\theta$. Similarly, 
$s^*\phi- t^*\phi= d\omega= t^{*}d\theta- s^{*}d\theta$ implies that $d\theta= -\phi$.
\end{proof}

We now look at what happens at points $x\in M$. The following is a first 
sign of the rigidity of multiplicative 2-forms.

\begin{lemma}
\label{lemma2}
If $\omega\in \Omega^2(G)$ is multiplicative, then, at points $x\in M$,
\[
\begin{split}
& \Ker(ds)_x + \Ker(\omega_x)  = (\Ker(dt)_x)^{\perp}\\
& T_xM + \Ker(\omega_x)  =(T_xM)^{\perp}\\
& \Ker(\omega_x) = \Ker(\omega_x)\cap \Ker(ds)_x \oplus \Ker(\omega_x)\cap T_xM
\end{split}
\] 
(The first and third identities also hold for $s$ and $t$ interchanged.)

Furthermore, the following identities hold:
\[ \left \{ \begin{array}{ll} 
          \dim (\Ker(\omega_x)\cap T_xM)= \frac{1}{2}(\dim(\Ker(\omega_x))+ 2\dim(M)- \dim(G))\\
         \dim (\Ker(\omega_x)\cap \Ker(ds)_x)= \frac{1}{2}(\dim(\Ker(\omega_x))- 2\dim(M)+ \dim(G))
           \end{array}
  \right. 
\]
\end{lemma}

\begin{proof}
Let us first prove the third equality. Let  $u_x\in \Ker(\omega_x)$, and write
\[ u_x= (u_x- (ds)_x(u_x))+ (ds)_x(u_x).\]
It then suffices to show
that $(ds)_x(u_x)\in \Ker(\omega_x)$. Using (\ref{eq1}) 
(the one involving $ds$) and the multiplicativity of $\omega$, we get
\[ \omega_x(u_x, (dm)_{x, x}(v_x, w_x))= \omega_x(u_x, v_x)+ \omega_x((ds)_x(u_x), w_x) \]
for all $(v_x, w_x)$ tangent to the graph of $m$ at $(x, x)$. This shows
that, indeed, $\omega_x((ds)_x(u_x), w_x)= 0$ for all $w_x\in T_xG$.
Now, for the first two equalities, note that the direct inclusions are consequences 
of Lemma \ref{lemma}. We now compute the dimensions of the spaces involved.
First recall that for any subspace 
$W$ of a linear presymplectic space $(V, \omega)$, we have 
\begin{equation}\label{eq:presymport}
\dim(W^{\perp})= \dim(V)- \dim(W)+ \dim(W\cap \Ker(\omega)).
\end{equation}
Then, comparing dimensions in  $\Ker(ds) + \Ker(\omega) \subset (\Ker(dt))^\perp$,  we get
\begin{eqnarray*}
\dim(\Ker(\omega))& \leq&  \dim(\Ker(\omega)\cap \Ker(ds))+ \dim(\Ker(\omega)\cap \Ker(dt)) +\\
                  & &   + 2 \dim(M) - \dim(G) 
\end{eqnarray*}
at all $g \in G$. At a point $x \in M$, $$\dim(\Ker(\omega)\cap
\Ker(ds)) = \dim(\Ker(\omega)\cap \Ker(dt))$$ 
(since $(di)_x$ is an isomorphism between these spaces), so
\begin{equation}
\label{strange}
\dim(\Ker(\omega_x)) \leq 2\dim(\Ker(\omega_x)\cap \Ker(ds)_x)
             + 2 \dim(M) - \dim(G).
\end{equation}
Comparing dimensions in $T_xM + \Ker(\omega_x) \subset (T_xM)^{\perp}$, we get
\[ 
\dim(\Ker(\omega_x))\leq \dim(\Ker(\omega_x)\cap T_xM)+ \dim(G)- 2\dim(M) .
\]
Since we already proved the third equality in the statement, we do know that 
\begin{equation}
\label{not-strange}
\dim(\Ker(\omega_x)\cap T_xM)=  \dim(\Ker(\omega_x))- \dim(\Ker(\omega_x)\cap \Ker(ds)_x
\end{equation}
Plugging this into the last inequality, we get precisely the opposite of (\ref{strange}). This shows that
the inequality (\ref{strange}), as well as the direct inclusions for the first two relations in the
statement, must become equalities. This immediately implies the dimension identities in the statement
of the lemma and  completes the proof.
\end{proof}

Note that the first and third identities in the lemma immediately imply that
\begin{equation}\label{eq:keromega}
(\Ker(dt)_x)^\perp = \Ker(ds)_x + \Ker(\omega_x)\cap T_xM, \;\; x\in M.
\end{equation}

\begin{corollary} 
\label{unique}
Two multiplicative forms $\omega, \omega^{'}\in \Omega^2(G)$
which have the same differential $d\omega= d\omega'$, and which coincide
at all $x\in M$, must coincide globally.
\end{corollary}

\begin{proof}
We may clearly assume that $\omega'= 0$. Equation (\ref{eq4})
implies that $i_v(\omega)= 0$ for all $v$ tangent to the $s$-fibers.
Since $\omega$ is closed, it follows that $\mathcal{L}_v(\omega)= 0$ for all
such $v$'s, and therefore $\omega= s^*\eta$ for some 
2-form $\eta$ on $M$. Since $\omega|_{TM}= 0$, we must have $\eta= 0$, hence $\omega= 0$.
\end{proof}

We now discuss the infinitesimal counterpart of multiplicative forms.
A multiplicative 2-form $\omega\in \Omega^2(G)$ induces a bundle map
$\rho^{*}_{\omega}: A\rmap T^*M $
characterized by the equation
\begin{equation}\label{eq:rho*} 
\rho^{*}_{\omega}(\alpha)\circ (dt)= i_{\alpha}(\omega).
\end{equation}
This equation is in fact equivalent to 
$\rho^{*}_{\omega}(\alpha)= i_{\alpha}(\omega)|_{M}$,
as a consequence of Lemma \ref{lemma}, part (ii).

\begin{proposition}
\label{rho} 
The bundle map $\rho_{\omega}$ associated to a multiplicative 2-form $\omega\in \Omega^2(G)$
has the following properties:
\begin{enumerate}
\item[(i)] for $\alpha, \beta \in \Gamma(A)$, we have
           $\SP{\rho^{*}_{\omega}(\alpha),\rho(\beta)} = 
                        -\SP{\rho^{*}_{\omega}(\beta),\rho(\alpha)}$;
\item[(ii)] if $\phi$ is a closed 3-form on $M$ and $d\omega = s^*\phi - t^*\phi$, then
\[ 
\rho^{*}_{\omega}([\alpha, \beta])=  
\mathcal{L}_{\alpha}(\rho^{*}_{\omega}(\beta)) - \mathcal{L}_{\beta}(\rho^{*}_{\omega}(\alpha))+ 
d\SP{\rho^*_{\omega}(\alpha),\rho(\beta)} + i_{\rho(\alpha)\wedge \rho(\beta)}(\phi),
\]
for all $\alpha, \beta\in \Gamma(A)$;

\item[(iii)] Two multiplicative forms $\omega_1$ and $\omega_2$ coincide if and only 
if $\rho^{*}_{\omega_1}= \rho^{*}_{\omega_2}$, and $d\omega_1= d\omega_2$.
\end{enumerate}
\end{proposition}

\begin{proof} 
By the comment preceding the proposition, (i) is clear, and (iii) is just Corollary \ref{unique}. 
For part (ii), we fix $v$ a vector field on $M$, and we prove that the right and left hand sides of (ii) 
applied to $v$ are the same. Let $\alpha, \beta\in \Gamma(A)$, and denote the vector
fields on $G$ tangent to the $s$-fibers induced by them by the same letters.
 We consider
\begin{eqnarray} 
d\omega(\alpha, \beta, \tilde{v})& = & 
-\omega([\alpha, \beta], \tilde{v})+ \omega([\alpha, \tilde{v}], \beta)- \omega([\beta, \tilde{v}], \alpha)+
\label{part1} \\
 &  & + \mathcal{L}_{\alpha}(\omega(\beta, \tilde{v}))- \mathcal{L}_{\beta}(\omega(\alpha, \tilde{v}))+
\mathcal{L}_{\tilde{v}}(\omega(\alpha, \beta)) \nonumber
\end{eqnarray}
at points $x\in M$, where 
$\tilde{v}$ is a vector field on $G$ to be chosen. 
We will pick $\tilde{v}$ so that it agrees with $v$ on $M$, and so that 
it is $t$-projectable, i.e. 
$$
(dt)_{g}(\tilde{v}_g)= \tilde{v}_{t(g)}= v_{t(g)}.
$$
Note that, at points $x\in M$,  
$$
d\omega(\alpha, \beta, \tilde{v})= (s^*\phi - t^*\phi)(\alpha, \beta, \tilde{v})= 
-\phi(\rho(\alpha),\rho(\beta),v).
$$ 
We claim that this observation and  (\ref{part1}) imply part (ii).
 
In order to check that, we need some remarks on $t$-projectable vector fields 
and $t$-projectable functions on $G$ (i.e. functions $f$ with $f(g)= f(t(g))$):

\begin{enumerate}
\item[1)] any vector field $v$ on $M$ admits (locally if $G$ is non-Hausdorff) a $t$-projectable
extension $\tilde{v}$ to $G$;
\item[2)] if $X$ is a vector field on $G$, and $f$ is a $t$-projectable function, then 
$\mathcal{L}_{X}(f)(x) = \mathcal{L}_{dt(X_x)}(f|_{M})(x)$ for all $x\in M$;
\item[3)] if $\alpha, \beta\in \Gamma(A)$, then $\omega(\alpha, \beta)$, as a function on $G$, is
$t$-projectable;
\item[4)] if $\alpha\in \Gamma(A)$ and $\tilde{v}$ is $t$-projectable, then $\omega(\alpha, \tilde{v})$ is
$t$-projectable;
\item[5)] if $\alpha\in \Gamma(A)$ and $\tilde{v}$ is $t$-projectable, then 
$(dt)_{x}[\alpha, \tilde{v}]_{x}= [\rho(\alpha), v]_x$ for all $x\in M$.
Note that, if $\tilde{v}$ extends $v$, then this implies that $\omega([\alpha, \tilde{v}], \beta)=
\omega([\rho(\alpha), v], \beta)$ at all $x\in M$,
and for all $\beta\in \Gamma(A)$ (as a result of $s$ and $t$-fibers being  orthogonal with respect to
$\omega$).
\end{enumerate}

Taking these remarks into account in (\ref{part1}), we immediately obtain
\begin{eqnarray*}
\rho^*_{\omega}([\alpha,\beta])(v) & = & i_{[\rho(\beta),v]}\rho^*_{\omega}(\alpha) - 
                                i_{[\rho(\alpha),v]}\rho^*_{\omega}(\beta) + 
                                \mathcal{L}_{\rho(\alpha)}i_v\rho^*_{\omega}(\beta)-\\
                           &&   \mathcal{L}_{\rho(\beta)}i_v\rho^*_{\omega}(\alpha) +
                                i_v d(\SP{\rho^*_{\omega}(\beta),\rho(\alpha)}) +
                                \phi(\rho(\alpha), \rho(\beta), v).
\end{eqnarray*}
(Note that, in the statement of (ii), we follow the convention that $\mathcal{L}_{\alpha}$
means $\mathcal{L}_{\rho(\alpha)}$, and the same for $\mathcal{L}_\beta$). 
Now, using the identity $i_{[X,Y]}=\mathcal{L}_X i_Y - i_Y \mathcal{L}_X$
on the first two terms of the right hand side, we obtain (ii).

Let us briefly discuss the remarks above: 

In order to prove 1), let $\sigma$ be a  splitting of the bundle map $dt : TG\rmap t^*TM$.
(In general, if $G$ is non-Hausdorff, we can only find such a splitting locally. Since
the formula to be proven is local, this is enough.) 
Then $u_{x}= v_{x}- \sigma_{x}(v_x)$ is in $\Ker(dt)_x$. Now we extend it to a vector 
field on $G$ tangent to the $t$-fibers by left translations,
and set $\tilde{v}_{g}= \sigma_g(v_{t(g)})+ u_{g}$; 

Remark  2) is clear, while 3)  and 4) follow from equation (\ref{eq4});

To check 5),
we first look at $(dt)_{x}[\alpha, \tilde{v}]_{x}$. We see that
\begin{equation}\label{eq:flow}
(dt)_x \left.\frac{d}{d\epsilon}\right |_{\epsilon= 0}
(d\Phi_{\alpha}^{\epsilon})_{\Phi_{\alpha}^{-\epsilon}(x)}(\tilde{v}_{\Phi_{\alpha}^{-\epsilon}(x)})=
\left.\frac{d}{d\epsilon}\right |_{\epsilon= 0} 
d(t\circ\Phi_{\alpha}^{\epsilon})_{\Phi_{\alpha}^{-\epsilon}(x)}(\tilde{v}_{\Phi_{\alpha}^{-\epsilon}(x)}) ,
\end{equation}
where $\Phi_{\alpha}^{\epsilon}$ is the flow of $\alpha$ viewed as a vector field on $G$
(extended by right translation). 
We have $t\circ \Phi_{\alpha}^{\epsilon}= \Phi_{\rho(\alpha)}^{\epsilon}\circ t$, and, since
$\tilde{v}$ is $t$-projectable,
$(dt)(\tilde{v}_{\Phi_{\alpha}^{-\epsilon}(x)})= v_{\Phi_{\rho(\alpha)}^{-\epsilon}(x)}$. Hence the last
term in (\ref{eq:flow}) equals to
\[ 
\left. \frac{d}{d\epsilon}\right |_{\epsilon= 0} 
d(t\circ\Phi_{\alpha}^{\epsilon})_{\Phi_{\alpha}^{-\epsilon}(x)}(v_{\Phi_{\rho(\alpha)}^{-\epsilon}(x)}) ,
\]
and this gives us $[\rho(\alpha), v]_{x}$. 
\end{proof}

\begin{remark}
\label{explicit-form}\rm
 Note that the bundle map $\rho^{*}_{\omega}$ (\ref{eq:rho*}) carries all the information
about $\omega$ at units: at $x\in M$ one has a canonical decomposition
\[ 
T_xG\cong T_xM \oplus A_{x},\;  v_x \mapsto (ds(v_x), v_x- (ds)_x(v_x)),
\]
and, with respect to this decomposition, it is easy to see that 
\begin{equation} \label{eq:decomp}
\omega((X, \alpha), (Y, \beta))= 
\rho^{*}_{\omega}(\alpha)(Y)- \rho^{*}_{\omega}(\beta)(X)+ \SP{\rho^*_{\omega}(\alpha),\rho(\beta)}.
\end{equation}
\end{remark}

%%%%%%%%%%%%%%%%%%%%%%%%%%%%%%%%%%%%%%%
%%%%%%%%%%%%%%%%%%%%%%%%%%%%%%%%%%%%%%%
%%%%%%%%%%%%%%%%%%%%%%%%%%%%%%%%%%%%%%%
%%%%%%%%%%%%%%%%%%%%%%%%%%%%%%%%%%%%%%%
\section{Multiplicative 2-forms and induced Dirac structures}
\label{sec:nondeg}
%%%%%%%%%%%%%%%%%%%%%%%%%%%%%%%%%%%%%%%
%%%%%%%%%%%%%%%%%%%%%%%%%%%%%%%%%%%%%%%
%%%%%%%%%%%%%%%%%%%%%%%%%%%%%%%%%%%%%%%
%%%%%%%%%%%%%%%%%%%%%%%%%%%%%%%%%%%%%%%

In this section we discuss the relationship between multiplicative 2-forms and Dirac structures,
proving in particular Theorem \ref{theorem1}. 

Let $G$ be a Lie groupoid over $M$, 
and let $\phi\in \Omega^3(M)$ be a closed 3-form.
Given $\omega\in \Omega^2(G)$ multiplicative, we first look at
 when the Dirac structure $L_{\omega}$ associated to $\omega$ can be 
linearly pushed forward by the target map $t$.  
For $g\in G$, let $t_*(L_{\omega, g})$ denote the push-forward of $L_{\omega, g}$ by $(dt)_{g}$:
\begin{equation}\label{eq:pushdownatg} 
t_*(L_{\omega, g})= \{ ((dt)_{g}(v_g), \xi_x): i_{v_g}(\omega)= \xi_{x}\circ (dt)_g\}\subset T_xM\oplus T_{x}^{*}M, 
\end{equation}
where $x= t(g)$. In particular, restricting \eqref{eq:pushdownatg} to points in $M$, we get a 
(possibly non-smooth, but of constant rank) bundle of 
linear Dirac structures $L_{M}$ on $M$: 
\[ 
L_{M, x}= \{ ((dt)_{x}(v_x), \xi_x): i_{v_x}(\omega)= \xi_{x}\circ (dt)_x\}\subset T_xM\oplus T_{x}^{*}(M). 
\]
The problem is to understand when this bundle agrees with \eqref{eq:pushdownatg} for all $g \in G$.

\begin{definition}
\label{def-diractype}
We say that a multiplicative 2-form $\omega$ on $G$ is of {\bf Dirac type} if 
$t_*(L_{\omega, g})= L_{M, t(g)}$ for all $g\in G$.
\end{definition}

The next result provides alternative characterizations of 2-forms of Dirac type.

\begin{lemma}
\label{lemma-Dirac}
Given a multiplicative 2-form $\omega$ on $G$, one has 
\begin{equation} 
\label{L-formula}
L_{M}= \{ (\rho(\alpha)+ u, \rho_{\omega}^*(\alpha)): \alpha\in A, u\in \Ker(\omega)\cap TM \} ,
\end{equation}
and the following are equivalent:
\begin{enumerate}
\item[(i)] $\omega$ is of Dirac type;
\item[(ii)] $\Ker(ds)_{g}^{\perp}  =  \Ker(dt)_{g}+ \Ker(\omega_{g})$ for all $g\in G$ ;
\item[(iii)] $(dt)_{g}: \Ker(\omega_g)\rmap \Ker(\omega_{t(g)})\cap TM$ is surjective  for all $g\in G$;
\item[(iv)] $\dim(\Ker(\omega_g))  = \frac{1}{2}(\dim(\Ker(\omega_x))+ \dim(\Ker(\omega_y)))$ for all $g\in G$.
\end{enumerate}
\end{lemma}

\begin{proof}
Let $x\in M$. Using the first and the third equalities in Lemma \ref{lemma2}, we see that
$\Ker(dt)_{x}^{\perp}= \Ker(ds)_{x}+ \Ker(\omega_{x})\cap TM$. Since the equation $i_{v}(\omega)= \xi\circ (dt)_x$
in the definition of $L_{M, x}$ implies that $v\in \Ker(dt)_{x}^{\perp}$, the elements of $L_{M, x}$
are pairs $((dt)_{x}(\alpha)+ u, \xi)$, where $\alpha\in \Ker(ds)_{x}= A_{x}$, $u\in \Ker(\omega)\cap TM$,
and  $i_{\alpha}(\omega)= \xi_{x}\circ (dt)_x$. Since this last equation is exactly the one defining  
$\rho_{\omega}^*(\alpha)$,
(\ref{L-formula}) is proven.
 
Let $g \in G$, with $s(g)= x$, $t(g)= y$. First, let us compute 
the codimension of $\Ker(dt)_{g}+ \Ker(\omega_{g})$ in $\Ker(ds)_{g}^{\perp}$. Using (\ref{eq:presymport}),
and the fact that $\Ker(ds)_g\cap \Ker(\omega_g)\cong \Ker(ds)_y\cap \Ker(\omega_y)$, and similarly for $t$
(cf. Lemma \ref{lemma}), we find that the codimension equals to
\begin{eqnarray*} 
 & 2\dim(M)- \dim(G)+ \dim(\Ker(\omega_x)\cap \Ker(dt)_x)+ \\
 & \dim(\Ker(\omega_y)\cap \Ker(ds)_y)- \dim(\Ker(\omega_g)).
\end{eqnarray*}
Using the last two formulas of Lemma \ref{lemma2}, we get
\begin{eqnarray}
\label{codim} 
  & \dim(\Ker(ds)_{g}^{\perp})- \dim(\Ker(dt)_{g}+ \Ker(\omega_{g})) =  - \dim(\Ker(\omega_g))+\\
\nonumber         &  \frac{1}{2}(\dim(\Ker(\omega_y))+ \dim(\Ker(\omega_x))).
\end{eqnarray}
Next, we compute the corank of the map in (iii) (note that 
$(dt)_g(\Ker(\omega_g)) \subseteq \Ker(\omega_y)\cap TM$ by \eqref{eq4}):
\[ 
\dim(\Ker(\omega_y\cap T_yM))- \dim(\Ker(\omega_g))+ \dim(\Ker(\omega_g)\cap \Ker(dt)_g) ,
\]
where, as above,  we may replace $g$ by $x$ in the last term; replacing the first and last terms 
by the last two formulas of Lemma \ref{lemma2}, we obtain
\begin{eqnarray} \label{codim2}
 & \dim(\Ker(\omega_y\cap T_xM))- \dim((dt)_{g}(\Ker(\omega_g)) = - \dim(Ker(\omega_g)) + \\
\nonumber & \frac{1}{2}(\dim(\Ker(\omega_y))+ \dim(\Ker(\omega_x))).
\end{eqnarray}
Equations (\ref{codim}) and \eqref{codim2} immediately imply the  equivalence of (ii), (iii), and (iv). 

To see that (i) implies (iii), assume  that $v\in \Ker(\omega_y)\cap TM$. Since  
$(v, 0)\in t_{*}(L_{\omega, y})$, we must have $(v, 0)\in t_{*}(L_{\omega, g})$,
i.e. $v= (dt)_{g}(w)$, with $w\in \Ker(\omega_g)$, and this proves (iii). 
For the converse, since $L_{M, y}$ and $t_{*}(L_{\omega, g})$ have the same dimension,
it suffices to prove that $L_{M, y}\subset t_{*}(L_{\omega, g})$.
Let $(\rho(\alpha)+ u, \rho_{\omega}^*(\alpha))$ be an element in $L_{M, y}$.
Write $\rho(\alpha)= (dt)_{g}(\tilde{\alpha})$, and $u= (dt)_{g}(\tilde{u})$,
where $\tilde{\alpha}= (dR_{g})_{y}(\alpha)$, and $\tilde{u}\in \Ker(\omega_g)$. 
By formula (\ref{eq4}) we have $i_{\tilde{\alpha}}(\omega)= i_{\alpha}(\omega)\circ (dt)_{g}=
\rho^*(\alpha)\circ (dt)_{g}$, and we now see that $(\rho(\alpha)+ u, \rho^*(\alpha))$
is in $t_{*}(L_{\omega, g})$.
\end{proof}

Note that, in contrast to
the Poisson bivectors on Poisson groupoids \cite{Wein3}, which always
satisfy a condition like that in Definition \ref{def-diractype},
a form $\omega$ can be multiplicative and closed without being of Dirac type.
(See Example \ref{non-dirac}.) 

It can also happen that $\omega$ is multiplicative and of Dirac type, but 
``the leaves of $L_{M}$'' do not coincide with the orbits of $G$ (take, e.g., $\omega$ to be zero on
a non-transitive groupoid).
When they do coincide, the situation becomes much more rigid, as we now discuss.
Recall that the orbits $S$ of $G$ are (twisted) presymplectic manifolds (see Lemma \ref{lemma}, (iv));
we denote the corresponding family of 2-forms by $\theta= \{ \theta_{S}\}$,
and we consider the bundle $L_{\theta}$ of linear Dirac structures on $M$ induced by $\theta$,
\[ 
L_{\theta, x}:= \{ (v_x, \xi_x): v_x\in T_xS, \; \xi_x|_{S}= i_{v_x}(\theta) \}\subset T_{x}M\oplus T^{*}_{x}M .
\]

\begin{lemma}
\label{lem-orbit}
Given a multiplicative 2-form $\omega\in \Omega^2(G)$, and $x\in M$, the following are equivalent:
\begin{enumerate}
\item[(i)] the range $\Ran(L_{M, x})$ equals to the tangent space to the orbit of $G$ through $x$;
\item[(ii)] the conormal bundle to the orbit through $x$, $(\Im(\rho_{x}))^{\circ}$, sits inside $\Im(\rho^{*}_{\omega})$;
\item[(iii)] $\Ker(\omega_x)\cap T_xM\subset \Im(\rho)$;
\item[(iv)] $\Ker(\omega_x)\cap T_xM= \Ker(\theta_x)$;
\item[(v)] $L_{M, x}= L_{\theta, x}$.
\end{enumerate}
Moreover, if these hold at all $x\in M$, then $\omega$ is of Dirac type. 
\end{lemma}

\begin{proof} 
Let us first show that
\begin{equation}
\label{kernels} 
\Ker(\theta_x)= \Ker(\omega_x)\cap \Im(\rho), \;\mbox{ and }\; \Ker(\omega_x)\cap T_xM= \Im(\rho^{*}_{\omega})^{\circ}.
\end{equation}
Given $v\in T_xS=\Im(\rho)$, we write $v= \rho(\alpha)$
with $\alpha\in \Ker(ds)_x$. From the defining property for $\theta= \theta_S$, we have
$\theta(v, w)= \theta(\rho(\alpha), \rho(\beta))= \omega(\alpha, \beta)$ for all $w= \rho(\beta)\in T_xS$.
Hence a vector $v$ is in $\Ker(\theta_x)$ if and only if $v= \rho(\alpha)$ has the property that
$\omega(\alpha, \beta)= 0$ for all $\beta\in \Ker(ds)_x$, i.e. 
$\alpha\in \Ker(ds)_x\cap (\Ker(ds)_{x})^{\perp}= 
\Ker(ds)_x\cap (\Ker(dt)_{x}+ \Ker(\omega_x)\cap T_xM)$, 
where for the last equality we used \eqref{eq:keromega}. Hence
\begin{eqnarray*}
\Ker(\theta_x)& = & (dt)_x(\Ker(ds)_x\cap (\Ker(ds)_{x})^{\perp}\\
              &=& (dt)_x(\Ker(ds)_x\cap (\Ker(dt)_{x}+ \Ker(\omega_x)\cap T_xM)), 
\end{eqnarray*}
which is easily seen to be $\Ker(\omega_x)\cap (dt_x)(\Ker(ds)_x)= \Ker(\omega_x)\cap \Im(\rho)$. 
The other identity
in (\ref{kernels}) easily follows from the explicit formula for $\omega_x$ in terms of $\rho$ and $\rho^{*}_{\omega}$
mentioned in Remark \ref{explicit-form}. Now note that (\ref{kernels}) proves the equivalence of (iii) and (iv),  
and also shows that (ii) and (iii) are just dual to each other. Hence we are left with proving that  
(i) is equivalent with (iii). But this is immediate since  
$\Ran(L_{M, x}) = \Im(\rho_{x})+ \Ker(\omega_x)\cap TM$ 
(e.g., cf. (\ref{L-formula})).

We now prove the last assertion of the lemma. We will do that by showing that $\Ker(ds)_g^{\perp} = \Ker(dt)_g +
\Ker(\omega_g)$ for all $g \in G$ (and using Lemma \ref{lemma-Dirac}, (ii)).
So let $g: y \leftarrow x$ be in $G$. 

{\it Claim $1$: Given $v\in T_{g}(G)$, one has $(dt)_{g}(v)\in \Ker(\omega_y)$ if and only if $v\in \Ker(ds)_{g}^{\perp}$.}

This follows immediately from (\ref{eq4}) and Lemma \eqref{lemma}, part (i).

{\it Claim 2: If (i)-(v) hold at all $x\in M$, then 
\[ 
(dt)_{g}(\Ker(dt)_{g}^{\perp})= \Im(\rho_{y}) .
\]}
Clearly, the left hand side contains the right hand side, and we will compare the
dimensions of the two spaces. Recalling that $\dim(M) = \dim(G) - \dim(\Ker(dt)_g)$ and using \eqref{eq:presymport},
we find that the dimension of the left hand side is
\begin{eqnarray*} 
 \dim(\Ker(dt)_{g}^{\perp})- \dim(\Ker(dt)_{g}^{\perp}\cap \Ker(dt)_{g})& = & 
\dim(\Ker(dt)_{g}\cap \Ker(\omega_g)) \\
 & &- \dim(\Ker(dt)_{g}^{\perp}\cap \Ker(dt)_{g})\\
 & &+ \dim(M). 
\end{eqnarray*}
Note that, by Lemma \eqref{lemma}, part (iii) (see Remark \eqref{rem:st}), 
we can replace $g$ in the right hand side of the last formula by $x= s(g)$.
Hence the dimension we are interested in equals the one of $(dt)_{x}(\Ker(dt)_{x}^{\perp})$.
Using again that $\Ker(dt)_{x}^{\perp}= \Ker(ds)_x+ \Ker(\omega_x)\cap T_xM$, together with (iii), we obtain
$(dt)_{x}(\Ker(dt)_{x}^{\perp})= \Im(\rho_x)$. But $ \dim(\Im(\rho_x)) = \dim(\Im(\rho_{y}))$
because $x$ and $y$ are on the same orbit, and this concludes the proof of the claim.

{\it Claim 3: If (i)-(v) hold at all $x\in M$, then the following is a short exact sequence:
\begin{equation}
\label{seq-dim} 
0\rmap \Ker(dt)_{g}\cap \Ker(dt)_{g}^{\perp} \rmap 
\Ker(ds)_{g}^{\perp}\cap \Ker(dt)_{g}^{\perp} \stackrel{(dt)_{g}}\rmap \Ker(\omega_{y})\cap TM\rmap 0 .
\end{equation}}
The claim easily follows if one checks the surjectivity of the last map.
To see that, note that
if $u \in \Ker(\omega_{y})\cap TM$, then $u \in \Im(\rho)$ by (iii); so Claim $2$ implies that
$u = dt_g(v)$, for some $v \in \Ker(dt)_g^\perp$. But by Claim $1$,  $v \in \Ker(ds)_g^\perp$.
Hence the last map is surjective.

As a result of Claim 3, we get 
\begin{equation}\label{eq:claim3}
\dim(\Ker(ds)_{g}^{\perp}\cap \Ker(dt)_{g}^{\perp}) = 
\dim(\Ker(dt)_{g}\cap \Ker(dt)_{g}^{\perp}) + \dim(\Ker(\omega_{y})\cap TM).
\end{equation}
Let us recall again that
$$
\dim(\Ker(dt)_{g}^{\perp}\cap \Ker(dt)_{g})= \dim(\Ker(dt)_{x}^{\perp}\cap \Ker(dt)_{x}).
$$
Using the short exact sequence 
$$
0\rmap \Ker(dt)_{x}^{\perp}\cap \Ker(dt)_{x}\rmap \Ker(dt)_{x}^{\perp}\stackrel{(dt)_{x}}{\rmap} \Im(\rho_{x})\rmap 0,
$$
and (\ref{eq:presymport}) to compute $\dim( \Ker(dt)_{x}^{\perp})$, we find
\begin{equation}
\label{dim1}  
\dim(\Ker(dt)_{g}^{\perp}\cap \Ker(dt)_{g})= \dim(M)+ \dim(\Ker(dt)_{x}\cap \Ker(\omega_x))- \dim(\Im(\rho_{x})) .
\end{equation}
Let us now compute the dimension of $\Ker(ds)_{g}^{\perp}\cap \Ker(dt)_{g}^{\perp}= (\Ker(ds)_{g}+ \Ker(dt)_{g})^{\perp}$.
Using (\ref{eq:presymport}) and replacing
$\dim(\Ker(ds)_{g}+ \Ker(dt)_{g})$ by $\dim(\Im(\rho_x)+ \dim(G)- \dim(M)$, which is possible due to the
exact sequence
\[ 
0\rmap \Ker(ds)_{g}\rmap \Ker(ds)_{g}+ \Ker(dt)_{g} \stackrel{(ds)_g}{\rmap} \Im(\rho_{x}) \rmap 0,
\]
we find 
\begin{eqnarray}
\label{dim2} 
\dim(\Ker(ds)_{g}^{\perp}\cap \Ker(dt)_{g}^{\perp}) & = & \dim(M)- \dim(\Im(\rho_{x})) +       \\
& & \dim( (\Ker(ds)_{g}+ \Ker(dt)_{g})\cap \Ker(\omega_{g})). \nonumber
\end{eqnarray}
Plugging  (\ref{dim1}) and (\ref{dim2}) into \eqref{eq:claim3}, we find
\begin{eqnarray*} 
 \dim( (\Ker(ds)_{g}+ \Ker(dt)_{g})\cap \Ker(\omega_{g}))  
  &=& \dim(\Ker(dt)_{x}\cap \Ker(\omega_x)) +\\
  &&  \dim(\Ker(\omega_{y})\cap TM)\\
 &=&  \frac{1}{2}(\dim(\Ker(\omega_x))+ \dim(\Ker(\omega_y))),
\end{eqnarray*}
where the last identity follows from the last two formulas in Lemma \ref{lemma2}.  
In particular,
\[ 
\frac{1}{2}(\dim(\Ker(\omega_x))+ \dim(\Ker(\omega_y))- \dim(\Ker(\omega_g))\leq 0 .
\]
On the other hand, we have seen in the proof of Lemma \ref{lemma-Dirac} that this number is the codimension
of $\Ker(dt)_{g}+ \Ker(\omega_{g})$ in $\Ker(ds)_{g}^{\perp}$, hence positive. 
So it must be zero, and this completes the proof.
\end{proof}

Note that, if $L_{\theta}$ is smooth, then it is automatically a $\phi$-twisted Dirac structure
since the 2-forms $\theta_{S}$ satisfy $d\theta_{S}= -\phi|_{S}$.
In particular, we obtain

\begin{corollary}
\label{proof-theorem} 
Let $\omega$ be a relatively $\phi$-closed, multiplicative 2-form. 
Then the following are equivalent:
\begin{enumerate}
\item[(i)] the bundle $L_M$ is smooth and $\Ker(\omega)\cap TM\subset \Im(\rho)$ 
(or any of the equivalent conditions in Lemma \ref{lem-orbit});
\item[(ii)] there is a $\phi$-twisted Dirac structure $L$ on $M$ so that $t: G\rmap M$ is a Dirac map,
and the presymplectic leaves of $L$ coincide with the orbits of $G$.
\end{enumerate}
If these conditions hold, then $L$ coincides with 
\[ 
L_{M}= \{ (\rho(\alpha)+ u, \rho_{\omega}^*(\alpha)): \alpha\in A, u\in (\Im(\rho^*_{\omega}))^{\circ} \},
\]
and $L_{M}$ is a $\phi$-twisted Dirac structure with presymplectic leaves $(S, \theta_{S})$ and 
kernel $\Ker(\omega)\cap TM= (\Im(\rho^*_{\omega}))^{\circ}$.
\end{corollary}

Another interesting consequence is the following:

\begin{proposition/definition} \label{defprop:oversymp}
Given a relatively $\phi$-closed, multiplicative 
2-form $\omega$ on $G$, the following are equivalent:
\begin{enumerate}
\item[(i)] there exists a $\phi$-twisted Poisson structure on $M$ so that $t: G\rmap M$ is a Dirac map;
\item[(ii)] $\dim(\Ker(\omega_x))= \dim(G)- 2 \dim(M)$ for all $x\in M$;
\item[(iii)] $\Ker(dt_{g})^{\perp}= \Ker(ds_{g})$; 
\item[(iv)] $\Ker(\omega_x)\subset \Ker(ds)_x\cap \Ker(dt)_x$ for all $x\in M$
(or, equivalently, $\Ker(\rho_{\omega}^*)\subset \Ker(\rho)$).
\end{enumerate}
In this case we say that $(G, \omega)$ is a {\bf twisted over-symplectic} groupoid. 
\end{proposition/definition}

\begin{proof}
Condition (i) is equivalent to $L_{M}$ being smooth and having zero
kernel, 
i.e. $\Ker(\omega)\cap TM= \{0\}$.
Hence, by \eqref{eq:keromega},
$(\Ker(dt)_{x})^{\perp}= \Ker(ds)_{x}$, for all $x\in M$. 
A simple dimension counting,  using \eqref{eq:presymport} and 
$\Ker(\omega_x) = \Ker(\omega_x)\cap \Ker(dt)_x$, directly shows (ii).
Using once again \eqref{eq:presymport} to compute $(\dim(dt)_g)^\perp$, and recalling that
$\dim(\Ker(\omega_g)\cap \Ker(dt)_g) = \dim(\Ker(\omega_x)\cap \Ker(dt)_x)$, for $x=s(g)$,
another dimension counting shows that (ii) implies (iii).
Note that (iii) implies that $\Ker(\omega_x) \subset \Ker(ds)_x$, and since (iii) also holds
for $t$ and $s$ interchanged, we get $\Ker(\omega_x) \subset \Ker(ds)_x\cap \Ker(dt)_x$
and (iv) follows. Finally, if (iv) holds, then (\ref{L-formula}) implies
that $(\rho, \rho_{\omega}^*): A\rmap L_{M}$ is surjective; thus $L_M$ is smooth.
Since $\Ker(\omega_x)\cap T_xM = \Ker(ds)_x\cap \Ker(dt)_x\cap T_xM =\{0\}$, 
$\omega$ is of Dirac type and the Dirac structure induced on $M$ is Poisson.
\end{proof}

\begin{definition}
\label{robust}
Let $G$ be a Lie groupoid over $M$, and let $\phi$ be a closed
3-form on $M$. 
A multiplicative 2-form $\omega$ on $G$ is 
called {\bf robust} if the isotropy bundle $\mathfrak{g}(\omega)$, defined by 
\[ 
\mathfrak{g}_{x}(\omega)= \Ker(\omega_{x})\cap \Ker(ds)_x\cap \Ker(dt)_x, \ \ \ x\in M, 
\]
has constant rank equal to $\dim(G)- 2\dim(M)$. If $\omega$ is also relatively $\phi$-closed,
we say that $(G, \omega)$ is a {\bf $\phi$-twisted over-presymplectic groupoid}.
\end{definition}

We will explain this ``over'' terminology in Remark \ref{Dirac-form} below.

A {\bf $\phi$-twisted presymplectic groupoid} (Def.~\ref{def:twistpresymp}) is a $\phi$-twisted 
over-presymplectic groupoid $(G, \omega)$ with  $\dim(G)= 2\dim(M)$. If $\omega$ is nondegenerate,
then $(G,\omega)$ is called a {\bf $\phi$-twisted symplectic groupoid}. (This coincides with
the notion of {\bf quasi-symplectic groupoid} of Cattaneo and Xu \cite{CaXu}.)

Note that, if $(G, \omega)$ is a $\phi$-twisted over-presymplectic groupoid, 
then $\mathfrak{g}(\omega)$ is a {\it smooth} bundle of Lie algebras.

%%%%%%%%%%%%%%%%%%%%%%%%%%%%%%%%%%%%%%%%%%%%%%%%%%%%%%%%%%%%%%%%%%%%%%%%%%%%%%%%%%%%%%%%%%%%%%%%%%%%%%%%%%%%%
% \begin{lemma}
% \label{lem:inducedirac1} 
% Let $\omega\in \Omega^2(G)$ be a multiplicative 2-form. The following  
% statements are equivalent:
% \begin{enumerate}
% \item[(i)] $\omega$ is robust;
% \item[(ii)] $(\rho, \rho^*_{\omega}): A\rmap L$ is surjective;
% \item[(iii)] $\Ker(dt)^{\perp}= \Ker(ds)+ \Ker(\omega)\cap \Ker(dt)$ at all $g \in G$;
% \item[(iv)] $\Ker(dt)^{\perp}= \Ker(ds)+ \Ker(\omega)\cap \Ker(dt)$ at all $x \in M$;
% \item[(v)] $\Ker(\omega)= \Ker(\omega)\cap \Ker(ds)+ \Ker(\omega)\cap \Ker(dt)$ at all $g \in G$;

% \item[(vi)] $\Ker(\omega)= \Ker(\omega)\cap \Ker(ds)+ \Ker(\omega)\cap \Ker(dt)$ at all $x \in M$;
% \item[(vii)] $\rho: \Ker(\omega)\cap \Ker(ds)\rmap \Ker(\omega)\cap TM$ is surjective at all $x \in M$.
% \end{enumerate} 
% In this case moreover, $\omega$ is of Dirac type, and $L_{M}$ will be a $\phi$-twisted Dirac structure on $M$.
% \end{lemma}

\begin{lemma}
\label{lem:inducedirac1} 
Let $\omega$ be a multiplicative 2-form on $G$. The following  
statements are equivalent:
\begin{enumerate}
\item[(i)] $\omega$ is robust;
\item[(ii)] $(\rho, \rho^*_{\omega}): A\rmap L_M$ is surjective;
\item[(iii)] $\Ker(dt)_{g}^{\perp}= \Ker(ds)_{g}+ \Ker(\omega_{g})\cap \Ker(dt)_{g}$;
\item[(iv)] $(dt)_{g}: \Ker(\omega_{g})\cap \Ker(ds)_{g}\rmap \Ker(\omega_{t(g)})\cap TM$ is
onto;
\item[(v)] $\Ker(\omega_{g})= \Ker(\omega_{g})\cap \Ker(ds)_{g}+ \Ker(\omega_{g})\cap \Ker(dt)_{g}$.
\end{enumerate} 
Note that (iii), (iv) and (v) are required to hold for all $g\in G$ or, equivalently, 
for all $g= x\in M$. 
In this case, $\omega$ is of Dirac type and $L_{M}$ is a $\phi$-twisted Dirac structure on $M$.
\end{lemma}

\begin{proof}
The equivalence of (i) and (ii) is immediate by a dimension counting.  
Using (\ref{eq:presymport}) and noticing that (iii) holds if and only if the dimensions
of the right and left hand sides coincide (due to the inclusion in
Lemma \ref{lemma}, part (ii)), we conclude that (iii) is equivalent to
\begin{equation}
\label{eq:dimcount1}
\dim(\Ker(ds)_{g})- \dim(M)= \dim(\mathfrak{g}_{g}(\omega)) 
\end{equation}
Also, $(\rho,\rho^*_{\omega})$ is surjective if and only if the dimension of its image 
(which equals $\dim(\Ker(ds)) - \dim(\Ker(\omega)\cap \Ker(ds)\cap \Ker(dt))$)
coincides with the dimension of $L_M$ (which equals $\dim(M)$), and this is
precisely (\ref{eq:dimcount1}) at points in $M$. On the other hand,
Lemma \ref{lemma}, part (iii), shows that (\ref{eq:dimcount1}) is equivalent
to the same relation at the point $t(g)\in M$. So we conclude that (i), (ii) and (iii) are equivalent to each other. A direct dimension counting, allied with Lemma \ref{lemma},
and then the last two formulas of Lemma \ref{lemma2}, shows that (iv) is also equivalent to (i).

Finally, a similar argument (i.e. dimension counting together with Lemma \ref{lemma}
and the last formula of Lemma \ref{lemma2}) shows that the equality in (v) is equivalent to
\begin{eqnarray}\label{eq:omegag}
\dim(\Ker\omega_{g}) & = & \frac{1}{2}(\dim(\Ker(\omega_{x})+ \dim(\Ker(\omega_{y}))+ 
\label{eq:dimcount2}\\ 
    \nonumber             &&  \dim(G)- 2\dim(M)- \dim(\mathfrak{g}_{x}(\omega))
\end{eqnarray}
at all $g \in G$ (where $x= s(g)$, $y= t(g)$). Evaluating this expression at  $g=x \in M$
immediately implies that $\omega$ is robust (i.e., (i)).
Conversely, if (iii) holds, then the equivalence of (ii) and (iv)
of Lemma \ref{lemma-Dirac} implies \eqref{eq:omegag}.

\end{proof}

\begin{corollary}\label{cor:multipform}
Let $G$ be a Lie groupoid over $M$, and let $\omega\in \Omega^2(G)$ 
be multiplicative. The following statements are equivalent:
\begin{enumerate}
\item[(i)] $\dim(G)= 2\dim(M)$, and $\omega$ is robust;
\item[(ii)] $(\rho, \rho^{*}_{\omega}): A\rmap L_M$ is an isomorphism;
\item[(iii)] $\Ker(dt)_{g}^{\perp}= \Ker(ds)_{g}\oplus \Ker(\omega_{g})\cap \Ker(dt)_{g}$;
\item[(iv)] $\Ker(\omega_{g})= \Ker(\omega_{g})\cap \Ker(ds)_{g}\oplus \Ker(\omega_{g})\cap \Ker(dt)_{g}$;
\item[(v)] $(dt)_{g}: \Ker(\omega_{g})\cap \Ker(ds)_{g}\rmap \Ker(\omega_{t(g)})\cap TM$ is an isomorphism,
\end{enumerate}
where (iii), (iv) and (v) are required to hold for all $g\in G$ or, equivalently, for all $g= x\in M$.

Moreover, if $\omega$ is relatively $\phi$-closed (i.e. $(G, \omega)$ is a $\phi$-twisted presymplectic groupoid),
then $L_M$ is a $\phi$-twisted Dirac structure on $M$, characterized
by the properties that $t$ is a Dirac map and that
$(\rho, \rho^*_{\omega})$ in (ii) above is an isomorphism of algebroids.
\end{corollary}

\begin{proof}
The equivalence of (i) and (ii) follows immediately from (\ref{eq:dimcount1}).
The rest is a direct consequence of Lemma \ref{lem:inducedirac1}. 
\end{proof}

Clearly, the last corollary proves Theorem \ref{theorem1}.

\begin{remark} \ \ \rm
\label{Dirac-form}

\begin{enumerate}
\item[(i)] For any $\phi$-twisted over-presymplectic groupoid $(G, \omega)$,
the isotropy bundle $\mathfrak{g}(\omega)$ 
is a smooth bundle of Lie algebras, and a
subbundle of the (possibly non-smooth) isotropy Lie algebra bundle of $G$,  
$\mathfrak{g}(G)= \Ker(\rho)$. Integrating $\mathfrak{g}(\omega)$ to a bundle of 
simply connected Lie groups $G(\omega)$, assuming that the quotient $G/G(\omega)$ is smooth, 
we see that $\omega$ reduces to a
multiplicative 2-form $\overline{\omega}$ on $G/G(\omega)$, and $(G/G(\omega), \overline{\omega})$
becomes a $\phi$-twisted presymplectic groupoid, over the same manifold $M$,
which induces the same $\phi$-twisted Dirac structure on $M$. 
% That is simply because $A(G/G(\omega))= A(G)/\mathfrak{g}(\omega)= A(G)/Ker(\rho, \rho^*)\cong L$.
(Of course, if $(G, \omega)$ is over-symplectic, then $(G/G(\omega), \overline{\omega})$
is a  symplectic groupoid.) 
Although this does not always work (i.e. $G/G(\omega)$ may be non-smooth),
this  explains our ``over'' terminology,
and it shows that $\mathfrak{g}(\omega)$ can be viewed as 
an obstruction to our final goal of obtaining a one-to-one correspondence between
$\phi$-twisted Dirac structures on $M$ and groupoids over $M$ equipped with a certain extra structure. 
Examples of over-presymplectic groupoids which cannot be reduced to presymplectic groupoids
will be given in Section \ref{subsec:over}.
\item[(ii)] Similarly, if $(G, \omega)$ is a $\phi$-twisted presymplectic groupoid over $M$ then,
again modulo global regularity issues, one can quotient it out by
$\Ker(\omega)$ to obtain a 
$\phi$-twisted syplectic groupoid over the $\phi$-twisted
Poisson manifold $M/\Ker(L)$. 
Details of this construction will be given elsewhere.
\end{enumerate}
\end{remark}

%%%%%%%%%%%%%%%%%%%%%%%%%%%%%%%%%%%%%%%
%%%%%%%%%%%%%%%%%%%%%%%%%%%%%%%%%%%%%%%
%%%%%%%%%%%%%%%%%%%%%%%%%%%%%%%%%%%%%%%
%%%%%%%%%%%%%%%%%%%%%%%%%%%%%%%%%%%%%%%
\section{Reconstructing multiplicative forms}
\label{Reconstructing multiplicative forms}
%%%%%%%%%%%%%%%%%%%%%%%%%%%%%%%%%%%%%%%
%%%%%%%%%%%%%%%%%%%%%%%%%%%%%%%%%%%%%%%
%%%%%%%%%%%%%%%%%%%%%%%%%%%%%%%%%%%%%%%
%%%%%%%%%%%%%%%%%%%%%%%%%%%%%%%%%%%%%%%

In this section we explain how multiplicative forms can be reconstructed from their infinitesimal counterpart.
In particular, we will complete the proof of Theorem \ref{theorem2}. 

Recall that
(cf. Proposition \ref{rho}) for $\omega\in \Omega^2(G)$ multiplicative and 
relatively $\phi$-closed, the associated bundle map $\rho_{\omega}^{*}$ satisfies
\begin{equation} \label{cond1}
\SP{\rho^*_{\omega}(\beta),\rho(\alpha)}= - \SP{\rho^*_{\omega}(\alpha),\rho(\beta)}
\, \mbox{ for\ all }\, \alpha, \beta \in \Gamma(A);
\end{equation}
\begin{equation}\label{cond2}
\rho^{*}_{\omega}([\alpha, \beta])= \mathcal{L}_{\alpha}(\rho^{*}_{\omega}(\beta)) - 
\mathcal{L}_{\beta}(\rho^{*}_{\omega}(\alpha))- 
d\SP{\rho^*_{\omega}(\beta),\rho(\alpha)}+ 
i_{\rho(\alpha)\wedge \rho(\beta)}(\phi).
\end{equation}
We will use the notation
\[ 
(d_{A}\rho^*_{\omega})(\alpha, \beta) := 
\rho^{*}_{\omega}([\alpha, \beta])- \mathcal{L}_{\alpha}(\rho^{*}_{\omega}(\beta)) + 
\mathcal{L}_{\beta}(\rho^{*}_{\omega}(\alpha))+ d\SP{\rho^*_{\omega}(\beta),\rho(\alpha)}.
\] 

The rest of this section is devoted to the proof of the next result (see also \cite{CaXu}).
\begin{theorem}
\label{rec1}
If $G$ is an $s$-simply connected Lie groupoid, and $\phi\in \Omega^3(M)$ is closed,
then the correspondence $\omega\mapsto \rho_{\omega}^{*}$ induces a bijection
between the space of relatively 
$\phi$-closed multiplicative 2-forms on $G$, and bundle maps $\rho^{*}: A\rmap T^*M$
satisfying conditions (\ref{cond1}) and (\ref{cond2}) above. 
\end{theorem}

By Prop.~\ref{rho}, part (iii), we
only have to check that the map $\omega\mapsto \rho^{*}_{\omega}$ is surjective, i.e., 
given $\rho^*: A\rmap T^*M$ satisfying (\ref{cond1}) and (\ref{cond2}), we must 
produce a relatively $\phi$-closed multiplicative 2-form $\omega$ on $G$. 
By Corollary \ref{cor:multipform}, if  $\Ker(\rho^*)\cap \Ker(\rho) = \{0\}$, then
the resulting 2-form makes $G$ into a twisted presymplectic groupoid;
note that if $L$ is a Dirac structure on $M$, then Theorem \ref{rec1}, applied
to $\rho^*= pr_2:L\rmap T^*M$, together with Corollary \ref{cor:multipform} imply
the correspondence in Theorem \ref{theorem2}.

Let us recall \cite{CrFe} how $G$ can be reconstructed from $A$. 
Let $I=[0,1]$ and $\proj:A \rmap M$ be the natural projection (we will also denote other bundle projections
by $\proj$ whenever the context is clear).
Consider the Banach manifold 
$\tilde{P}(A)$ 
consisting of paths $a:I\rmap A$ 
of class $C^1$, whose base path $\gamma= \proj\circ a:I\rmap M$ is of class $C^2$,
and the submanifold $P(A)$ defined by the equation $\rho\circ a= \frac{d}{dt}\gamma$ (i.e. $a$ is an $A$-path).
The manifold $P(A)$ comes endowed with an infinitesimal action of the infinite dimensional Lie algebra $\mathfrak{g}$
consisting of time dependent sections $\eta_{t}$ ($t\in [0, 1]$)\footnote{In this section, $t$ will denote a 
``time'' parameter.} of $A$, with $\eta_0= \eta_1= 0$. To define 
the Lie algebra map
\[ 
\mathfrak{g}\ni \eta\mapsto X_{\eta}\in \mathcal{X}(P(A)) 
\]
describing the action, it will be more convenient to introduce the flows of the vector fields
$X_{\eta}$. One advantage of this approach is that $X_{\eta}$ will be defined on the entire $\tilde{P}(A)$.
Given $a_0\in P(A)$, we construct  the flow 
$a_{\epsilon}= \Phi_{X_{\eta}}^{\epsilon}(a_0)$ in such a way that $a_{\epsilon}$
are paths above $\gamma_{\epsilon}(t)= \Phi_{\rho(\eta_t)}^{\epsilon}\gamma_{0}(t)$, where $\gamma_{0}$ is the
base path of $a_0$, and $\Phi_{\rho(\eta_t)}^{\epsilon}$ is the flow of the vector field $\rho(\eta_t)$.
We choose a time dependent section $\xi_{0}$ of $A$ with $\xi_{0}(t, \gamma_0(t))= a_0(t)$, 
and  consider the $(\epsilon, t)$-dependent section of $A$,
$\xi= \xi(\epsilon, t)$, solution of
\begin{equation}
\label{def-xi} 
\frac{d\xi}{d\epsilon}- \frac{d\eta}{dt}= [\xi, \eta], \; \; \xi(0, t)= \xi_0(t).
\end{equation}
Then $a_{\epsilon}(t)= \xi_{\epsilon}(t, \gamma_{\epsilon}(t))$.
This defines the desired vector fields $X_{\eta}$, the action of $\mathfrak{g}$, and the foliation
on $P(A)$. It is clear from the definition that
\begin{equation} 
\label{action}
(d\proj)(X_{\eta})= \rho(\eta)\circ \proj \;\; \mbox{ and } \;
\proj\circ \Phi_{X_{\eta}}^{\epsilon}(a)= \Phi_{\rho(\eta)}^{\epsilon}\circ \proj(a) .
\end{equation} 
Now, $G(A)= P(A)/\sim$ is a topological groupoid for any $A$: the source (resp. target) map
is obtained by taking the  starting (resp. ending) 
point of the base paths, and the multiplication is defined by concatenation of paths.
Moreover, $G$ must be isomorphic to $G(A)$.

We will construct forms on $G(A)$ by constructing forms on $P(A)$ which are basic with respect to 
the  action (i.e. $\Lie_{X_{\eta}}\omega=0$ and $i_{X_{\eta}}\omega= 0$ for all $\eta\in \mathfrak{g}$).
First, the 3-form $\phi$ on $M$ induces a 2-form
on $TM$ that, at $X\in T_xM$, is $i_{X}(\proj_{M}^{*}\phi_{x})$, where $\proj_{M}:TM\rmap M$ is the projection. 
We pull it back by $\rho$ to $A$ and lift it to $\tilde{P}(A)$ to get a 2-form $\omega_{\phi}$
on $\tilde{P}(A)$:
\[ 
\omega_{\phi, a}(V, W)= \int_{0}^{1} \phi(\rho(a(t)), (d\proj)_{a(t)}(V(a(t)), (d\proj)_{a(t)}(W(a(t)))) dt.
\] 
In order to produce basic forms, we must understand
the behavior of $\omega_{\phi, a}$ when we take derivatives along, or interior products by, 
vector fields of the form $X_{\eta}$. 

\begin{lemma}
\label{lema1}
For any closed 3-form $\phi$ on $M$,
\[ 
\omega_{\phi}(X_{\eta, a}, X_{a})= \int_{0}^{1} \phi(\rho a(t), \rho\eta(t, \gamma(t)), (d\proj)_{a(t)}(X_{a}(t))) dt,
\]
\[ 
d\omega_{\phi}= \t^{*}\phi- \s^{*}\phi, 
\]
where $\eta\in \mathfrak{g}$, $X_{a}$ is a vector tangent to $\tilde{P}(A)$ at $a \in P(A)$, 
and $\s, \t: \tilde{P}(A)\rmap M$ take the start/end points of the base path.
\end{lemma}

\begin{proof}
The first formula is immediate from the definition of $\omega_{\phi}$
 and the 
 first equation in (\ref{action}), 
while the second formula follows from  Stokes' Theorem.
\end{proof}

Let us now consider the Liouville one-form $\sigma^{c}$ on $T^*M$, and the associated
canonical symplectic form $\omega^c$.
Recall that, on $T^*M$, 
\[ 
\sigma^{c}_{\xi_x}(X_{\xi_x})= 
\SP{\xi_x, (d\proj)_{\xi_{x}}(X_{\xi_x})  } ,\ \ (\xi_x\in T_{x}^{*}M, X_x\in T_{\xi_x}T^*M)
\]
and $\omega^{c}= -d\sigma^{c}$. Using $\rho^*: A\rmap T^*M$, we pull
$\sigma^{c}$ and $\omega^{c}$  back to $A$, and we denote by $\tilde{\sigma}$
and $\tilde{\omega}$ the resulting forms on $\tilde{P}(A)$. Hence $\tilde{\omega}= - d\tilde{\sigma}$,
and 
\[ 
\tilde{\sigma}_{a}(X_a)= \int_{0}^{1} \SP{ \rho^*(a(t)), (dp)_{a(t)}(X_a(t))} dt .
\]

\begin{lemma}
\label{lema2} 
If $\rho^*$ satisfies (\ref{cond1}), then, for any $A$-path $a$  and any 
vector field $X$ on $\tilde{P}(A)$, 
\[ 
i_{X_{\eta, a}}(\tilde{\sigma}) = -\int_{0}^{1} \SP{{\rho^{*}\eta(t, \gamma(t))},\rho(a(t))} dt,
\;\; \mbox{ and }
\]
\[ 
\mathcal{L}_{X_{\eta}}(\tilde{\sigma})(X_{a})= 
\int_{0}^{1} \SP{{(d_{A}\rho^*)(a(t),\eta(t, \gamma(t))}, (d\proj)(X_{a}(t))} dt - 
i_{X_a}d(\int_{0}^{1} \SP{ {\rho^{*}\eta(t, \gamma(t))}, \frac{d\gamma}{dt} } dt),
\]
where the last term is the differential of the function 
$a\mapsto \int_{0}^{1} \SP{\rho^{*}\eta(t, \gamma(t)),\frac{d\gamma}{dt}} dt$,
$\gamma$ is the base path of $a$ and $\eta \in \mathfrak{g}$.
\end{lemma}

\begin{proof}
For the first formula, we use the definition of $\tilde{\sigma}$:
\[ 
\tilde{\sigma}(X_{\eta, a})= \sigma^{c}_{\rho^*(a)}((d\rho^*)_{a}(X_{\eta, a}))= 
\int_{0}^{1} 
\SP{\rho^*(a), (d\proj)(X_{\eta, a})} ,
\]
the first formula in (\ref{action}), and (\ref{cond1}).
To prove the second formula, we use the definition of $\tilde{\sigma}$ 
and the second formula in (\ref{action}) to rewrite the Lie derivative 
\[
\begin{split}
\mathcal{L}_{X_{\eta}}(\tilde{\sigma})(X_{a})& = \left.\frac{d}{d\epsilon}\right |_{\epsilon= 0} 
\tilde{\sigma}(d\Phi_{X_{\eta}}^{\epsilon})_{a}(X_a) =\\
 & =  \int_{0}^{1}
 \left.\frac{d}{d\epsilon}\right |_{\epsilon= 0} 
\SP{\rho^*\xi_{\epsilon}(t, \gamma_{\epsilon}(t)), 
   (d\Phi_{\rho(\eta_t)}^{\epsilon})_{\gamma(t)}(X^{\,'}(t))} dt
\end{split}
\]
where $\gamma_{\epsilon}(t)$ and $\xi_{\epsilon}$ are as in the construction above of the vector fields $X_{\eta}$,
and $X^{\,'}(t)= (d\proj)_{a(t)}(X_a(t))$. 
To compute this expression, we use a connection $\nabla$ on $M$. The  expression in the last integral
has the following two terms:
\begin{equation}
\label{term1}
\SP{\rho^*\xi_{\epsilon}(t, \gamma_{\epsilon}(t)), 
\partial_{\epsilon}(d\Phi_{\rho(\eta_t)}^{\epsilon})_{\gamma(t)}(X^{\,'}(t))    }
\end{equation}
\begin{equation}
\label{term2}
\SP{\partial_{\epsilon} \rho^*\xi_{\epsilon}(t, \gamma_{\epsilon}(t)),X^{\,'}(t)} 
\end{equation}
where $\partial_{\epsilon}$ is the derivation of paths in $TM$ and $T^*M$ induced by the connection. 
On the other hand, for any
vector fields $V$ and $W$ on $M$, 
$\partial_{\epsilon}(d\Phi_{W}^{\epsilon})_{x}(V_{x})= \overline{\nabla}_{V_x}(W)$,
where $\overline{\nabla}_{V}(W)= \nabla_{W}(V)+ [V, W]$ is the conjugated connection 
(this is a simple check in local coordinates). Hence (\ref{term1}) equals to
\[ 
\SP{\rho^*(\xi), [X^{\,'}, \rho(\eta)]+ \nabla_{\rho(\eta)}(X^{\,'})  }_{\gamma(t)}= 
\SP{\mathcal{L}_{\rho(\eta)}(\rho^*(\xi))- \nabla_{\rho(\eta)}(\rho^*(\xi)),X^{\,'}(t)}_{\gamma(t)},
\]
where, for a moment, we have made $X^{\,'}$ into a vector field extending $X^{\,'}(t)$ (for each fixed $t$). 
On the other hand, (\ref{term2}) equals to
\[ 
\SP{\nabla_{\frac{d\gamma}{dt}}(\rho^*(\xi))+ \frac{d\rho^*(\xi_{\epsilon})}{d\epsilon},X^{\,'}(t)}
= \SP{\nabla_{\rho(\eta)}(\rho^*(\xi))+ \rho^*([\xi, \eta])+ 
\frac{d\rho^*(\eta_{t})}{dt},X^{\,'}(t)},  
\]
(at the point $\gamma(t)$), where we have used the defining equation (\ref{def-xi}) for $\xi$. 
Adding the two expressions we obtained for (\ref{term1}) and (\ref{term2}) (at $\epsilon= 0$), we get 
\[
\begin{split}
\mathcal{L}_{X_{\eta}}(\tilde{\sigma})(X_{a}) & = \int_{0}^{1} 
\SP{\mathcal{L}_{\rho(\eta)}(\rho^*(\xi_0))+ \rho^*([\xi_0, \eta])+ \frac{d\rho^*(\eta_{t})}{dt},
X^{\,'}(t)}_{\gamma(t)}dt\\
 & = \int_{0}^{1} \SP{(d_A\rho^*)(\xi_0, \eta)+ i_{\rho(\xi_0)}(d\rho^*(\eta)) 
+ \frac{d\rho^*(\eta_{t})}{dt},X^{\,'}(t)} dt
\end{split}
\]
where, in the last equality, we used the definition of $d_{A}(\rho^*)$. 
At this point, the computation is transfered to $M$, since the expression 
\[ 
\mathcal{L}_{X_{\eta}}(\tilde{\sigma})(X_{a})- \int_{0}^{1} 
\SP{(d_A\rho^*)(a(t), \eta(t, \gamma(t))),(d\proj)_{a(t)}(X_{a}(t))} dt 
\]
equals to
\begin{equation}
\label{onM} 
\int_{0}^{1} \SP{i_{\frac{d\gamma}{dt}}(du^t) + \frac{du^{t}}{dt},X^{\,'}(t)}_{\gamma(t)} dt
\end{equation}
where $u^{t}= \rho^*(\eta^{t})$. To finish the proof, we will use the next lemma. 
\end{proof}

\begin{lemma}
\label{lema3}
For any path $\gamma$ on $M$, any path  $X^{\,'}:I \rmap TM$ above $\gamma$, and any time-dependent 
$1$-form $u^{t}$ on $M$, we have
\begin{equation*}
\label{onM2} 
\mathcal{L}_{X^{\,'}}(\int_{0}^{1} \SP{u(t, \gamma(t)), \frac{d\gamma}{dt}} dt) + 
\int_{0}^{1} \SP{i_{\frac{d\gamma}{dt}}(du^t) + \frac{du^{t}}{dt},X^{\,'}(t)}_{\gamma(t)} dt
= \left.\SP{u(t, \gamma(t)),X^{\,'}(t)}\right|_{0}^{1}
\end{equation*}
\end{lemma}

(The function on which $\mathcal{L}_{X^{\,'}}$ acts is defined as in Lemma \ref{lema2}.)

\begin{proof} We may assume that there
is a vector field $Z$ such that $Z(\gamma(t))= X^{\,'}(t)$ (otherwise one just brakes $\gamma$ into smaller paths,
and note that the formula to be proven is additive with respect to concatenation of paths). 
We first compute the first integral 
using, as above, a connection $\nabla$ and the formula 
$\partial_{\epsilon}(d\Phi_{W}^{\epsilon})_{x}(V_{x})= \overline{\nabla}_{V_x}(W)$:
\[
\begin{split}
\left.\frac{d}{d\epsilon}\right |_{\epsilon= 0}
\int_{0}^{1} \SP{u(t, \Phi_{Z}^{\epsilon}(\gamma(t))),
\frac{d}{dt}(\Phi_{Z}^{\epsilon}(\gamma(t)))} dt = \\
=\int_{0}^{1} 
\SP{u^t,\partial_{\epsilon} (d\Phi_{Z}^{\epsilon})_{\gamma(t)}(\frac{d\gamma}{dt}) } dt
+ 
\SP{\partial_{\epsilon}(u(t, \Phi_{Z}^{\epsilon}(\gamma(t))),\frac{d\gamma}{dt}}=\\
= 
\int_{0}^{1}(\SP{u^t,\overline{\nabla}_{\frac{d\gamma}{dt}}(X^{\,'})}+ \SP{u^t,\frac{d\gamma}{dt}}) dt
\end{split}
\]
Now, it is easy to see that the sum of the term in the last integral with the 
term appearing in the second integral in the statement
is precisely 
\[ 
\SP{u,\nabla_{\frac{d\gamma}{dt}}(X^{\,'})}_{\gamma(t)}
+ 
\SP{\nabla_{\frac{d\gamma}{dt}}(u)+ \frac{du^t}{dt},X^{\,'}}_{\gamma(t)}
= \frac{d}{dt}\SP{ u(t, \gamma(t)),X^{\,'}(t)} 
\]
\end{proof}
Using Cartan's formula $\mathcal{L}_{X}= di_{X}+ i_{X}d$, 
the next result follows directly from Lemma \ref{lema1}.

\begin{lemma}
\label{lema4} 
If $\rho^*$ satisfies (\ref{cond1}), then, for any $A$-path $a$ and any 
vector field $X$ on $\tilde{P}(A)$, we have 
\[ 
i_{X_{\eta, a}}(\tilde{\omega})(X_{a}) = 
i_{X_a}d(\int_{0}^{1}\SP{\rho^{*}(\eta),\frac{d\gamma}{dt}-\rho( a)} dt)- 
\int_{0}^{1} \SP{(d_{A}\rho^*)(a, \eta),(d\proj)(X_{a})} dt .
\]
\end{lemma}

We can now complete the proof of  Theorem \ref{rec1}, i.e., reconstruct the 2-form $\omega$ out of $\rho^*$.
Let us assume that $\rho^*$ satisfies both conditions (\ref{cond1}) and (\ref{cond2}),
and put $\tilde{\omega}_{\phi}= \tilde{\omega}+ \omega_{\phi}$. Then the first 
equation in Lemma \ref{lema1} and the equation in Lemma \ref{lema4} show that 
$i_{X_{\eta}}(\tilde{\omega}_{\phi})= 0$
on $P(A)$. On the other hand, the second equation in Lemma \ref{lema1} and the fact that
$d\s$ and $d\t$ vanish on $X_{\eta}$'s (since $\eta$ vanishes at end-points), imply  that
$$
i_{X_{\eta}}d \tilde{\omega}_{\phi}= i_{X_{\eta}}(\t^*\phi- \s^*\phi)= 0.
$$ 
Hence $\tilde{\omega}_{\phi}$
is basic, and induces a 2-form $\omega_{0}$ on $G(A)$. The multiplicativity of $\omega_{0}$ follows from the
additivity of integration. To compute the associated $\rho^{*}_{\omega_{0}}(\alpha)(X)= \omega(\alpha, X)$, 
one has to look
at the identification of the Lie algebroid of $G(A)$ with $A$ (see \cite{CrFe}). After straightforward 
computations, we find that 
\[ 
\rho^{*}_{\omega_{0}}(\alpha)(X)= \omega_{can}((\rho(\alpha),\rho^*(\alpha)), (X,0)) .
\]
Here, $\omega_{can}$ is the linear version of the symplectic form, 
\[ 
\omega_{can}((X_1,\eta_1), (X_2,\eta_2))= \eta_2(X_1)- \eta_1(X_2) .
\]
We see that $\rho_{\omega_{0}}^{*}= -\rho^*$. On the other hand $d\omega_{0}= t^*\phi- s^*\phi$,
as it follows from the similar formula satisfied by $\omega_{\phi}$ (cf. Lemma \ref{lema1}).
Hence $\omega= -\omega_0$ will have the desired properties. 

The construction of $\tilde{\omega}_{\phi}$ is inspired by the one in \cite{CaXu}, which we recover when $\rho^*$ is an isomorphism.

%%%%%%%%%%%%%%%%%%%%%%%%%%%%%%%%%%%%%%%
%%%%%%%%%%%%%%%%%%%%%%%%%%%%%%%%%%%%%%%
%%%%%%%%%%%%%%%%%%%%%%%%%%%%%%%%%%%%%%%
\section{Examples}
\label{examples}
%%%%%%%%%%%%%%%%%%%%%%%%%%%%%%%%%%%%%%%
%%%%%%%%%%%%%%%%%%%%%%%%%%%%%%%%%%%%%%%
%%%%%%%%%%%%%%%%%%%%%%%%%%%%%%%%%%%%%%%
%%%%%%%%%%%%%%%%%%%%%%%%%%%%%%%%%%%%%%%
We discuss in this section some examples of multiplicative 2-forms, presymplectic groupoids
and their corresponding Dirac structures.

%%%%%%%%%%%%%%%%%%%%%%%%%%%%%%%%%%%%%%%
%%%%%%%%%%%%%%%%%%%%%%%%%%%%%%%%%%%%%%%
%%%%%%%%%%%%%%%%%%%%%%%%%%%%%%%%%%%%%%%
%%%%%%%%%%%%%%%%%%%%%%%%%%%%%%%%%%%%%%%
\subsection{Multiplicative 2-forms: first examples}
%%%%%%%%%%%%%%%%%%%%%%%%%%%%%%%%%%%%%%%
%%%%%%%%%%%%%%%%%%%%%%%%%%%%%%%%%%%%%%%
%%%%%%%%%%%%%%%%%%%%%%%%%%%%%%%%%%%%%%%
%%%%%%%%%%%%%%%%%%%%%%%%%%%%%%%%%%%%%%%
We now discuss some basic examples of multiplicative 2-forms on Lie groupoids.

\begin{example} ({\it Lie groups})\rm
 
If $H$ is a Lie group (so the base $M$ is just the one-point space, 
consisting of the identity in $H$),
then the zero form is the only multiplicative form on $H$. This follows from Lemma \ref{lemma}
(part (ii), or part (iv)).
\end{example}

\begin{example}({\it Lie groupoids integrating tangent bundles})\rm

Let $M$ be a $\phi$-twisted presymplectic manifold. Hence $\phi\in \Omega^3(M)$ is closed,
and $M$ is equipped with a 2-form $\omega_{M}$ with $d\omega_M+ \phi= 0$. 
Consider the pair groupoid $G= M\times M$ with the product $(x, y)\circ (y, z)= (x, z)$
(hence $s= pr_2$, $t= pr_1$). A simple computation shows that $M \times M$ equipped with
the 2-form $\omega=pr_1^*\omega_M - pr^*_2\omega_M$ is a $\phi$-twisted presymplectic groupoid,
and that the $\phi$-twisted Dirac structure induced on $M$ (identified with the diagonal in $M \times M$) 
is just the one associated with $\omega_M$. As usual, one obtains the
$s$-simply connected
$\phi$-twisted presymplectic groupoid corresponding to $\omega_M$ by pulling
 $\omega \in \Omega^2(M\times M)$ 
back to $\Pi(M)$, the fundamental groupoid of $M$, using the natural covering map
$\Pi(M) \rmap M\times M$ (which is also a groupoid morphism).
\end{example}

\begin{example}
\label{Pull-backs}
({\it Pull-backs})\rm 

Let $L$ be a  $\phi$-twisted Dirac structure on $M$, 
and let $(G(L),\omega_L)$ be a presymplectic groupoid integrating it.
If $f: P\rmap M$ is a submersion (we will see that weaker conditions
are possible), we can form the {\it pull-back}
groupoid $f^*G(L):= P\times_{M}G(L)\times_{M}P$ consisting of
triples $(p, g, q)$ with $g: f(p)\leftarrow f(q)$, $s= pr_3$, $t= pr_1$,
and $(p, g, q)\cdot (q, h, r)= (p, gh, r)$. It is simple to check that
$\dim(f^*G(L))= 2 \dim(P)$, and that the form $pr_2^*\omega_L \in \Omega^2(f^*G(L))$
is multiplicative, robust and relatively $f^*\phi$-closed. So
$(f^*G(L),pr_2^*\omega_L)$ is a $f^*\phi$-twisted presymplectic groupoid over
$P$. Infinitesimally, it corresponds to the {\it pull-back} Dirac structure (see e.g. \cite{BuRa})
\[ 
f^*L= \{ (X, f^*(\xi)): ((df)(X), \xi)\in L \}.
\]
We remark that the construction of the pull-backs $f^*L$ and $f^*G(L)$
is also possible in situations where $f$ is not a submersion.
In this case, $(f^*G(L), pr_{2}^{*}\omega_{L})$ is often just 
a (twisted) {\it over}-presymplectic groupoid, but not presymplectic. Such examples
arise for instance when one considers inclusions of submanifolds (see Example 
\ref{Dirac submanifolds of Poisson manifolds}). 
\end{example}

\begin{example}
\label{non-dirac} ({\it Multiplicative 2-forms of non-Dirac type})\rm

There are closed multiplicative 2-forms which are not of Dirac type.
In order to provide an explicit example, we start with a general observation.
Let $G$ be a groupoid over $M$, and let $\theta$ be a closed 2-form
on $M$.

{\it Claim: If the multiplicative 2-form $\omega=t^*\theta
  -s^*\theta$ on $G$ is of Dirac type, then
\[ 
\Im(\rho_x)+ \Im(\rho_x)^{\perp_{\theta}} 
\]
has constant dimension along points $x$ in a fixed orbit of $G$. (Here ``$\perp_\theta$'' is the orthogonal
with respect to $\theta$.)}

Let us prove the claim, following the notation of Section \ref{sec:nondeg}.

If $\omega$ is of Dirac type, then the ranges of $t_{*}(L_{\omega, g})$
and $L_{M, t(g)}$ must coincide. Let us denote $s(g)= x$, $t(g)= y$.
Since $\omega=t^*\theta -s^*\theta \in \Omega^2(G)$, we can write 
$\rho^{*}_{\omega}(\alpha)= i_{\rho(\alpha)}(\theta)$. Now the second
formula in (\ref{kernels}) implies that 
$\Ker(\omega)\cap T_y M= \Im(\rho_{y})^{\perp_{\theta}}$ for all $y\in M$,
and, using (\ref{L-formula}), we see that
\[ 
\dim(\Ran(L_{M, y}))= \dim( \Im(\rho_y)+ \Im(\rho_y)^{\perp_{\theta}}) .
\]
On the other hand, it is easy to see from the definition of $t_{*}(L_{\omega, g})$ that
$\Ran(t_{*}(L_{\omega, g}))= (dt)_{g}(\Ker(dt_{g})^{\perp})$ which,
as shown in the proof of Claim 2, Lemma \ref{lem-orbit}, 
has dimension equal to $(dt)_{x}(\Ker(dt_{x})^{\perp})=(dt)_{x}(\Ker(ds_{x})+ \Ker(\omega_x))$. Hence
\[ 
\dim(\Ran(t_{*}(L_{\omega, g}))= \dim(\Im(\rho_x)+ \Ker(\omega_x)\cap T_xM)=
\dim(\Im(\rho_x)+ \Im(\rho_x)^{\perp_{\theta}}) ,
\]
and this proves the claim.

To find our example, let $v$ be a vector field on $M$, 
and let $\Phi_{v}$ denote its flow.
The domain $G(v) \subset \mathbb{R}\times M$ of $\Phi_{v}$ 
is a groupoid over $M$
(the source is the second projection, the target 
is $\Phi_v$, and the multiplication is defined by 
$(t_1, y)(t_2, x)= (t_1+ t_2, x)$) integrating
the action Lie algebroid defined by $v$ (the underlying vector bundle
is the trivial line bundle, the anchor is multiplication by $v$, and the bracket
is $[f, g]= fv(g)- v(f)g$). On such a groupoid, the image of $\rho$ is either zero or one-dimensional, so
$\Im(\rho_x)+ \Im(\rho_x)^{\perp_{\theta}}= \Im(\rho_x)^{\perp_{\theta}}$ 
for any closed $\theta \in \Omega^2(M)$, 
and this need not be constant along orbits of $v$. For instance, one can take
\[ 
M= \mathbb{R}^2, \;\; v= x\frac{\partial}{\partial y}- 
y\frac{\partial}{\partial x},\;\;  \theta= y dx dy .
\]
Then the circle $S^1$ is an integral curve, and
$\Im(\rho_x)^{\perp}$ is one-dimensional everywhere on $S^1$,
except for the points $(1, 0)$ and $(-1, 0)$. 
\end{example}

\subsection{Examples related to Poisson manifolds}
%%%%%%%%%%%%%%%%%%%%%%%%%%%%%%%%%%%%%%%
%%%%%%%%%%%%%%%%%%%%%%%%%%%%%%%%%%%%%%%
%%%%%%%%%%%%%%%%%%%%%%%%%%%%%%%%%%%%%%%
%%%%%%%%%%%%%%%%%%%%%%%%%%%%%%%%%%%%%%%

\begin{example}({\it Symplectic groupoids})\rm 

Let $(G,\omega)$ be a presymplectic groupoid over $M$, and let $L$ be the corresponding
Dirac structure on $M$. Recall that $L$ comes from a Poisson
structure if and only if 
$\Ker(L)= \Ker(\omega)\cap TM = \{0\}$. But this condition is equivalent to
$\Ker(\omega_x)=0$ for all $x \in M$ (by Corollary \ref{cor:multipform}, part (v), 
and Lemma \ref{lemma2}), which is in turn equivalent
to $\Ker(\omega_g)=0$ for all $g \in G$ (by Lemma \ref{lemma-Dirac}, (iv)).
Hence $L$ comes from a Poisson structure if and only if $\omega$ is nondegenerate. 
We see in this way that our main result, restricted to Poisson structures, recovers  
the well-known correspondence between Poisson manifolds and symplectic groupoids.
In the presence of a closed 3-form, we recover the twisted version 
of this correspondence,
which was conjectured in \cite{WeSe} and proved in \cite{CaXu}. 
\end{example}

\begin{example}({\it Gauge transformations of Poisson manifolds})\rm 

Following \cite{BuRa}, we now explain how to 
produce twisted presym\-plec\-tic group\-oids
out of symplectic groupoids through gauge 
transformations associated to
2-forms.

Let $M$ be a smooth manifold equipped with a 
$\phi$-twisted Poisson structure $\pi$. We denote the
corresponding $\phi$-twisted Dirac structure 
by $L_{\pi}= \mbox{graph}(\widetilde{\pi})$.
The gauge transformation of $L_{\pi}$ 
associated to a 2-form $B$ on $M$
is given by
$$
L_{\pi} \mapsto \tau_B(L_{\pi}) = 
\{(\widetilde{\pi}(\eta),\eta + \widetilde{B}(\widetilde{\pi}(\eta))),\;\;
\eta \in T^*M \}
$$
As explained in  \cite{WeSe}, 
$\tau_B(L_{\pi})$ is a $(\phi-dB)$-twisted Dirac structure which may fail to be Poisson.

Let $(G,\omega)$ be a $\phi$-twisted symplectic groupoid integrating $\pi$.
Since $L_{\pi}$ and $\tau_B(L_{\pi})$ have isomorphic Lie algebroids (see \cite{WeSe}),
$\tau_B(L_{\pi})$ is integrable as an algebroid to a groupoid isomorphic to $G$.
The 2-form
\begin{equation}\label{eq:gaugeomega}
\tau_B(\omega) := \omega + t^*B - s^*B \in \Omega^2(G)
\end{equation}
is easily seen to be multiplicative, robust and relatively $(\phi-dB)$-closed, and,
as remarked in \cite[Thm~2.16]{BuRa}, it induces the Dirac structure $\tau_B(L_{\pi})$
on $M$. So $(G,\tau_B(\omega))$ is the presymplectic groupoid corresponding to
$\tau_B(L_{\pi})$.

Our results also show that this is true more generally: If $L$ is a (twisted) Dirac 
structure on $M$ associated with a presymplectic groupoid $(G(L),\omega_L)$,
then $\rho^{*}_{\tau_B(\omega_L)}= \rho^*_{\omega}+ i_{\rho(\alpha)}B$; hence the image
of $(\rho, \rho^*_{\tau_B(\omega_L)})$ is $\tau_B(L)$, and $(G(L),\tau_B(\omega_L))$
is a presymplectic groupoid integrating $\tau_B(L)$. 
\end{example}

\begin{example}({\it Dirac submanifolds of Poisson manifolds})\rm
\label{Dirac submanifolds of Poisson manifolds}

In this example, we relate our results to those in \cite[Sec.~9]{CrFe2}.
We describe how certain submanifolds of 
Dirac manifolds carrying an induced Dirac structure 
(such submanifolds of Poisson manifolds were studied in  \cite{CrFe2,Xu01}) 
give rise to over-presymplectic groupoids,
 whose reduction (in the sense of Remark \ref{Dirac-form}, (i)) produce presymplectic groupoids of
the submanifolds. For simplicity, we restrict the discussion to the untwisted case.

Let $L_M$ be a Dirac structure on $M$, let $N \hookrightarrow M$ be a submanifold,
and suppose that the pull-back Dirac structure induced at each point by inclusion,
\[ 
L_{N}:= \{ (X, \xi|_{TN}): X\in T_xN, (X, \xi)\in L_M\}\subset TN\oplus T^*N, 
\]
is a smooth bundle; it is not difficult to check that $L_{N}$ defines a Dirac structure on $N$. 
In the particular case of $L_M$ coming from a Poisson structure $\pi_M$ on $M$, 
and $L_{N}$ coming from a Poisson structure $\pi_{N}$
on $N$, $N$ is called a {\bf Poisson-Dirac submanifold}\footnote{A Poisson submanifold
is a Poisson-Dirac submanifold for which the inclusion 
$(N,\pi_N) \hookrightarrow (M,\pi_M)$ is a Poisson map; this is equivalent to 
$\Im(\tilde{\pi}_N)=\Im(\tilde{\pi}_M)$.} of $(M, \pi_M)$ \cite[Sec.~9]{CrFe2}. 

Let us consider the vector bundle  
\[ 
\mathfrak{g}_{N}(M):= TN^{\circ} \cap(L_M\cap T^*M) = (\Ran(L_{M})+ TN)^{\circ} 
\]
over $N$, which is in fact a bundle of Lie algebras, and assume that it has constant rank. 
(Here $^\circ$ denotes the annihilator.)
In this case, the restriction of the Lie algebroid $L_{M}$ to $N$ \cite{HiMa} is well-defined
and determines a  Lie subalgebroid
$L_{N}(M)$ whose underlying vector bundle is
\[ 
\{ (X, \xi): X\in TN, (X, \xi)\in L_{M} \}.
\]
This Lie algebroid fits into the following exact sequence of Lie algebroids:
\[ 
0\rmap \mathfrak{g}_{N}(M)\rmap L_{N}(M)\rmap L_N \rmap 0.
\]
Let us assume that $L_M$ is integrable, and let $(G(L_M), \omega_{M})$
be the associated presymplectic groupoid. 
Then $L_{N}(M)$ is also integrable (as it sits inside $L_M$ as a Lie subalgebroid),
and the associated groupoid $G(L_{N}(M))$ is a subgroupoid
of $G(L_M)$. Moreover, the restriction $\omega_{N, M}$ of $\omega_{M}$ to $G(L_{N}(M))$
makes $G(L_{N}(M))$ into an over-presymplectic groupoid over $N$, corresponding to the
pull-back Dirac structure $L_N$. 

We observe, however, that the reduction
procedure of Remark \ref{Dirac-form}, (i), produces a smooth presymplectic groupoid  
if and only if $L_{N}$ is also 
integrable. In this case, the reduced groupoid will be precisely 
the pull-back (see Example \ref{Pull-backs}) 
presymplectic groupoid $(G(L_N), \omega_N)$ of $L_N$. 
Hence, the presymplectic groupoids 
of a Dirac structure and of a Dirac submanifold are related by reduction of an  
intermediary over-presymplectic groupoid, as illustrated below: 
\[ \newdir{ (}{{}*!/-5pt/@^{(}}
\xymatrix{ (G(L_N(M)), \omega_{N, M})\ar@{ (->}[r]\ar[dr] &(G(L_M), \omega_M)\\
& (G(L_N), \omega_N)}\]

In general, this quotient space may not
be a manifold, but it is the same as the 
Weinstein groupoid introduced in \cite{CrFe} (see also
\cite{CaFe}).

If $L_{M}$ comes from a Poisson structure $\pi_M$, the discussion above shows that
$N$ is a Poisson-Dirac submanifold of 
$(M, \pi_M)$ if and only if $G(L_{N}(M))$ is an over-symplectic groupoid, and the reduction
procedure just described  recovers Proposition~9.13 in \cite{CrFe2}. 
%Changing the point of interest, and given the Dirac manifold $N$, 
%this gives a possible strategy for describing
%the presymplectic groupoid $G(L_N)$ in terms of the more familiar symplectic groupoids.
%Of course, the first step is to embed $N$ into a Poisson manifold $(M, \pi)$ so
%that $L_{N}= L_{\pi}|_{N}$. 
\end{example}

%%%%%%%%%%%%%%%%%%%%%%%%%%%%%%%%%%%%%%%
%%%%%%%%%%%%%%%%%%%%%%%%%%%%%%%%%%%%%%%
%%%%%%%%%%%%%%%%%%%%%%%%%%%%%%%%%%%%%%%
%%%%%%%%%%%%%%%%%%%%%%%%%%%%%%%%%%%%%%%
\subsection{Over-presymplectic groupoids and singular presymplectic groupoids}\label{subsec:over}
%%%%%%%%%%%%%%%%%%%%%%%%%%%%%%%%%%%%%%%
%%%%%%%%%%%%%%%%%%%%%%%%%%%%%%%%%%%%%%%
%%%%%%%%%%%%%%%%%%%%%%%%%%%%%%%%%%%%%%%
%%%%%%%%%%%%%%%%%%%%%%%%%%%%%%%%%%%%%%%

In this subsection we discuss examples of over-presymplectic groupoids that
cannot be reduced to presymplectic groupoids, as already mentioned in Example 
\ref{Dirac submanifolds of Poisson manifolds}. In particular, we will provide 
concrete examples of non-integrable Poisson structures admitting over-symplectic 
groupoids inducing them.

Let us consider the following particular case of the construction in 
Example \ref{Dirac submanifolds of Poisson manifolds}.
Let $(M,\pi)$ be a Poisson manifold, with 
$s$-simply connected symplectic groupoid $(G,\omega)$.
Let $C$ be a Casimir function on $M$, and let $a$ be a regular value of $C$.
Then the level manifold $N = C^{-1}(a) \hookrightarrow M$ is a Poisson submanifold
of $M$. 
As we saw in the previous example, the restriction $G_{N,M}=t^*(N)=s^*(N)$ of $G$
to $N$ is an $s$-simply connected groupoid over $N$, and the pull-back
$\omega_{N,M}$ of $\omega$ to $G_{N,M}$ makes it into an over-symplectic groupoid
over $N$, inducing on $N$ its Poisson structure. Indeed, the kernel of $\omega_{N,M}$ is
spanned by the hamiltonian vector field of $t^*C=s^*C$, and since this
vector field projects to zero on $M$, it is easy to check that
condition (iv) of Lemma \ref{defprop:oversymp} is satisfied. To pass
to a symplectic groupoid for $N$, we simply form the quotient of
$G_{N,M}$ by the hamiltonian flow.  

\begin{example}(Over-symplectic groupoids of nonintegrable Poisson submanifolds) \rm
\label{Casimir}

To obtain an interesting class of examples of the above construction,
take  $G$ to be the cotangent bundle $T^*H \cong H\ltimes \hhstar$ 
of a simply-connected Lie group $H$,
$M=\hhstar$, and 
$$
C:\hhstar\to \Rr, \;\; u \mapsto \frac{1}{2}(u,u)_{\hhstar}
$$ 
the kinetic energy function  of a bi-invariant (possibly
indefinite) metric $(\cdot,\cdot)_{\hhstar}$ on $H$. 
The hamiltonian flow of $t^*(C)$ is then the geodesic
flow on $T^*H$.  If $N$ is the unit ``sphere'' $\hhstar_{1/2}=C^{-1}(1/2)$,
then the unit co``sphere'' bundle $G_{N,M} = (T^*H)_{1/2}$ is
an oversymplectic   groupoid over it.
The quotient $Q$ of $(T^*H)_{1/2}$
by the geodesic flow is then the canonical symplectic groupoid for
the Poisson manifold  $\hhstar_{1/2}$. The elements of $Q$ are
geodesics in $H$ considered as oriented submanifolds. One can view
them as cosets of an open subset of the connected oriented
one-dimensional subgroups of $H$. From this point of view, the
groupoid structure of $Q$ is just the one induced from the Baer
groupoid (see \cite{WeGeom}) consisting of all the cosets in $H$.

A specific case of the construction in the previous paragraph
recovers the ``pathology'' of symplectic groupoids discussed in
Section 6 of \cite{CaFe}.  We let $H$ be the product of a ``space
manifold'' $SU(2)=S^3$ with its usual riemannian metric and a
``time manifold'' $\Rr$ carrying the negative of its usual metric.
The unit ``sphere'' in $\hhstar$ for this Lorentz metric may then
be identified with the product $S^2\times \Rr$, with the Poisson
structure for which the $S^2$ slice over each $\tau\in \Rr$ is a
symplectic leaf with symplectic structure equal to $1+\tau^2$
times the standard symplectic structure.  The critical point of
$1+\tau^2$ at $\tau=0$ is responsible for a singularity in the
space of unit-speed (equivalently, space-like) geodesics.  Most of
these geodesics are diffeomorphic to $\R$, but the ones which are
perpendicular to the $\tau$ direction are circles.  As a result,
the 6-dimensional quotient groupoid, i.e. the space of geodesics
in $S^3\times\R$, is a singular fiber bundle over
$S^2\times S^2$---the smooth 4-dimensional manifold of oriented
geodesics\footnote{The projection $T^*H \to T^*S^3$ is equivariant
with respect to the geodesic flows.} in $S^3$.
The fiber, which is neither Hausdorff nor locally Euclidean, is
 the quotient
space $\R^2/(x, y) \sim (x, x+y)$  It can be
obtained from the standard cone in $\Rr^3$ by removing the vertex
and replacing it with a line, with a topology such that any
sequence of points on the cone converging to where the vertex used to
be now converges to every point on the line.
\end{example}

An alternative way of obtaining over-presymplectic groupoids (which may not be reducible)
inducing a given Dirac structure $L$ is by means of extensions of the Lie
algebroid of $L$ by 2-cocycles (see e.g. \cite{Mack}). 

Let $L$ be a Dirac structure on $M$, and let $u\in \Gamma(\Lambda^2L^*)$ be a 2-cocycle
of the Lie algebroid of $L$. Let us assume, for simplicity, that $L$ is not twisted.
Then there is an associated algebroid $L\ltimes_{u}\mathbb{R}= L\oplus \mathbb{R}$
which fits into an exact sequence 
\begin{equation}\label{eq:extension}
0\to \R \to L\ltimes_{u}\mathbb{R} \to L \to 0.
\end{equation}
The anchor of $L\ltimes_{u}\mathbb{R}$ is just $(X,a)\mapsto \rho(X)$, where $\rho$ is the anchor of $L$,
while the Lie bracket on $\Gamma(L)$ is
\[ 
[(X, a), (Y, b)]=([X, Y], \mathcal{L}_X b-\mathcal{L}_Y a + u(X, Y)).
\]
The interesting point of this construction is that $L\ltimes_{u}\mathbb{R}$ 
may be integrable even when
$L$ is not. 
Moreover, the associated groupoid $G(L\ltimes_{u}\mathbb{R})$ admits a canonical 
$1$-form $\sigma_u$ (whose construction is similar to the construction of the form $\sigma$
in Proposition \ref{ex:monodromy} (iii)), and $(G(L\ltimes_{u}\mathbb{R}), d\sigma_u)$ will be an
over-presymplectic groupoid over $M$ inducing $L$ on the base.
In the case of a Poisson manifold $(M, \pi)$ ($L= L_{\pi}$),
the groupoid corresponding to the extension of the Lie algebroid $T^*M$ by $u=\pi$, together
with the $1$-form $\sigma_u$, is an example of a {contact groupoid}. 

\begin{example}
\label{Contact groupoids}
({\it Contact groupoids})\rm 

Contact groupoids are groupoids associated with Jacobi manifolds, as we now
briefly explain (see \cite{crainiczhu} and references therein for details).

A {\bf Jacobi manifold} is a manifold equipped with a bivector $\Lambda$
and a vector field $E$ satisfying $[\Lambda, \Lambda]= 2\Lambda\wedge E$, and
$[\Lambda, E]= 0$. A particular example is given by Poisson manifolds $(M,\pi)$,
in which case $\Lambda = \pi$ and $E=0$; more generally, if $g \in C^\infty(M)$,
then $g \pi$ may fail to be Poisson, but $\Lambda=g\pi$ and $E=X_{g}$ define
a Jacobi structure on $M$.

Any Jacobi manifold $(M,\Lambda,E)$ determines a Lie algebroid structure on
$T^*M\oplus \mathbb{R}$, which in the case of a Poisson manifold $(M,\pi)$
is precisely $T^*M\ltimes_{\pi}\mathbb{R}$ defined in \eqref{eq:extension}.
Just as Poisson structures correspond to symplectic groupoids (or Dirac structures
correspond to presymplectic groupoids), Jacobi manifolds are associated with
{\bf contact groupoids}. 
These are groupoids $G$ endowed with a contact $1$-form $\sigma$
and a smooth function $f\in C^\infty(G)$ such that $\sigma$ is $f$-multiplicative:
\[ 
m^* \sigma= pr_2^* f \cdot pr_1^* \sigma + pr^*_2 \sigma.
\]
(Here $pr_j: G\times G \to G$ is the natural projection onto the $j$-th factor.)
When $M$ is a Poisson manifold, the function $f$ is constant equal to $1$, so
$\sigma$ is multiplicative. One can therefore
associate two groupoids to a Poisson
manifold $(M,\pi)$, its symplectic groupoid $(G_s,\omega)$ and
the contact groupoid $(G_c,\sigma)$ obtained by regarding $(M,\pi)$ 
as a Jacobi manifold. 
Since $\sigma$ is multiplicative,
so is the 2-form $d\sigma$, and one can check that $(G_c,d\sigma)$ is  an over-symplectic
groupoid over $M$ inducing  on $M$ its Poisson structure. As a result, one can see $(G_s,\omega)$
as a reduction of  $(G_c,d\sigma)$.

So, whenever a Poisson manifold $(M,\pi)$ is integrable as a Jacobi manifold,
it automatically admits an over-symplectic groupoid inducing $\pi$.
As explained in \cite{crainiczhu}, $(M,\pi)$ can be a nonintegrable Poisson manifold, 
and yet be integrable as a Jacobi manifold (and vice-versa).
For example, the co``sphere'' bundle $(T^*H)_{1/2}$ in Example \ref{Casimir} is
actually a contact groupoid with contact form given by the restriction of
the canonical 1-form on $T^*H$. For more concrete examples and a detailed
comparison of the obstructions to Poisson and Jacobi integrability, we refer
the reader to \cite{crainiczhu}.
\end{example}

\subsection{Lie group actions and equivariant cohomology}
\label{Lie group actions}

Let $H$ be a connected Lie group acting on a manifold $M$. 
We consider the {\it action groupoid} $H\ltimes M$ over $M$, with 
\[ 
s(g, x)= x, \;\;\;\;\; t(g, x)= g x, \;\;\; \forall (g, x)
\in H \times M,
\]
and multiplication of composable pairs given by
\[
m( (g_1, x_1), (g_2, x_2) )= (g_1 g_2, x_2).
\]
It was pointed out in \cite{bxz} that
 the space of twisted multiplicative 2-forms on
$H\ltimes M$ is closely related to the equivariant cohomology of $M$ in degree three. 
Our main results provide the following description of this
 relationship at the infinitesimal level. 

The Lie algebroid of $H\ltimes M$ (the {\it action Lie algebroid} $\mathfrak{h} \ltimes M$)
is the trivial bundle $\mathfrak{h}_M:= \mathfrak{h}\times M \to M$; the anchor
is defined by the infinitesimal action, and the bracket is uniquely determined by the
Leibniz rule and the Lie bracket on $\mathfrak{h}$. Following Theorem \ref{theorem3}, 
the infinitesimal counterpart of 
twisted multiplicative 2-forms on $H\ltimes M$ are pairs $(\rho^*, \phi)$, where $\phi$ is
a closed 3-form on $M$, and $\rho^*:\mathfrak{h}_M\rmap T^*M$ is a bundle map
satisfying conditions (\ref{cond1}), (\ref{cond2}). We denote by 
$\omega_{\rho^*, \phi} \in \Omega^2(H\ltimes M)$ the
relatively $\phi$-closed, multiplicative 2-form associated with $(\rho^*,\phi)$.

When $H$ is a compact Lie group,
the equivariant cohomology of $M$ can be computed by Cartan's complex of equivariant 
differential forms on $M$, denoted $\Omega_{H}^{*}(M)$ (see e.g. \cite{Vergne}). This complex 
consists of $H$-invariant, $\Omega^*(M)$-valued
polynomials on $\mathfrak{h}$, with degree twice the polynomial degree plus the form-degree:
\[ 
\Omega^{k}_{H}(M)= (\bigoplus_{2i+ j= k} S^{i}(\mathfrak{h}^*)\otimes \Omega^{j}(M))^H .
\]
The differential is $d_H:= d_1- d_2$ where $d_{1}(P)(v)= d(P(v))$, $d_2(P)(v)= i_{v}(P(v))$; 
if $\alpha$ is an $\Omega^*(M)$-valued
polynomial on $\mathfrak{h}$, invariance
means that 
\begin{equation}\label{eq:invariance}
g^*(\alpha(\Ad_g(v)))=\alpha(v),
\end{equation} 
for all $v \in \mathfrak{h}$ and $g \in H$.

Note that equivariantly closed 3-forms can be written as
\[ 
\rho^*+ \phi \in \Omega^{3}_{H}(M),
\] 
where $\phi\in \Omega^3(M)$ is closed, and $\rho^* \in \mathfrak{h}^*\otimes \Omega^1(M)$, 
which are both invariant (as in \eqref{eq:invariance})  and satisfy
\[ \left \{ \begin{array}{ll} 
          i_{v}(\rho^*(v)) = 0\\
          i_{v}(\phi)- d(\rho^*(v))= 0
           \end{array}
  \right. \]
for all $v\in \mathfrak{h}$. 
Infinitesimally, the invariance condition for $\rho^*$ reads
\[ 
\rho^*([v, w])= \Lie_v(\rho^*(w)) ,
\]
for all $v, w\in\mathfrak{h}$. Using this equation, one can easily check that 
conditions (\ref{cond1}), (\ref{cond2}) are satisfied, and hence
there is a corresponding multiplicative 2-form $\omega_{\rho^*, \phi}$. 
Assuming that $\rho^*+ \phi \in \Omega^{3}_{H}(M)$, 
we will now describe how  one can obtain a simple explicit formula for $\omega_{\rho^*, \phi}$
just using general properties of multiplicative forms.

Let $pr_g$ and $pr_x$ be the natural projections of
$H \times M$ onto $H$ and $M$, respectively, and let $\lcart, \rcart$ denote
the left and right-invariant Maurer-Cartan forms on $H$ (i.e., 
$\lcart_g(V)= (dL_{g^{-1}})_{g}(V), \; \rcart_g(V)= (dR_{g^{-1}})_{g}(V)$).

\begin{proposition}
\label{recover-AMM} 
Suppose that 
$\rho^*$ and $\phi$ satisfy conditions (\ref{cond1}), (\ref{cond2}), and let
$\omega= \omega_{\rho^*, \phi} \in \Omega^2(H\ltimes M)$  be the corresponding 2-form. 
The following are equivalent:
\begin{enumerate}
\item[(i)] $\rho^*+ \phi\in\Omega_{H}^{3}(M)$;
\item[(ii)] the restriction of $\omega$ to all slices $\{g\}\times M$ vanishes, for all $g\in H$;
\item[(iii)]  $\omega$ is given by the formula
\[ 
\omega_{g,x} = \SP{\rho_{x}^{*}pr_{g}^{*}\lcart, pr_{x}^{*}+ \frac{1}{2} \rho_{x}pr_{g}^{*}\lcart},
\]
or, explicitly,
\begin{equation} 
\omega_{g, x}((V, X), (V', X')) =\SP{\rho_{x}^{*}(\lcart_g(V)), \rho_{x}(\lcart_g(V'))} + 
\SP{\rho_{x}^{*}(\lcart_g(V)),X'}
 -\SP{\rho_{x}^{*}(\lcart_g(V')), X} \label{explicit}
\end{equation}
\end{enumerate}
\end{proposition}

\begin{proof} 
We first observe a few facts.
The defining formula for $\rho^*$ implies that
\begin{equation}\label{def-star}
\omega((v, 0), (0, X))= \SP{\rho^*(v), X},
\end{equation}
for $v\in \mathfrak{h}$, $X \in TM$.
Using (\ref{eq4}), we can write
\begin{eqnarray*}
\omega_{g, x}((V_g, 0), (V_{g}', X_{x}'))&=& \omega_{gx}(((dR_{g}^{-1})_{g}(V_{g}), 0), dt_{g,x}(V_{g}', X_{x}')) \\
                   &=& \omega_{gx}(((dR_{g}^{-1})_{g}(V_{g}), 0), dt_{gx}((dR_{g}^{-1})_{g}V_{g}', gX_{x}'))  \\
                   &=& \SP{\rho^{*}_{gx}(dR_{g}^{-1})_{g}(V_{g}), \rho_{gx}((dR_{g}^{-1})_{g}(V_{g}'))+ g X_{x}'} 
\end{eqnarray*}
for all $V_g, V_{g}'\in T_gH$, $X_{x}'\in T_xM$. (Here $gX$ denotes the infinitesimal action of $H$ on $TM$.)
Hence, for general pairs
$((V_g, X_x), (V_g', X_x'))$, we have
\begin{equation}\label{coc-c}
\omega_{g, x}((V_g, X_x), (V_g', X_x'))= \omega^{0}_{g, x}((V_g, X_x), (V_g', X_x')) + 
\omega_{g, x}((0_g, X_x), (0_g, X_{x}')),
\end{equation}
where
\begin{eqnarray} 
\omega^{0}_{g, x}((V_g, X_x), (V_g', X_x')) & = & 
         \SP{\rho^{*}_{gx}dR_{g^{-1}}(V_g), \rho_{gx}dR_{g^{-1}}(V_g')}
         + \SP{\rho^{*}_{gx}dR_{g^{-1}}(V_g), g X_x'}- \nonumber\\
 & & \SP{\rho^{*}_{gx}dR_{g^{-1}}(V_g'), gX_x} 
\label{for-0}
\end{eqnarray}
Let $\omega^1= \omega - \omega^0$. We make two simple remarks: First, 
$\omega^1$ encodes precisely the restrictions of $\omega$ to the slices $\{g\}\times M$ (see \eqref{coc-c});
Second, if $\rho^*$ is invariant, then (\ref{for-0}) coincides with the formula (\ref{explicit}) in the statement.
Hence it suffices to show that $\omega = \omega^0$ if and only if  $(\rho^*, \phi)$ is an equivariant form.
One possible route to prove that 
is as follows: one can show that $\omega^1$ (or, equivalently, $\omega^0$) is multiplicative if and only if 
$\rho^*$ is invariant, and $\omega^1$ is closed if and only if $i_{v}(\phi)- d(\rho^*(v))= 0$. 
Since $\rho^{*}_{\omega^1}= 0$, the uniqueness of Corollary \ref{unique} implies the proposition.
We will present an alternative argument instead. 

Let us rewrite $\omega^1=\omega - \omega^0$  as
$\SP{c(g), X_x\wedge X_x'}$, defining a smooth function $c\in C^{\infty}( H; \Omega^2(M))$. The
multiplicativity of $\omega$, applied on vectors $((0, X), (0, X'))$, reads
\begin{equation}
\omega_{hg, x}((0, X_{x}), (0, X_{x}'))=
\omega_{h, gx}((0, gX_{x}), (0, gX_{x}'))+ \omega_{g, x}((0, X_{x}), (0, X_{x}')),
\end{equation}
which precisely means that 
\begin{equation}\label{eq:cocycle} 
c(hg)= g^{*}c(h)+ c(g),\ \ \forall \ g, h, \in H
\end{equation}
(i.e., $c$ is an $\Omega^2(M)$-valued $1$-cocycle on $H$). 
Note that, in order to prove that $c=0$, it suffices to show that
$\Lie_{v}(c)= 0$ for all $v\in \mathfrak{h}$. 
Indeed, by differentiating \eqref{eq:cocycle}, we obtain
$\Lie_{V_g}(c)= 0$ for all $V_g= dR_g(v)\in T_g H$ and all $g\in H$, and $c$ must be
constant. Since, again by \eqref{eq:cocycle}, $c(1)$ is clearly zero, 
$c$ must vanish (see also Remark \ref{rem:cocycle}).

We now claim that  
\begin{equation}
\label{form-to-prove}
\Lie_{v}(c)= d_H(\rho^{*}+ \phi)(v) = d(\rho^{*}(v))- i_{v}(\phi)
\end{equation}
for all $v\in \mathfrak{h}$. 
In order to prove \eqref{form-to-prove}, let $V$ be a vector field on $H$ extending $v$, 
and let $X$ and $X'$ be vector fields  on $M$. We evaluate 
$d\omega= s^*\phi-t^*\phi$ on $(V, 0), (0, X), (0, X')$:
\begin{eqnarray*} 
d\omega((V, 0), (0, X), (0, X'))&=&
\Lie_{(V, 0)}(\omega((0, X), (0, X')))- \Lie_{(0, X)}(\omega((V, 0), (0, X'))) +  \\
& &  \Lie_{(0, X')}(\omega((V, 0), (0, X)))+ \omega((V, 0), (0, [X, X']))\\
&= & -\phi(\rho(\bar{\lcart}(V)), X, X'), 
\end{eqnarray*}
where $\bar{\lcart}(V)_{g}= (dR_g^{-1})_{g}(V_g)$. 
Using (\ref{def-star}) and evaluating the previous formula at $g=1\in H$, we find
\[
\Lie_{v}(c)(X, X')= \Lie_X(\rho^*(v)(X'))- \Lie_{X'}(\rho^*(v)(X))- \rho^*(v)([X, Y])- \phi(\rho(v), X, Y),
\]
which is just \eqref{form-to-prove}. This proves the proposition.
\end{proof}

\begin{remark}\label{rem:cocycle}\rm
We observe that the proof of Proposition \ref{recover-AMM} indicates how to express
$\omega_{\rho^*, \phi}$ for general pairs $(\rho^*, \phi)$ (i.e. which only satisfy
the conditions  (\ref{cond1}), (\ref{cond2})). 
More precisely, the cocycle condition \eqref{eq:cocycle} for $c$ is equivalent to
saying that $g \mapsto (g, c(g))$ is a group homomorphism
from $H$ into the group $H\ltimes \Omega^2(M)$ defined by
$(h, a)(g, b)= (hg, g^{*}a+ b)$. The proof of Proposition \ref{recover-AMM}
then shows that the induced Lie algebra
map $\mathfrak{h}\rmap \mathfrak{h}\ltimes \Omega^2(M)$ is 
\begin{equation}\label{eq:liealgcoc}
v\mapsto (v, d(\rho^{*}(v))- i_{v}(\phi)).
\end{equation}
So, if $H$ is simply connected, the Lie algebra cocycle
\eqref{eq:liealgcoc} integrates uniquely to a group cocycle $c$,
and $\omega$ will be given by (\ref{for-0}) plus $c(g)(X_x, Y_x)$.
\end{remark}

\begin{remark}\rm  In general, there is a natural map
\[ 
H^{*}(\Omega_{H}^{*}(M))\rmap H^{*}_{H}(M) 
\]
from the cohomology of the Cartan complex into the equivariant cohomology of $M$,
which is an isomorphism if $H$ is compact \cite{AtBo}. 
The equivariant cohomology groups can be obtained from
 a double complex $\Omega^{p}(H^{q}\times M)$,
with de Rham differential increasing the degree $p$, 
and a group-cohomology differential increasing $q$; see e.g. \cite{BSS}. 
Our result gives both an explicit description of this map in degree three, 
\[ 
\rho^*+ \phi\mapsto \omega_{\rho^{*}, \phi}+ \phi ,
\]
as well as an interpretation of this map in terms of multiplicative forms.
\end{remark}

If $\rho^*+ \phi\in \Omega_{H}^{3}(M)$ satisfies the non-degeneracy
condition $\dim(\Ker(\rho)\cap \Ker(\rho^*))=\dim(H)-\dim(M)$, then  
$(H\ltimes M, \omega_{\rho^{*}, \phi})$ becomes
an over-presymplectic groupoid, which is presymplectic if $\dim(M)= \dim(H)$. 
The associated $\phi$-twisted Dirac structure can be described directly as
\[ 
L= \{ (\rho(v), \rho^*(v)): v\in \mathfrak{h} \}\subset TM\oplus T^*M .
\]
A simple example is $M= \mathfrak{h}^{*}$ with the coadjoint action
of $H$, $\phi= 0$, and $\rho^*$ given by $\rho^{*}_{\xi}(v)= v$. The associated 
groupoid will be $H\ltimes \mathfrak{h}^* \cong T^*H$, with the canonical symplectic form.
A more interesting example will be the AMM groupoid of Subsection \ref{ex:AMM}.

%%%%%%%%%%%%%%%%%%%%%%%%%%%%%%%%%%%%%%%
%%%%%%%%%%%%%%%%%%%%%%%%%%%%%%%%%%%%%%%
%%%%%%%%%%%%%%%%%%%%%%%%%%%%%%%%%%%%%%%
%%%%%%%%%%%%%%%%%%%%%%%%%%%%%%%%%%%%%%%
\section{Presymplectic realizations of Dirac structures}
\label{sec-realizations}
%%%%%%%%%%%%%%%%%%%%%%%%%%%%%%%%%%%%%%%
%%%%%%%%%%%%%%%%%%%%%%%%%%%%%%%%%%%%%%%
%%%%%%%%%%%%%%%%%%%%%%%%%%%%%%%%%%%%%%%
%%%%%%%%%%%%%%%%%%%%%%%%%%%%%%%%%%%%%%%
Let $(M,\pi)$ be a Poisson manifold. Recall that a {\bf symplectic realization}
of $M$ is a Poisson map from a symplectic manifold $(P,\eta)$ to $M$ (see e.g. \cite{CaWe}). 
The following important property of symplectic realizations brings them close to
the theory of hamiltonian actions: any symplectic realization $\mu: P \to M$ induces
a canonical action of the Lie algebroid $T^*M$, induced by $\pi$, on $P$ by assigning
to each $\alpha \in \Omega^1(M)$ the vector field $X \in \mathcal{X}(P)$ defined by
$$
i_X\eta = \mu^*\alpha.
$$
When $\mu$ complete (i.e., the hamiltonian vector field $X_{\mu^*f}$ is complete
whenever $f\in C^\infty(M)$ has compact support), $M$ is integrable \cite{CrFe2} and
this action extends to a symplectic action of $G$, the $s$-simply connected symplectic groupoid
of $M$ (see \cite{CDW,CrFe2}), with moment map $\mu$ \cite{MiWe}. 
In this way, we get a natural correspondence between symplectic actions of $G$ and 
complete symplectic realizations of $M$. 
In particular, if $M = \mathfrak{h}^*$ is the dual of a Lie algebra, the action of
the associated symplectic groupoid $T^*H=H\ltimes \mathfrak{h}^*$ factors through 
an $H$-action, and complete symplectic realizations of $\mathfrak{h}^*$ become
hamiltonian $H$-spaces. In this section, we will extend this picture to
twisted Dirac manifolds.

\subsection{Presymplectic realizations}

We recall the definition introduced  in Section \ref{Dirac-pre}.
\begin{definition}\label{def:presympreal}
 A {\bf presymplectic realization}
of a $\phi$-twisted Dirac manifold $(M, L)$ 
is a Dirac map $\mu: (P, \eta)\rmap (M, L) $, where $\eta$ is a $\mu^*\phi$-closed 
2-form (i.e., $d\eta+ \mu^{*}\phi= 0$), such that $\Ker(d\mu)\cap \Ker(\eta)=\{0\}$.
\end{definition}

The following results explain this definition.

\begin{lemma}
\label{lemma-real} 
Let $(M, L)$ be a $\phi$-twisted Dirac manifold, let
$\mu: P\rmap M $ be a smooth map, and let $P$ be equipped with a 2-form $\eta$ 
satisfying $d\eta+ \mu^{*}\phi= 0$. The following are equivalent:
\begin{itemize}
\item[(i)] The map $\mu$ is a presymplectic realization of $L$;
\item[(ii)] For all $p\in P,\ (w, \xi)\in L_{\mu(p)}$,
there exists a unique $X \in T_pP$ satisfying the equations:
\[ \left \{ \begin{array}{ll} 
         w = d\mu(X)\\
         \mu^*(\xi) = i_{X}(\eta)
           \end{array}
  \right. ,\]
\item[(iii)] The map $\mu$ is Dirac and $d\mu$ maps $\Ker(\eta)$ isomorphically onto $\Ker(L)$.
\end{itemize}

\end{lemma}

\begin{proof}
Note that $\mu$ being a Dirac map is equivalent to the equations in (ii) 
having a solution for $X$. The uniqueness of the solutions is equivalent
to $\Ker(d\mu) \cap \Ker(\eta)=\{0\}$, so (i) and (ii) are equivalent.
Note that $\Ker(L)=\{w = d\mu(X) \; |\;i_X\eta =0\}=d\mu(\Ker(\eta))$,
so $d\mu:\Ker(\eta) \to \Ker(L)$ is an isomorphism if and only 
$\Ker(d\mu)\cap \Ker(\eta) = \{0\}$. Hence all the conditions are equivalent.
\end{proof}

Note that if the conditions in Lemma \ref{lemma-real} hold, then (ii) defines a map 
$\rho_{P}: L_{\mu(p)}\rmap T_{p}P$, $(w,\xi)\mapsto X$.
A direct computation shows that
\begin{enumerate}
\item the induced map $ \rho_{P}: \Gamma(L)\rmap \mathcal{X}(P)$ is a map 
of Lie algebras ($\Gamma(L)$ is equipped with twisted Courant bracket);
\item $d\mu(\rho_{P}(l)) = \rho(l)$ for all $l\in L$, 
\end{enumerate}
which precisely means that $ \rho_{P}$ is an infinitesimal action of the Lie 
algebroid $L$ on $P$ (and this is what we were after!).

\begin{corollary} 
Any presymplectic realization $\mu: (P, \eta)\rmap (M, L)$ of a $\phi$-twisted Dirac
structure is canonically equipped with an infinitesimal action of the Lie algebroid of
$L$. 
\end{corollary}

We call a presymplectic realization $\mu:P\to M$ {\bf complete} 
if $ \rho_{P}(l)$ is a complete vector field 
whenever $l\in \Gamma(L)$ has compact support. 
As with symplectic realizations of Poisson manifolds, a complete realization
defines a complete Lie algebroid action \cite{MoMr}, 
which can be integrated to a global action of the groupoid
$G(L)$ associated with $L$ on $P$: 
indeed, an  algebroid action of $L$ defines a map
$\nabla: \Gamma(L)\otimes C^{\infty}(P)\rmap C^{\infty}(P)$ which behaves like
a flat ($L$-) connection; then parallel transport defines the desired action
of $G(L)$ on $P$ (see \cite[pp.26--27]{CrFe2} for details). For the integration of
general Lie algebroid actions, see \cite{MoMr}.

For complete symplectic realizations of Poisson manifolds, the
induced action of the symplectic groupoid is {\it symplectic} \cite{MiWe}.
The following property generalizes this fact.
Let $(M,L)$ be a $\phi$-twisted Dirac structure, and let $(G(L),\omega_L)$
be the associated groupoid.

\begin{corollary} \label{cor:realiz}
If the realization $\mu:P \to M$ is complete (e.g., if $P$ is compact), then there is an induced
action of $G(L)$ on $P$, $m_{P}: G(L)\times_{M}P\rmap P$. Moreover, if $G(L)$ is smooth (i.e., if $L$ is integrable), then the
action is smooth and
\begin{equation}
\label{mult-on-P} 
m_{P}^{*}\eta= pr_{G}^{*}\omega_{L}+ pr_{P}^{*}\eta 
\end{equation}
(where $pr_{G}: G(L)\times_{M}P\rmap G(L)$ and $pr_{P}: G(L)\times_{M}P\rmap P$
are the natural projections). 
\end{corollary}

\begin{proof}
In order to check \eqref{mult-on-P}, note that the 2-forms 
$\omega_1 := m_{P}^{*}\eta$ and $\omega_2:= pr_{G}^{*}\omega_{L}+ pr_{P}^{*}\eta$ 
are both multiplicative in
the semi-direct product groupoid $G(L)\ltimes P$.
A direct computation shows that $\rho^*_{\omega_1} = \rho^*_{\omega_2}$, so
it follows from Theorem \ref{rec1} that the forms must coincide.
\end{proof}

\begin{remark}\rm By the same arguments as in \cite[Thm.~8.2]{CrFe2}, the existence
of a complete presymplectic realization $\mu:P \to M$ which is a surjective submersion implies
the integrability of $L$. Note also that such realizations $\mu$ can be used to compute
$(G(L), \omega_L)$ (though we do not know how to use this to give a direct proof of the integrability of $L$):
First, we note that
the groupoid $G(L)\ltimes P$ over $P$ is isomorphic
to the monodromy groupoid $G(\Im( \rho_{P}))$ of the (regular) foliation
$\Im( \rho_{P})$, so that $G(L)$ is a quotient of $G(\Im( \rho_{P}))$;
second, \eqref{mult-on-P} says that the
form $t^*\eta- s^*\eta$ on $G(\Im( \rho_{P}))$ descends to
$\omega_L$ on $G(L)$. 
\end{remark}

%%%%%%%%%%%%%%%%%%%%%%%%%%%%%%%%%%%%%%%
%%%%%%%%%%%%%%%%%%%%%%%%%%%%%%%%%%%%%%%
%%%%%%%%%%%%%%%%%%%%%%%%%%%%%%%%%%%%%%%
%%%%%%%%%%%%%%%%%%%%%%%%%%%%%%%%%%%%%%%
\subsection{Realizations of Cartan-Dirac structures and quasi-hamiltonian spaces}
\label{ex:AMM}
%%%%%%%%%%%%%%%%%%%%%%%%%%%%%%%%%%%%%%%
%%%%%%%%%%%%%%%%%%%%%%%%%%%%%%%%%%%%%%%
%%%%%%%%%%%%%%%%%%%%%%%%%%%%%%%%%%%%%%%
%%%%%%%%%%%%%%%%%%%%%%%%%%%%%%%%%%%%%%%

As observed in \cite[Example~4.2]{WeSe}, any Lie group with a bi-invariant metric
carries a $\phi$-twisted Dirac structure, where $\phi$ is the associated bi-invariant
Cartan form. We call it a {\bf Cartan-Dirac structure}. In this section, we will discuss
presymplectic realizations  and groupoids of Cartan-Dirac structures. We recover, in this framework,
quasi-hamiltonian spaces \cite{AMM} and the AMM-groupoid of \cite{bxz},
proving the following result.

\begin{theorem}\label{thm:AMM}
Let $H$ be a connected Lie group, and let 
$(\cdot,\cdot)_{\mathfrak{h}}$ be an invariant inner product on its Lie algebra $\mathfrak{h}$.
Let $L$ denote the associated Cartan-Dirac structure on $H$.
Then 
\begin{itemize}
\item[(i)] There is a one-to-one correspondence between presymplectic realizations of $(H,L)$
and quasi-hamiltonian $\mathfrak{h}$-spaces (which are infinitesimal versions of
quasi-hamiltonian $H$-spaces introduced in \cite{AMM}); 
\item[(ii)] The AMM-groupoid
of \cite{bxz} is a presymplectic groupoid inducing the Cartan-Dirac structure $L$ on $H$.
\end{itemize}
\end{theorem}

Before we prove this theorem, let us recall some definitions and fix our notation.

A {\bf quasi-hamiltonian $H$-space} \cite{AMM} 
is a manifold $P$ endowed with a smooth action of $H$, an invariant 2-form $\eta\in \Omega^2(P)$,
and an equivariant map $\mu: P\rmap H$ (the moment map), such that 
\begin{enumerate}
\item the differential of $\eta$ is given by
\[ 
d\eta= -\mu^*\phi ;
\]
\item the map $\mu$ satisfies 
\[ 
i_{\rho_{P}(v)}(\eta)= \frac{1}{2}\mu^* (\lcart+ \bar{\lcart}, v)_{\mathfrak{h}};
\]
\item at each $p\in P$, the kernel of $\eta_p$ is given by
\[ 
\Ker(\eta_{p})=  \{\rho_{P, p}(v): v \in \Ker(\Ad_{\mu(p)}+ 1)\} .
\]
\end{enumerate}
Here, $\rho_{P}: \mathfrak{h}\rmap TP$ is the induced infinitesimal action of 
$\mathfrak{h}$ on $P$, $\lcart$ (resp. $\bar{\lcart}$) is the left- (resp. right-)
invariant Maurer-Cartan form on $H$, and $\phi\in \Omega^3(H)$ is the bi-invariant Cartan form:
\[
\phi= \frac{1}{12} ( \lcart, [\lcart,
\lcart])_{\mathfrak{h}}=\frac{1}{12}(\bar{\lcart},[\bar{\lcart},\bar{\lcart}])_{\mathfrak{h}}.
\]
On the Lie algebra, we have
$\phi(u,v,w)=\frac{1}{2}(u,[v,w])_{\mathfrak{h}}$.
The equivariance
of $\mu$ is with respect to the action of $H$ on itself by conjugation. 
Infinitesimally, equivariance becomes
\begin{equation}\label{eq:infequiv} 
(d\mu)_{p}(\rho_{P}(v))= \rho_{H}(v) 
\end{equation}
for all $v\in \mathfrak{h}$, where $\rho_{H}: \mathfrak{h}\rmap TH$ 
is the infinitesimal conjugation action
(explicitly, $\rho_{H}(v)= v_r- v_l$, where $v_l$ and $v_r$ 
are the vector fields obtained from $v \in \mathfrak{h}$ by
left and right translations). 

\begin{definition} \label{def:infquasi}
A {\bf quasi-hamiltonian $\mathfrak{h}$-space} is a manifold $P$ carrying an $\mathfrak{h}$-action
$\rho_{P}: \mathfrak{h}\rmap TP$,
together with an $\mathfrak{h}$-invariant 2-form $\eta \in \Omega^2(P)$ and
and equivariant $\mu: P \to H$ (as in \eqref{eq:infequiv}), satisfying conditions
(1), (2) and (3).
\end{definition}

%%%%%%%%%%%%%%%%%%%%%%%%%%%%%%%%%%%%%%%%%%%%%%%%%%%%%%%%%%%%%%%%%%%%%%%%%%%%%%%%%%%%%%%%
Conditions (1),(2) and (3) in the definition of quasi-hamiltonian spaces
strongly resemble the conditions we used to define presymplectic realizations,
see Lemma \ref{lemma-real}, (ii).
In order to find the  underlying  Dirac structure $L$ on $H$ making
quasi-hamiltonian $\mathfrak{h}$-spaces into presymplectic realizations,
recall that an $\mathfrak{h}$-action on $P$ together with an
equivariant map $\mu:P \to H$ is equivalent to an action of the action
Lie algebroid $\mathfrak{h}\ltimes H$, with moment $\mu$ \cite{MiWe}.
Hence the Lie algebroid of $L$ is isomorphic to $\mathfrak{h}\ltimes H $, with anchor $\rho_H$;
in other words, there is a map $\rho^*:\mathfrak{h}\to T^*H$ such that
$(\rho_H,\rho^*):\mathfrak{h}\ltimes H\to L $ is an isomorphism.
To find $\rho^*$, we compare condition (2) for quasi-hamiltonian spaces
and the second equation in Lemma \ref{lemma-real}, (ii), and obtain
\[ 
\rho^*(v) = \frac{1}{2}(\lcart+ \bar{\lcart}, v)_{\mathfrak{h}}= \frac{1}{2}(v_r+ v_l),
\]
where in the last equality we use the metric to identify $T^*H$ with $TH$. 
More explicitly, the Dirac structure we obtain on $H$ is
\[ 
L=\{(v_r-v_l, \frac{1}{2}(v_r+v_l)):\;v\in\mathfrak{h}\}\subset TH \oplus TH,
\]
which is precisely the
$\phi$-twisted Dirac structure discussed in \cite[Example~4.2]{WeSe}.
We call $L$ the {\bf Cartan-Dirac structure}
on $H$ associated with $(\cdot,\cdot)_{\mathfrak{h}}$. 

We now proceed to the proof of Theorem \ref{thm:AMM}.

\begin{proof}
Suppose that
$\mu:P \to H$ is a presymplectic realization of the Cartan-Dirac structure $L$ on $H$.
Let $\rho_{P}^L$ be the induced infinitesimal action of $L$ on $P$, with moment $\mu$.
Since the Lie algebroid of $L$ is isomorphic to $\mathfrak{h}\ltimes H$, it immediately
follows that $\rho_{P}^L$ determines an $\mathfrak{h}$-action $\rho_P$ on $P$, for which
$\mu$ is $\mathfrak{h}$-equivariant. Explicitly,
\begin{equation}
\label{rho-versus} 
\rho_{P}(v)= \rho_{P}^{L}(v_r- v_l, \frac{1}{2}(\lcart+ \bar{\lcart}, v)_{\mathfrak{h}}),\;\; 
v\in \mathfrak{h}. 
\end{equation}
Since $d\eta + \mu^*\phi =0$, (1) in Definition \ref{def:infquasi} holds.
Condition (2) is just the second equation in Lemma \ref{lemma-real}, (ii). 
Since $d\mu:\Ker(\eta) \to \Ker(L)$ is an isomorphism, and 
$\Ker(L_g)= \{ \rho_{H}(v): v\in \Ker(\Ad_g+ 1)\}$, the equivariance
of $\mu$, $d\mu(\rho_P(v))= \rho_H(v)$ implies condition (3).
Finally, we note that $\eta$ is $\mathfrak{h}$-invariant:
$$
\Lie_{\rho_P(v)}(\eta) = di_{\rho_P(v)}\eta + i_{\rho_P(v)}d\eta
= \frac{1}{2}d(\mu^*(\lcart+ \bar{\lcart}, v)_{\mathfrak{h}})) - i_{\rho_{P}(v)}\mu^*\phi = 0,
$$
where the last equality follows from the Maurer-Cartan equations for $\lcart$ and $\bar{\lcart}$.
So $P$ is a quasi-hamiltonian $\mathfrak{h}$-space.

Conversely,  if $P$ is a quasi-hamiltonian $\mathfrak{h}$-space, then
$d\eta + \mu^*\phi =0$, and we must check that condition (ii) in
Lemma \ref{lemma-real} holds. If $(w,\xi) = 
(v_r- v_l, \frac{1}{2}(\lcart+ \bar{\lcart}, v)_{\mathfrak{h}}) \in L$, then 
$w=\rho_H(v)=d\mu(\rho_P(v))$, so the 
first equation in (ii) has a solution $X=\rho_P(v)$. The second equation in (ii)
is just (2) in Definition \ref{def:infquasi}. The uniqueness
of this solution follows from (3): if $i_X\eta =0$, then $X=\rho_P(v)$ for some
$v$ with $v_l+v_r=0$; since $d\mu(X)=\rho_H(v)= v_r-v_l=0$, we must have $v=0$.
So $\mu:P \to H$ is a presymplectic realization. This proves part (i) of
the theorem.

In order to prove part (ii), 
we describe the presymplectic groupoids associated with the Cartan-Dirac
structure $L$. Since $L$, as a Lie algebroid, is isomorphic to the action Lie algebroid
$\mathfrak{h}\ltimes H$,
the action groupoid $H\ltimes H$ (with action given by conjugation,  
$g\cdot x= gxg^{-1}$) integrates it.
(On $H\ltimes H$, $s(g, x)= x$, $t(g, x)= g x g^{-1}$, $(g_1, x_1)\cdot (g_2, x_2)= (g_1 g_2, x_2)$.)
We are exactly in the situation of Example \ref{Lie group actions}. 
As observed in \cite{AMM},
$\rho^* + \phi\in \Omega^{3}_{H}(H)$. 
(In particular, it follows from this observation that $L$ is indeed a $\phi$-twisted Dirac structure.)
Applying Proposition \ref{recover-AMM}, we immediately obtain the formula for the 
multiplicative 2-form on $H\times H$ corresponding to $L$:
\[ 
\omega_{(g, x)} = \frac{1}{2} \left((\mathrm{Ad}_{x} p_g^* \lambda, p_g^*
\lambda)_{\mathfrak{h}} +
( p_g^*\lambda, p_x^*(\lambda+\bar{\lambda}))_{\mathfrak{h}}\right).
\]
As in Proposition \ref{recover-AMM}, $p_g$ and $p_x$ denote the projections onto the first and second
components of $H\times H$.
This is precisely the 2-form in the ``double'' $D(H)=H\times H$ introduced
 in \cite{AMM}; the groupoid $(H\ltimes H, \omega)$ also appears
in \cite{bxz}, where it is called the {\bf AMM groupoid}.
So the AMM groupoid is a presymplectic groupoid
associated with the Cartan-Dirac structure on $H$, though it is not necessarily $s$-simply connected.
If $H$ is simply connected, then  $(G(L), \omega_L)= (H\ltimes H, \omega)$,
while, in general, one must pull-back $\omega$ to $\tilde{H}\ltimes H$, where $\tilde{H}$ is
the universal cover of $H$.
\end{proof} 

Finally, if the infinitesimal action can be integrated to a global action,
then the notion of quasi-hamiltonian $H$-space coincides with that of 
quasi-hamiltonian $\mathfrak{h}$-space.
For instance, if $H$ is simply connected, there is a one-to-one correspondence 
between  quasi-hamiltonian $H$-spaces 
and complete 
presymplectic realizations of $L$ (which, in turn, are equivalent to
quasi-hamiltonian $\mathfrak{h}$-spaces for which the $\mathfrak{h}$-action
is by complete vector fields).
In particular,

\begin{corollary} If $H$ is simply connected, then there is a 
one-to-one correspondence between compact quasi-hamiltonian $H$-spaces
and compact presymplectic realizations of the Cartan-Dirac structure $L$ on $H$. 
\end{corollary}

%%%%%%%%%%%%%%%%%%%%%%%%%%%%%%%%%%%%%%%
%%%%%%%%%%%%%%%%%%%%%%%%%%%%%%%%%%%%%%%
%%%%%%%%%%%%%%%%%%%%%%%%%%%%%%%%%%%%%%%
%%%%%%%%%%%%%%%%%%%%%%%%%%%%%%%%%%%%%%%
\section{Multiplicative 2-forms, foliations and regular Dirac structures}
\label{sec-foliations}
%%%%%%%%%%%%%%%%%%%%%%%%%%%%%%%%%%%%%%%
%%%%%%%%%%%%%%%%%%%%%%%%%%%%%%%%%%%%%%%
%%%%%%%%%%%%%%%%%%%%%%%%%%%%%%%%%%%%%%%
%%%%%%%%%%%%%%%%%%%%%%%%%%%%%%%%%%%%%%%

In this section, we explain connections  between our results and
some aspects of foliation theory: we will see that multiplicative 2-forms
on monodromy groupoids of foliations are directly related to 
foliated cohomology, and
they are relevant for the explicit description of presymplectic groupoids
associated to {\it regular} Dirac structures. 

\subsection{Foliations}\label{sec:foliation}
Let us recall some basic facts of foliations theory. The reader is referred to
\cite{KmTd} and references therein for details.

By the Frobenius theorem, a foliation on $M$ can be viewed
as a subbundle $\F$ of $TM$ (of vectors tangent to the leaves) for which $[\F,
\F]\subset \F$; alternatively,
foliations are the same thing as algebroids with injective anchor
map. The {\bf monodromy groupoid} of $\F$ consists of leafwise homotopy classes of leafwise
paths in $M$ (i.e., each $s$-fiber $s^{-1}(x)$ is
the universal cover of the leaf through $x$, constructed with
$x$ as base point). This groupoid is the same as the one described in
Section \ref{Reconstructing multiplicative forms}, i.e., it is the unique
$s$-simply connected Lie groupoid integrating $\F$ viewed as an algebroid.
We denote it by $G(\F)$.
The space of foliated forms on $M$,
$\Omega^\bullet(\F)=\Gamma(\wedge^\bullet \F^*)$,
carries a foliated de Rham operator
\begin{align}\label{eq:foldiff}
d_{\F}\omega(X_1,\dots,X_{p+1})=&
\sum_{i}(-1)^{i}\mathcal{L}_{X_i}(\omega(X_1,\dots,\hat{X_i},\dots,X_{p+1}))\\\nonumber
+\sum_{i<j}(-1)&^{i+j-1}\omega([X_i, X_j],
X_1,\dots,\hat{X_i},\dots, \hat{X_j},\dots,X_{p+1})),
\end{align}
and we denote by $H^\bullet(\F)$ the resulting cohomology (which is just the cohomology of $\F$
as an algebroid).
One defines, in a similar way, the foliated cohomology $H^\bullet(\F; E)$ with
coefficients in a foliated bundle $E$,
i.e., a bundle $E$ over $M$ endowed with a flat $\F$-connection 
$\nabla: \Gamma(\F)\times \Gamma(E)\to\Gamma(E)$.
The corresponding complex is now $\Omega^\bullet(\F; E)=\Gamma(\wedge^\bullet \F^*\otimes
E)$, and the differential
is given just as in \eqref{eq:foldiff}, with $\mathcal{L}_{X}$ replaced by
$\nabla_X$. The basic example of a
foliated bundle is the normal bundle $\nu=TM/\F$, with $\F$-connection given by the
well-known Bott connection, $\nabla: \Gamma(\F)\times \Gamma(\nu)\to\Gamma(\nu)$, 
$$ 
\nabla_{V} \overline{X}= \overline{[V, X]},
$$
where $X \mapsto \overline{X}$ is the projection from $TM$ onto $\nu$.
As usual, the Bott connection induces connections on the dual $\nu^*$ and on the associated
tensor bundles. 

Here we will deal with the cohomology spaces $H^\bullet(\F)$
and $H^\bullet(\F; \nu^*)$ 
(in degrees one and two). These spaces are relate by a transversal de Rham
operator
\begin{equation}\label{eq:transv} 
d_{\nu}: H^\bullet(\F)\rmap H^\bullet(\F; \nu^*)
\end{equation}
(which can be extended to higher exterior powers of $\nu^*$). 
Let us give a direct description of this map in degree two, since
higher degrees can be treated analogously. 
Given a class $[\theta]\in H^2(\F)$ represented by a foliated 2-form $\theta$,
let $\widetilde{\theta}$ be a 2-form on $M$ with 
$\theta= \widetilde{\theta}|_{\F}$. Since
$d\widetilde{\theta}|_{\F}= 0$, it follows that the map
$\Gamma(\wedge^2\F)\to\Gamma(\nu^*)$ defined by
\[ 
(V,W)\mapsto d\widetilde{\theta}(V, W, -), 
\] 
gives a closed foliated 2-form with coefficients in $\nu^*$; we set $d_{\nu}([\theta])$
to be its class in
$H^2(\F; \nu^*)$.

\subsection{Multiplicative 2-forms on monodromy groupoids}
\label{two-forms on monodromy}

In this example we relate the space of closed 
multiplicative 2-forms on monodromy groupoids to cohomology spaces which are
well known in foliation theory.

Let $Mult^2(G)$ denote the space of closed multiplicative 2-forms
on a Lie groupoid $G$.
 
\begin{proposition}
\label{ex:monodromy} Let $\F$ be a foliation on $M$, and let $G= G(\mathcal{F})$
be the monodromy groupoid
of $\F$. Then
\begin{enumerate}
\item[(i)] any $\omega\in Mult^2(G)$ induces a cohomology class $c(\omega) \in
H^2(\F)$.
\item[(ii)] a foliated cohomology class $c\in H^2(\F)$  is of type $c(\omega)$
if and only if $d_{\nu}(c)= 0$.
\item[(iii)] $c(\omega)= 0$ if and only if $\omega$ is multiplicatively exact,
i.e. $\omega= d\sigma$ with $\sigma\in \Omega^1(G)$ multiplicative.
\end{enumerate}
\end{proposition}

Before we prove this proposition, let us recall a more conceptual way to describe
the operator \eqref{eq:transv}. There is a  spectral sequence associated to
the foliation  $\F$ (see e.g \cite{KmTd}),
converging to $H^\bullet(M)$, with
\begin{equation}
\label{E1-term} 
E^{p, q}_{1}= H^{p}(\F; \Lambda^q \nu^*) ,
\end{equation}
and so that $d_{\nu}$ is just the boundary map $d_{1}^{p, q}: E^{p, q}_1\rmap E^{p, q+1}_1$.
The spectral sequence is
associated to the filtration $F_{p}\Omega^\bullet(M)$ of $\Omega^\bullet(M)$ with
\[ 
F_{q}\Omega^n(M)= \{ \eta\in\Omega^n(M): i_{V_{1}}\ldots i_{V_{n-q+1}}\eta =
0,\ \text{for\ all}\ V_{i}\in \Gamma(\F)\} 
\]
($F_{0}\Omega^\bullet(M)= \Omega^\bullet(M)$, and $F_{q}\Omega^n(M)= 0$ 
for $q>n$). It is easy to see that
\[ 
E^{p, q}_{0}= F_{q}\Omega^{p+q}(M)/F_{q+1}\Omega^{p+q+1}(M)\cong
\Omega^{p}(\F; \Lambda^{q}\nu^*) 
\]
and that the boundary $d^{p, q}_{0}$ is precisely the leafwise de Rham operator
$d_{\F}$. 
This shows that the $E_{1}$-terms are indeed given by (\ref{E1-term}), and a
standard computation shows
that $d^{p, q}_{1}$ has the explicit description mentioned above.\\

\begin{proof}
Note that we have an isomorphism
\[ 
\frac{F_0\Omega^2(M)}{F_2\Omega^2(M)}\stackrel{\sim}{\longrightarrow} \{\rho^*:\F\rmap T^*M: \rho^* \
\text{satisfies}\ (\ref{cond1})\},
\]
sending $[\eta]$ to $\rho^*$, defined by $\langle \rho^*(V), X\rangle=
\eta(V, X)$. Moreover, 
the closedness of $[\eta]$ in the complex
$F_0\Omega^\bullet(M)/F_2\Omega^\bullet(M)$ corresponds to
(\ref{cond2}) for $\rho^*$. On the other hand, $\rho^*$ corresponds
to an exact $[\eta]$ if and only if
$\SP{\rho^*(V), X}= d\sigma(V, X)$ for some
$\sigma\in F_0\Omega^1(M)/F_2\Omega^1(M)= \Omega^1(M)$.
But then the closed multiplicative 2-form $\omega$ associated with
$\rho^*$ is $\omega_0=d(t^*\sigma- s^*\sigma)$
(note that $\omega_0$ is
multiplicative and closed, and it is easy to see that 
$\rho^{*}_{\omega_0} = \rho^*$). As a result, we get
an isomorphism
\[ 
H^2\left(\frac{F_0\Omega^\bullet(M)}{F_2\Omega^\bullet(M)} \right) \stackrel{\sim}{\longrightarrow}
\frac{Mult^2(G)}{\{ d(t^*\sigma- s^*\sigma): \sigma\in \Omega^1(M)\}} .
\]
Now, using the short exact sequence of complexes
\[ 
0\rmap \frac{F_1\Omega^\bullet(M)}{F_2\Omega^\bullet(M)}\rmap
\frac{F_0\Omega^\bullet(M)}{F_2\Omega^\bullet(M)}\rmap
\frac{F_0\Omega^\bullet(M)}{F_1\Omega^\bullet(M)}\rmap 0 ,
\]
we get an exact sequence in cohomology
\begin{equation}
\label{ex-seq}
H^{1}(\mathcal{F})\stackrel{d_{\nu}}{\rmap} H^{1}(\mathcal{F};\nu^*)\rmap
\frac{Mult^{2}(G(\mathcal{F}))}{\{ d(t^*\sigma- s^*\sigma): 
\sigma\in \Omega^1(M)\}} \stackrel{c}{\rmap}
H^{2}(\mathcal{F})\stackrel{d_{\nu}}{\rmap} H^{2}(\mathcal{F};\nu^*) .
\end{equation}
This immediately implies statements (i) and (ii). 
Note that the map $c$ in \eqref{ex-seq} is given by 
$c(\omega)= [c_{\omega}]$, where
$c_{\omega}\in \Omega^2(\F)$ is defined by $c_{\omega}(V, W)= \SP{\rho^{*}_{\omega}(V), W}$. 

We now prove (iii). The fact that $c_{\omega}= 0$ for multiplicative 2-forms 
of type $\omega= d\sigma$, with $\sigma$ multiplicative, follows by showing
that the restriction
of $\sigma$ to $\F$ (a foliated $1$-form) gives, after differentiation, precisely
the foliated 2-form $c_{\omega}$ induced by $\omega$. This can be checked
by an argument similar to the one used to prove the
formula of Proposition \ref{rho}, part (ii) (but the argument is simpler, following from Remarks
2) and an analogue of 3) in that proof). For the converse, we fix $\omega$
with $[c_{\omega}]= 0$. From the exact sequence \eqref{ex-seq}, we may assume
that  $c_{\omega}= 0$. 
But then Lemma \ref{lema2} shows that the $1$-form $\tilde{\sigma}=
-\rho^*\sigma^c$ is basic, so it descends to a
multiplicative $1$-form $\sigma$ on $G(\F)= G$. Since
$\omega$ is induced by $\tilde{\omega}= d\tilde{\sigma}$, we conclude that
$\omega= d\sigma$.
\end{proof}

\subsection{Dirac structures associated to foliations}
\label{Dir-fol}

Any regular foliation $\mathcal{F}$ on $M$ defines a Dirac structure 
$L_{\mathcal{F}}$ whose presymplectic leaves are precisely the leaves of
$\mathcal{F}$,
with the zero form. In other words,
\[ 
L_{\mathcal{F}}= \mathcal{F}\oplus \nu^{*}\subset TM\oplus T^*M .
\]
The Dirac structure $L_{\mathcal{F}}$ is always integrable, and we now describe
the associated presymplectic groupoid. As above, we denote by $G(\mathcal{F})$ the 
monodromy groupoid of $\mathcal{F}$, and by $\nu$ the normal bundle. The
parallel
transport with respect to the Bott connection $\nabla$ (see Section \ref{sec:foliation}) 
is well defined along
leafwise paths (because $\nabla$ is an $\mathcal{F}$-connection), and is
invariant under
leafwise homotopy (because $\nabla$ is flat). Hence it defines an action of
$G(\mathcal{F})$ on
$\nu$, which we dualize to an action on $\nu^*$. We form the 
semi-direct product groupoid
\[ 
 G(\mathcal{F})\ltimes \nu^* 
\]
consisting of pairs $(g, v)$ with $v\in \nu_{s(g)}^{*}$, source and
target maps induced by those of $G(\mathcal{F})$, and  multiplication
\[ 
(g, v)(h, w)= (gh, h^{-1}v+ w).
\]
Restricting the canonical symplectic form $\omega_{can}$ on $T^*M$ to $\nu^*$,
its pull-back
\[ 
\omega_{\mathcal{F}}= pr_{2}^{*}\omega_{can} 
\]
by the second projection is a multiplicative 2-form on $G(\mathcal{F})\ltimes \nu^*  $.
It is not difficult to check the following

\begin{lemma}\label{lem:folgr}
$(G(\mathcal{F})\ltimes \nu^*, \omega_{\mathcal{F}})$ is the presymplectic groupoid associated with
$L_{\mathcal{F}}$.
\end{lemma}

%\begin{example}
%\label{Dir-fol-tw}
%({\it A simple twist})\rm

A slight ``twist'' of this result will yield more examples of multiplicative 2-forms
which are not of Dirac type.
 
Let us consider a closed 3-form $\phi$ on $M$
with the property that
\[ 
i_{V}i_{W}\phi = 0, \ \ \forall \ V, W\in \Gamma(\mathcal{F}).
\]
Using the filtration of Section \ref{sec:foliation}, this condition means that
$\phi\in F_{2}\Omega^3(M)$.
By Theorem \ref{rec1}, applied with $\rho^*=0$, 
there exists a unique multiplicative 2-form $\omega_{\phi}$ on
$G(\mathcal{F})$ such that 
\[ 
d\omega_{\phi}= s^*\phi- t^*\phi, \ \ \ \omega_{\phi, x}= 0 \ \ \forall\ x\in
M .
\]

\begin{proposition} The following are equivalent:
\begin{enumerate}
\item[(i)] $\omega_{\phi}= 0$;
\item[(ii)] $\omega_{\phi}$ is of Dirac type;
\item[(iii)] $\phi\in F_{3}\Omega^3(M)$ (or, equivalently, $\phi$ is basic).
\end{enumerate}
\end{proposition}

\begin{proof} 
By Lemma \ref{lemma},
both $\Ker(ds)$ and $\Ker(dt)$ sit inside $\Ker(\omega_{\phi})$ at all points $g\in
G(\mathcal{F})$.  This implies that $\Ker(ds)_{g}^{\perp}=
T_{g}G(\mathcal{F})$  and
$\Ker(dt)_{g}+ \Ker(\omega_{\phi, g})= \Ker(\omega_{\phi, g})$, and
the equivalence of (ii) and (i) follows from Lemma
\ref{lemma-Dirac}.

Next, note that, while (i) is equivalent to $s^*\phi- t^*\phi= 0$ at all $g \in G(\mathcal{F})$, 
(iii) is equivalent to the same condition at all $x\in M$ (by Corollary~\ref{unique}). 
Since  $s^*\phi-t^*\phi$ is a multiplicative 3-form which has zero differential,
the equivalence of (i) and (iii) follows from a degree three version of
Corollary \ref{unique}
(proven in the same way). 
\end{proof}

Let us point out that, although $\omega_{\phi}$ is not of Dirac type in general,
this form  is still relevant
for the construction of forms of Dirac type and of presymplectic groupoids.
In order to see that, note that 
the Dirac structure $L_{\mathcal{F}}$ is a
$\phi$-twisted Dirac structure, since $\phi$ vanishes along the leaves of the 
foliation.
By the properties of $\omega_{\phi}$ in the proof above, adding $\omega_{\phi}$
does not affect the isotropy bundle (see Definition \ref{robust}) of a
closed 2-form. Thus we get

\begin{corollary} Viewing $L_{\mathcal{F}}$ as a $\phi$-twisted Dirac structure,
the associated
$\phi$-twisted presymplectic groupoid is $G(\mathcal{F})\ltimes \nu^*$ with the
2-form 
$\omega_{\mathcal{F}}+ pr_1^{*}\omega_{\phi}$. 
\end{corollary}

%\end{example}

%%%%%%%%%%%%%%%%%%%%%%%%%%%%%%%%%%%%%%%
%%%%%%%%%%%%%%%%%%%%%%%%%%%%%%%%%%%%%%%
%%%%%%%%%%%%%%%%%%%%%%%%%%%%%%%%%%%%%%%
%%%%%%%%%%%%%%%%%%%%%%%%%%%%%%%%%%%%%%%
\subsection{Presymplectic groupoids of regular Dirac structures}
\label{regular-int}
%%%%%%%%%%%%%%%%%%%%%%%%%%%%%%%%%%%%%%%
%%%%%%%%%%%%%%%%%%%%%%%%%%%%%%%%%%%%%%%
%%%%%%%%%%%%%%%%%%%%%%%%%%%%%%%%%%%%%%%
%%%%%%%%%%%%%%%%%%%%%%%%%%%%%%%%%%%%%%%

We will call a Dirac structure {\bf regular} if its presymplectic leaves
have constant dimension.  To begin, we will restrict ourselves to the
untwisted case, with $\phi=0$. If $L$ is regular, then it determines

\begin{enumerate} 
\item a regular foliation $\mathcal{F}$ (whose leaves are the presymplectic
leaves of $L$); 
\item a closed foliated 2-form $\theta\in \Omega^2(\mathcal{F})$
(defined by the leafwise presymplectic forms of $L$).
\end{enumerate} 
Conversely, we can recover $L$ from $\mathcal{F}$ and $\theta$:
\[ 
L= \{ (X, \xi): X\in \mathcal{F}, \xi|_{\mathcal{F}}= i_{X}(\theta) \}.
\]
In this section, we discuss examples of regular Dirac structures
for which $G(L)$ admits a simplified description in terms of this data, 
$\mathcal{F}$ and $\theta$. (Note that the case $\theta= 0$ has been treated in Section
\ref{Dir-fol}.)

A simplified description of $G(L)$ depends on the {\bf classifying class} of
$L$, denoted by $c(L)$, that we now discuss. 
Using the transversal de Rham operator
$d_{\nu}$ defined in \eqref{eq:transv}, $c(L)$ is defined as
\[ 
c(L)= d_{\nu}(\theta)\in H^2(\mathcal{F}; \nu^*).
\]
While $\theta$ carries all the information of $L$ as a Dirac structure, $c(L)$
characterizes  $L$ as a Lie algebroid. As suggested by the
exact sequence
\[ 
0\rmap \nu^*\rmap L\rmap \mathcal{F}\rmap 0 ,
\]
the relationship between $L$ and $c(L)$ is the same as the one of
extensions and 2-cocycles, as briefly discussed in Section  \ref{subsec:over}. 
Let us make it more explicit. First of all, any 
closed $u\in \Omega^2(\mathcal{\F};\nu^*)$
defines an algebroid $\mathcal{F}\ltimes_{u} \nu^*$, with underlying vector bundle 
$\mathcal{F}\oplus \nu^*$,
projection on the first factor as anchor map, and  bracket
\[ 
[(X, v), (Y, w)]=([X, Y], \nabla_X(w)-\nabla_Y(v)+ u(X, Y)).
\]
When $u=0$, we simplify the  notation to $\mathcal{F}\ltimes \nu^*$. (Note
that this is the Lie algebroid underlying
the Dirac structure $L_{\mathcal{F}}$ of the Section \ref{Dir-fol}.)
The isomorphism class of the Lie algebroid $\mathcal{F}\ltimes_{u}
\nu^*$ 
depends only on the cohomology class of $u$: 
if $u'= u+ dv$ with $v\in \Omega^1(\mathcal{F}; \nu^*)$, then
$$
(X, \xi)\mapsto (X, \xi + v(X))$$ 
is an isomorphism between
$\mathcal{F}\ltimes_{u'} \nu^*$ and 
$\mathcal{F}\ltimes_{u} \nu^*$.

In order to see that $c(L)$ is the class corresponding to $L$, we 
choose a linear splitting $\sigma$ of the map $L\rmap \mathcal{F}$. 
On one hand, 
\begin{equation}
\label{u-sigma}
u_{\sigma}(X, Y)= [\sigma(X), \sigma(Y)]- \sigma([X, Y])
\end{equation}
is a representative
of $d_{\nu}(\theta)$ (this follows from the explicit description of $d_{\nu}$ given in 
Section \ref{two-forms on monodromy}); on the other hand, $\sigma$
induces a linear isomorphism $L\cong \mathcal{F}\oplus \nu^*$ which maps the
brackets on $L$ into the
brackets of $\mathcal{F}\ltimes_{u} \nu^*$.

\begin{example}(The case $c(L)= 0$) \rm
We now describe the presymplectic groupoid of $L$ when $c(L)=0$.
This case is closely related to our discussion 
in Section \ref{two-forms on monodromy}, which we now extend. 
Any multiplicative 2-form $\omega$ on the monodromy groupoid $G(\mathcal{F})$
defines a foliated form $c_{\omega}= \omega|_{\mathcal{F}}$ (where  
we view $\mathcal{F}\subset TM\subset TG(\mathcal{F})$), whose cohomology class is
precisely the $c(\omega)$
defined in Section \ref{two-forms on monodromy}.  In particular, we have an induced regular  
Dirac structure $L$ (namely the one defined by $\mathcal{F}$ and $c_{\omega}$). 
In this case, we say that $L$ {\bf comes from} $\omega$, and we write
$L= L(\omega)$.

\begin{corollary}
For a regular Dirac structure $L$ on $M$ the following are equivalent:
\begin{enumerate}
\item[(i)] $c(L)= 0$;
\item[(ii)] $L$ comes from a closed multiplicative 2-form on the monodromy
groupoid of $\mathcal{F}$;
\item[(iii)] the underlying algebroid of $L$ is isomorphic to
$\mathcal{F}\ltimes \nu^*$.
\end{enumerate}
In this case $L$ is integrable. Moreover, if one chooses $\omega$ as in (ii) (in which case $L= L(\omega)$),
then 
\[ 
(G(L),\omega_L) \cong (G(\mathcal{F})\ltimes \nu^* ,\omega_{\mathcal{F}} - pr_{1}^{*}\omega)
\]
(Here  $G(\mathcal{F})\ltimes \nu^*$ and $\omega_{\mathcal{F}}=
pr_{2}^{*}\omega_{can}$ 
are as in Section \ref{Dir-fol}).
\end{corollary}

\begin{proof} Using the map $\rho^{*}_{\omega}: \mathcal{F}\rmap T^*M$ induced from
$\omega$, we have
\[ 
L\cong \mathcal{F}\ltimes \nu^{*}, (v, \xi)\mapsto (v, \xi-\rho^{*}_{\omega}(v)),
\]
which is an isomorphism of Lie algebroids. Hence $G(L)\cong
G(\mathcal{F})\ltimes \nu^{*}$.

To find the 2-form $\omega_{L}$ on $G(\mathcal{F})\ltimes \nu^{*}$, we
look at its infinitesimal counterpart $\rho^{*}: \mathcal{F}\ltimes \nu^{*}\rmap
T^*M$. This is obtained
by transporting $pr_{2}: L\rmap T^*M$ (which defines $\omega_{L}$ on $G(L)$) by
the isomorphism above.
Hence $\rho^{*}(v, \xi)= \xi- \rho^{*}_{\omega}(v)$. Now, $\rho^{*}_{0}(v, \xi)=
\xi$
is precisely the infinitesimal counterpart of the multiplicative 2-form
$\omega_{\mathcal{F}}$,
while $\rho^{*}_{1}(v, \xi)= \rho^{*}_{\omega}(v)$ comes from the multiplicative
form $\omega$ on $G(\mathcal{F})$ 
and the projection on the first factor. Hence the form induced by $\rho^{*}$ is
$\omega_{\mathcal{F}} - pr_{1}^{*}\omega$.
\end{proof}
\end{example}

Another case in which we can make $G(L)$ more explicit is when 
$c(L)$ is integrable as a foliated cohomology class; as we will see,
this is similar to Van Est's approach to Lie's third theorem for Lie algebras (see \cite{Cra} and references therein).

\begin{example}(The case of integrable $c(L)$) \rm
 Recall that, if a Lie groupoid $G$ acts on a vector bundle $E$, we can define
differentiable cohomology groups
$H^{*}_{diff}(G; E)$, and the Van Est map maps these cohomology groups into
Lie algebroid cohomology with coefficients in $E$. We refer the reader to
\cite{WeXu, Cra} for a general discussion.
Here we only deal with the Van Est
map for $G(\mathcal{F})$, with coefficients in $\nu^*$, and in degree two:
\[ 
\Phi: H^{2}_{diff}(G(\mathcal{F}); \nu^*)\rmap H^{2}(\mathcal{F}; \nu^*).
\]
We now recall its definition. A differentiable 2-cocycle
on $G(\mathcal{F})$ with coefficients in $\nu^{*}$ is a smooth function $c$
which associates to any 
composable pair $(g, h)$ an element $c(g, h)\in \nu^{*}_{t(g)}$, 
which vanishes whenever $g$ or $h$ is a unit. We say that $c$ is closed if
\[ 
gc(h, k)- c(gh, k)+ c(g, hk)- c(g, h)= 0 ,
\]
for all triples $(g, h, k)$ of composable arrows in $G(\mathcal{F})$.
Two cocycles $c$ and $c'$ are said to be cohomologous if their difference is
of type $(g, h)\mapsto g d(h)- d(gh)+ d(g)$ for some section $d\in \Gamma(G;
t^{*}\nu^*)$.
This defines $H^{2}_{diff}(G(\mathcal{F}))$. 

Any closed $c$ defines a foliated form $\Phi(c)\in \Omega^2(\mathcal{F};
\nu^*)$:
roughly speaking, $\Phi(c)$ is obtained from $c$ by taking derivatives 
along leafwise vector fields (for the precise formulas, see \cite{Cra,WeXu}).
For our purpose, it will be useful to give a more abstract description of
$\Phi(c)$ using extensions. ;
Analogously to Lie algebroid extensions by  algebroid 2-cocycles,
differentiable 2-cocycles induce groupoid structures on $G(\mathcal{F})\times \nu^*$ 
with the multiplication extending the one
in $G(\mathcal{F})\ltimes \nu^*$:
\[ 
(g, v) (h, w)= (gh, h^{-1}v+ w+ (gh)^{-1}c(g, h)) .
\]
The resulting groupoid is denoted by $G(\mathcal{F})\ltimes_{c} \nu^*$. 
Still as in the infinitesimal case, the isomorphism class
of $G(\mathcal{F})\ltimes_{c} \nu^*$
depends  only on the cohomology class of $c$, and this groupoid fits into an
exact sequence of groupoids
\begin{equation}\label{eq:seqgr} 
1\rmap \nu^*\rmap G(\mathcal{F})\ltimes_{c} \nu^*\rmap G(\mathcal{F})\rmap 1.
\end{equation}
Passing to Lie algebroids, this induces an extension of
$\mathcal{F}$ by $\nu^*$,
hence a cohomology class in $H^{2}(\mathcal{F}; \nu^*)$. 
This defines $\Phi([c])$, and determines $\Phi$ at the cohomology level.
Note that this cohomology class has a canonical representative, and that
will define $\Phi(c)$, i.e, the map $\Phi$ at the chain level.
The exact sequence \eqref{eq:seqgr} has a canonical splitting, which induces a linear
splitting $\sigma$ at the algebroid level; the associated foliated form (\ref{u-sigma})
defines $\Phi(c)$.

\begin{corollary} If the characteristic class $c(L)$ comes from a 
differentiable cocycle $c$ (i.e. if $c(L)$ is integrable), 
then $L$ is integrable and $G(L)\cong G(\mathcal{F})\ltimes_{c} \nu^*$. 
\end{corollary}

At first sight, this corollary is just the definition of the
integrability of $c(L)$.
This is due to our definition of $\Phi$ in terms of extensions. 
However, \cite{Cra} gives us a precise description of
when $c(L)$ is integrable, related to the monodromy
groups of $L$. 
In this way  the result becomes meaningful. 

More precisely, by \cite{Cra} we know that $u\in H^2(\mathcal{\F}; \nu^*)$ is
integrable if and only if
all its leafwise periods vanish. This means that, for each leaf $S$ and any 
2-sphere $\gamma$
in $S$, $\int_{\gamma} u|_{S}= 0$. 
On the other hand, by the very definition of the monodromy groups $N_{x}(L)$
\cite{CrFe},
and by the description of $c(L)$ above in terms of a splitting $\sigma$, we have
\[ 
N_{x}(L)= \{ \int_{\gamma} c(L)|_{S}: \gamma\in \pi_2(S, x)\} ,
\]
where $S$ is the leaf through $x$. As in \cite{CrFe2}, this groups can be interpreted (or
defined) as the groups defined by the variations of the presymplectic areas. 
And, still completely analogous to the
Poisson case treated in \cite{CrFe2}, 
we state the conclusion without further details.

%This is not very surprising now,
%since we already know that the vanishing of the leafwise periods of $c(L)$
%implies integrability,
%and that the integrability is controlled by the monodromy groups $N_x(L)$.
%Note also that, since $c(L)= d_{\nu}([\theta])$,
%it follows from the functoriality of $d_{\nu}$ that the integrals above are the
%derivatives (variations) of the integrals of $\theta$, i.e. 
%\[ \int_{\gamma} c(L)|_{S}= \mathcal{A}_{\gamma}^{'}(\theta) .\]
%Hence the description of the monodromy groups given in Example
%\ref{integrability}.
%Again, we do not give more details since they are basically the same as in the
%Poisson case \cite{CrFe2}.
%In conclusion

\begin{corollary} The following are equivalent
\begin{enumerate} 
\item[(i)] $c(L)$ is integrable;
\item[(ii)] all the leafwise periods of $c(L)$ vanish;
\item[(iii)] all the monodromy groups $N_{x}(L)$ vanish.
\end{enumerate}
\end{corollary}
\end{example}

The discussion of general regular Dirac structures (i.e., when $c(L)$ is not necessarily
integrable)
can be treated, again, exactly as in the Poisson case \cite{CrFe2}. 

%So, let us briefly mention the highlights, completing the outline we started 
%in Example \ref{integrability}. 

%Following the Poisson case \cite{CrFe2}, we denote by $\mathcal{A}^{'}(L)$ the
%bundle of variations of the
%presymplectic areas. Hence, the fiber at $x$ is precisely the monodromy group
%$N_x(L)\subset \nu^*$.
%Define also {\it the structural bundle of $L$}, $\mathcal{S}(L)$, as the quotient
%$\nu^*/\mathcal{A}^{'}(L)$. 
%The global translation of the exact sequence of Lie algebroids above is that the
%topological groupoid
%$G(L)$ fits into an exact sequence of groupoids
%\[ 1\rmap \mathcal{S}(L)\rmap G(L)\rmap G(\mathcal{F})\rmap 1.\]
%We have seen above that $\mathcal{S}(L)= \nu^*$ if and only if $c(L)$ is
%integrable, and that
%implied that $G(L)= G(\mathcal{F})\times \nu^*$ as manifolds. In general, the
%sequence above gives 
%%us an idea of how $G(L)$ looks like. In particular, this implies that the
%connected component of the
%identity in the isotropy groups of $G(L)$ are isomorphic to $\nu^*{*}/N_x(L)$,
%hence the smoothness of $G(L)$
%forces the monodromy groups to be discrete (and that is part of the
%integrability obstructions).

%And, finally, as in the Poisson case, the following are equivalent:
%\begin{enumerate}
%\item[(i)] $L$ is integrable;
%\item[(ii)] $\mathcal{S}(L)$ is smooth;
%\item[(iii)] $\mathcal{A}^{'}(L)$ is an \'etale bundle of groups (i.e. smooth
%bundle, with discrete fibers).
%\end{enumerate}
%\end{remark}

So far, we have dealt with the case of $\phi=0$. However,
much of the discussion in this section extends to 
$\phi$-twisted regular Dirac structures.
For instance, one should use the $\phi$-twisted Courant bracket
in the construction of the foliated form 
(\ref{u-sigma}). This defines the classifying class 
of $L$, $c(L)\in H^2(\mathcal{F};\nu^*)$, and its integrals over leafwise 2-loops defines
the monodromy groups $N_x(L)$. 

We finish the section with remarks on the particular 
twisted case discussed in Section \ref{Dir-fol}.

\begin{example}
({\it Twisted regular Dirac structures})\rm

Let $\phi\in F_1\Omega^3(M)$, i.e. 
\[ 
i_Ui_Vi_W\phi= 0, \ \ \forall\ U, V, W\in \mathcal{F} .
\]

Let $L$ be a regular Dirac structure $L$ (with presymplectic foliation $\mathcal{F}$).
Then $L$ is automatically
a $\phi$-twisted regular presymplectic Dirac structure. To distinguish these two structures,
we write $(L, \phi)$ for $L$ viewed as a twisted Dirac structure. We also write
 $c(L, \phi)$ for its class, $G(L, \phi)$ for the associated groupoid, etc.

As in Section \ref{Dir-fol}, if $\phi\in F_2\Omega^3(M)$, then $G(L, \phi)= G(L)$;
if moreover $\phi\in F_3\Omega^3(M)$, then the corresponding 
multiplicative 2-forms coincide. 

Let us now consider
the class induced by $\phi$, $\bar{\phi}\in F_1\Omega^3(M)/F_2\Omega^3(M)\cong 
\Omega^2(\mathcal{F}; \nu^*)$, i.e.
\[ \bar{\phi}(V, W)= i_Vi_W(\phi) .\]
Note that, if $\bar{\phi}= d\psi$ for some $\psi\in \Omega^1(\mathcal{F}; \nu^*)$,
then Theorem \ref{rec1} applied to $G(\mathcal{F})$, with the twist by $\phi$, and
to $\rho^*: \mathcal{F}\rmap T^*M$ which is just $\psi$ viewed as a bundle map. 
In this way we get an  induced $\phi$-twisted multiplicative 2-form
\[ 
\omega_{\psi}\in \Omega^2(G(\mathcal{F})).
\]
Let us denote by $pr: G(L)\rmap G(\mathcal{F})$ the projection induced by the first 
projection $pr_1: L\rmap \mathcal{F}$. 

\begin{corollary} If $[\bar{\phi}]\in H^2(\mathcal{F}; \nu^*)$ 
vanishes, then $G(L, \phi)= G(L)$.

More precisely, any choice of $\psi\in \Omega^1(\mathcal{F}; \nu^*)$ such that $\bar{\phi}= d\psi$
induces an isomorphism $G(L, \phi)\cong G(L)$ which maps the multiplicative 2-form
on $G(L, \phi)$ to the form $\omega_L+ pr^*\omega_{\psi}$.
\end{corollary}

\begin{proof} First of all, by the description of the classifying classes in terms of
splittings $\sigma$ and of formula (\ref{u-sigma}), it immediately follows that 
the forms $u_{\sigma}$ representing $c(L)$, and $u_{\sigma, \phi}$ representing $c(L, \phi)$,
satisfy $u_{\sigma, \phi}= u_{\sigma}+ \bar{\phi}$. In particular,
\[ c(L, \phi)= c(L)+ [\bar{\phi}].\] 
Hence, by the classifying properties of 
$c(L)$, the first assertion follows. 

The second assertion follows from the general properties mentioned
above: $\sigma$ induces algebroid isomorphisms $L\cong \mathcal{F}\ltimes_{u_{\sigma}} \nu^*$,
$(L, \phi)\cong \mathcal{F}\ltimes_{u_{\sigma, \phi}} \nu^*$, while $\psi$ induces an isomorphism 
between the
algebroids associated to the two cocycles. Actually the resulting isomorphism does not depend on
$\sigma$, and is just $(X, \xi)\mapsto (X, \bar{\psi}(X)+ \xi)$. It is clear now that the infinitesimal
counterpart of the twisted multiplicative form is, on the untwisted $L$, just the sum of
$pr_{2}: L\rmap T^*M$ (defining $\omega_L$), and the composition of $pr_1:L\rmap \mathcal{F}$
with $\bar{\phi}$ (which is the infinitesimal counterpart of $pr^*\omega_{\psi}$).
This concludes the proof.
\end{proof}

\end{example}

% -----------------------------------------------------------------------
\bibliographystyle{amsplain}
\def\lllll{}

\end{document}